\documentclass[final]{siamltex}

\usepackage{amsfonts}
\newtheorem{remark}[theorem]{Remark}
\newtheorem{thm}{Theorem}

\def\N{{\mathbb N}}
\def\R{{\mathbb R}}
\def\A{{\bm A}}
\def\B{{\bm B}}
\def\E{{\bm E}}
\def\H{{\bm H}}
\def\I{{\bm I}}
\def\Q{{\bm Q}}
\def\S{{\bm S}}
\def\U{{\bm U}}
\def\W{{\bm W}}
\def\V{{\bm V}}
\def\ess{{\rm ess}}

\def\b{{\bm b}}
\def\f{{\bm f}}
\def\g{{\bm g}}
\def\n{{\bm n}}
\def\q{{\bm q}}
\def\u{{\bm u}}

\def\0{{\bf0}}

\def\alphaa{{\bm\alpha}}
\def\thetaa{{\bm\theta}}
\def\epsi{{\bm\varepsilon}}

\font\amb=cmmib10

\font\ambi=cmmib10 scaled 700

\newfam\mbfam \def\mb{\textfont1=\amb\fam\mbfam\amb\scriptfont1=\ambi}

\textfont\mbfam\amb \scriptfont\mbfam\ambi

\def\bm#1{\mathchoice
{\hbox{\mb\textfont1=\amb$#1$}}%
{\hbox{\mb\textfont1=\amb$#1$}}%
{\hbox{\mb$\scriptstyle\textfont1=\ambi#1$}}%
{\hbox{\mb$\scriptscriptstyle\textfont1=\ambi#1$}}}
\def\m@th{\mathsurround=0pt}
\def\lc#1{\hbox to .89\hsize{$\displaystyle{#1}\hfill$}\cr}
\def\rc#1{\hbox to .89\hsize{$\displaystyle\hfill{#1}$}\cr}
\def\eqal#1{\null\,\vcenter{\openup\jot\m@th
  \ialign{\strut\hfil$\displaystyle{##}$&&$\displaystyle{{}##}$\hfil
      \crcr#1\crcr}}\,}
\title{Global regular solutions to a~Kelvin-Voigt type 
thermoviscoelastic system\thanks{Partially supported 
by Polish Grant NN $201\,396\,937$}} 

\author{Irena Paw\l ow\thanks{Systems Research Institute, Polish Academy of Sciences,
Newelska 6, 01-447 Warsaw, Poland\hfil\break 
\indent\ \ \ and\hfil\break
\indent\ \ \ Institute of Mathematics and Cryptology, Cybernetics Faculty,
Military University\hfil\break 
\indent\ \ \ of Technology,S. Kaliskiego 2, 00-908 Warsaw, Poland
E-mail: pawlow@ibspan.waw.}\and Wojciech M. Zaj\c aczkowski\thanks{Institute of Mathematics, Polish Academy of Sciences, \'Sniadeckich 8, 00-950 Warsaw, Poland\hfil\break 
\indent\ \ \ and\hfil\break
\indent\ \ \ Institute of Mathematics and Cryptology, Cybernetics Faculty,
Military University\hfil\break 
\indent\ \ \ of Technology,S. Kaliskiego 2, 00-908 Warsaw, Poland E-mail: E-mail: wz@impan.gov.pl}}

\begin{document}
\maketitle

\begin{abstract}
A classical 3-D thermoviscoelastic system of Kelvin-Voigt 
type is considered. The existence and uniqueness of a global regular solution 
is proved without small data assumption. The existence proof is 
based on the successive approximation method. The crucial part constitute 
a~priori estimates on an arbitrary finite time interval, which are derived 
with the help of the theory of anisotropic Sobolev spaces with a~mixed norm.
\end{abstract}



\pagestyle{myheadings}
\thispagestyle{plain}
\markboth{I. Paw\l ow and W. M. Zajczkowski}{Global regular solutions to a~Kelvin-Voigt type 
thermoviscoelastic system}

\section{Introduction}
\subsection{Motivation and goal}

This article is concerned with the existence and uniqueness of global 
regular solutions to a~classical 3-D thermoviscoelastic system at small 
strains. The system describes materials which have the properties both of 
elasticity and viscosity. Such materials are usually referred to as 
Kelvin\--Voigt type.

As noted in the recent paper on this subject by Roub\'i\v cek [21] -- 
and according to our best knowledge as well -- the existence of global 
solutions to a~thermoviscoelastic system with constant both specific heat and 
heat conductivity is, in spite of great effort through many decades, still 
open in dimensions $n\ge2$. In dimension $n=1$ it was established in the 
pioneering papers by Slemrod [22], Dafermos [5] and Defermos-Hsiao [6].

The local in time existence and global uniqueness of a~weak solution to 3-D 
thermoviscoelastic system with constant specific heat and heat conductivity 
has been proved by Bonetti-Bonfanti [3]. Other known results on 
multidimensional thermoviscoelasticity deal with a~modified energy equation. 
Modifications involve either nonconstant specific heat or nonconstant heat 
conductivity. A thermoviscoelastic system with temperature-de\-pen\-dent 
specific heat has been addressed by Blanchard\--Guib\'e [2] where the 
existence of global, weak-renormalized solutions has been proved, and 
recently in [21] where the existence of a~very weak solution has been 
established.
We mention also that the framework of renormalized solutions has been 
applied in [26] for 3-D thermoviscoelastic system arising in structural 
phase transitions.

In a~more general setting allowing for large strains a~3-D thermoviscoelastic 
system has been studied under small data assumption by Shibata [23] and 
recently by Gawinecki-Zaj\c aczkowski [11].

For thermoviscoelastic problems with a~modified heat conductivity we refer to 
Eck-Jaru\v sek-Krbec [8] and the references therein.

In the present paper we consider a~thermoviscoelastic system with specific 
heat linearly increasing with temperature and with constant heat conductivity. 
Such setting is a~particular case of systems addressed in [2] and [21].

The novelty of the existence result presented in this paper concerns the 
regularity of a~3-D global solution corresponding to sufficiently  
smooth but arbitrary in size initial data. The proof of the existence 
theorem is based on the successive approximation method. The key regularity 
estimates are derived with the help of the parabolic theory in anisotropic 
Sobolev spaces $W_{p,p_0}^{2.1}(\Omega^T)$, $\Omega^T=\Omega\times(0,T)$, 
$p,p_0\in(1,\infty)$, with a~mixed norm with respect to space and time 
variables. Such framework has been previously applied by the authors [19] 
to the thermoviscoelastic system arising in shape memory alloys. It allowed 
to generalize the former results on this subject in [27].

As known, in deriving a~priori estimates for a~solution of a~system of 
balance laws it is common to begin with estimates arising from the 
conservation of a~total energy. Such estimates provide $L_\infty$-time 
regularity for the conserved quantities. To take advantage of such time 
regularity in deriving subsequent regularity estimates it is desirable to 
work in Sobolev spaces with a mixed norm, for example 
$W_{p,p_0}^{2,1}(\Omega^T)$, where the space exponent $p$ is determined by the 
energy structure and the time exponent $p_0$ may be arbitrarily large. This is 
the idea behind using the framework of Sobolev spaces with a~mixed norm to 
the thermoviscoelastic system under considerations.
The theory of IBVP's in Sobolev spaces with a~mixed norm is the subject of 
recent theoretical studies. We apply the general results due to Krylov [13] 
and Denk-Hieber-Pr\"uss [7].

\subsection{Thermoviscoelastic system}

The system under consideration has the following form:
$$
\u_{tt}-\nabla\cdot[\A_1\epsi_t+\A_2(\epsi-\theta\alphaa)]=\b,
\leqno(1.1)
$$
$$
c_v\theta\theta_t-k\Delta\theta=-\theta(\A_2\alphaa)\cdot\epsi_t+
(\A_1\epsi_t)\cdot\epsi_t+g\quad {\rm in}\ \ \Omega^T=\Omega\times(0,T),
\leqno(1.2)
$$
where
$$
\epsi\equiv\epsi(\u)={1\over2}(\nabla\u+(\nabla\u)^T),\quad 
\epsi_t\equiv\epsi(\u_t)={1\over2}(\nabla\u_t+(\nabla\u_t)^T).
$$
Here $\Omega\subset\R^3$ is a~bounded domain occupied by a~body in a~fixed 
reference configuration, and $(0,T)$ is the time interval. The system is 
completed by appropriate boundary and initial conditions. Here we assume
$$
\u=\0,\ \ \n\cdot\nabla\theta=0\quad{\rm on}\ \ S^T=S\times(0,T),
\leqno(1.3)
$$
$$
\u|_{t=0}=\u_0,\ \ \u_t|_{t=0}=\u_1,\ \ \theta|_{t=0}=\theta_0\quad {\rm in}
\ \ \Omega,
\leqno(1.4)
$$
where $S$ is the boundary of $\Omega$ and $\n$ is the unit outward normal 
to $S$.
\goodbreak

The field $\u:\ \Omega^T\to\R^3$ is the displacement, 
$\theta:\ \Omega^T\to\R_+=(0,\infty)$ is the absolute temperature, the second 
order tensors $\epsi=(\varepsilon_{ij})_{i,j=1,2,3}$ and 
$\epsi_t=((\varepsilon_t)_{ij})_{i,j=1,2,3}$ denote respectively the 
linearized strain and the strain rate.

Equation (1.1) is the linear momentum balance with the stress tensor given 
by a~linear thermoviscoelastic law of the Kelvin-Voigt type (cf. [8], 
Chap. 5.4)
$$
\S=\A_1\epsi_t+\A_2(\epsi-\theta\alphaa).
$$
The fourth order tensors $\A_1=((A_1)_{ijkl})_{i,j,k,l=1,2,3}$ and\\ 
$\A_2=((A_2)_{ijkl})_{i,j,k,l=1,2,3}$ are respectively the linear viscosity 
and the elasticity tensors, defined by
$$
\epsi\mapsto\A_p\epsi=\lambda_ptr\epsi\I+2\mu_p\epsi,\quad p=1,2,
\leqno(1.5)
$$
where $\lambda_1$, $\mu_1$ are the viscosity constants, and $\lambda_2$, 
$\mu_2$ are the Lam\'e constants, both $\lambda_1$, $\mu_1$ and $\lambda_2$, 
$\mu_2$ with the values within the elasticity range
$$
\mu_p>0,\quad 3\lambda_p+2\mu_p>0,\quad p=1,2;
\leqno(1.6)
$$
$\I=(\delta_{ij})_{i,j=1,2,3}$ is the identity tensor.

The second order symmetric tensor $\alphaa=(\alpha_{ij})_{i,j=1,2}$ with 
constant $\alpha_{ij}$, represents the thermal expansion. The vector field 
$\b:\ \Omega^T\to\R^3$ is the external body force.

Above and hereafer the summation convention over the repeated indices is 
used, vectors and tensors are denoted by bold letters, and the dot denotes 
the inner product of tensors, e.g.
$$
(\A\epsi)\cdot\epsi=A_{ijkl}\varepsilon_{kl}\varepsilon_{ij}.
$$
Moreover,
$$
\A\epsi=(A_{ijkl}\varepsilon_{kl})_{i,j=1,2,3}\quad {\rm and}\quad 
\nabla\cdot(\A\epsi)=\bigg({\partial\over\partial x_j}
(A_{ijkl}\varepsilon_{kl})\bigg)_{i=1,2,3}.
$$

Equation (1.2) is the energy balance in which the linear Fourier law for 
the heat flux, $\q=-k\nabla\theta$ with constant heat conductivity $k>0$, 
and temperature-dependent specific heat, $c_v\theta$ with $c_v>0$, have been 
adopted. The first two terms on the right-hand side of (1.2) represent heat 
sources created by the deformation of the material and by the viscosity. 
The field $g:\ \Omega^T\to\R$ is the external heat source.

The boundary conditions in (1.3) mean that the body is fixed at the boundary 
$S$ and thermally isolated. The initial conditions (1.4) prescribe 
displacement, velocity and temperature at $t=0$.

The system (1.1)--(1.2) can be derived by various arguments of 
thermodynamics, see e.g. [9, 3]. In Section 2 we summarize its 
thermodynamic basis.

\subsection{Linear elasticity and viscosity operators. Assumptions}

For further analysis we formulate problem (1.1)--(1.4) in terms of the 
linear viscosity and elasticity operators, $\Q_1$ and $\Q_2$, defined by
$$
\u\mapsto\Q_p\u=\nabla\cdot(\A_p\epsi(\u))=\mu_p\Delta\u+
(\lambda_p+\mu_p)\nabla(\nabla\cdot\u),\quad p=1,2,
\leqno(1.7)
$$
with domains $D(\Q_p)=\H^2(\Omega)\cap\H_0^1(\Omega)$.

Then system (1.1), (1.2) takes the form
$$\eqal{
&\u_{tt}-\Q_1\u_t=\Q_2\u-\nabla\cdot(\theta\A_2\alphaa)+\b,\cr
&c_v\theta\theta_t-k\Delta\theta=-\theta(\A_2\alphaa)\cdot\epsi_t+
(\A_1\epsi_t)\cdot\epsi_t+g\quad {\rm in}\ \ \Omega^T,\cr}
\leqno(1.8)
$$
with boundary and initial conditions (1.3), (1.4).

Throughout we shall assume that
\vskip6pt

\noindent
{\bf(A1)} $\Omega\subset\R^3$ is a~bounded domain with the boundary $S$ of 
class at least $C^2$; $T>0$ is an arbitrary finite number;
\vskip6pt

\noindent
{\bf(A2)} $\alphaa=(\alpha_{ij})_{i,j=1,2,3}$ is a~second order symmetric 
tensor with constant $\alpha_{ij}$;
\vskip6pt

\noindent
{\bf(A3)} The fourth order tensors $\A_1$ and $\A_2$ are defined by (1.5) 
with the coefficients $\mu_p$, $\lambda_p$, $p=1,2$, satisfying (1.6).
\vskip6pt

We list the implications of assumption {\bf(A3)} which are used in further 
analysis. The condition (1.6) ensures the symmetry of tensors $\A_p$:
$$
(A_p)_{ijkl}=(A_p)_{jikl}=(A_p)_{klij},\quad p=1,2,
\leqno(1.9)
$$
and their coercivity and boundedness
$$
a_{p*}|\epsi|^2\le(\A_p\epsi)\cdot\epsi\le a_p^*|\epsi|^2,\quad p=1,2,
\leqno(1.10)
$$
where
$$
a_{p*}=\min\{3\lambda_p+2\mu_p,2\mu_p\},\quad 
a_p^*=\max\{3\lambda_p+2\mu_p,2\mu_p\}.
$$

\noindent
Moreover, (1.6) ensures the following properties of operators $\Q_p$, $p=1,2$:
\begin{itemize}

\item[---] $\Q_p$ are strongly elliptic (property holding true under weaker 
assumption $\mu_p>0$, $\lambda_p+2\mu_p>0$, (see [20], Sect. 7)) and satisfy 
the estimate [17], Lemma 3.2:
$$
{\bm c}_p\|\u\|_{\H^2(\Omega)}\le\|\Q_p\u\|_{{\bm L}_2(\Omega)}\quad 
{\rm for}\ \ \u\in D(\Q_p),\ \ p=1,2,
\leqno(1.11)
$$
with positive constants ${\bm c}_p$ depending on $\Omega$. Since clearly,
$$
\|\Q_p\u\|_{{\bm L}_2(\Omega)}\le\bar c_p\|\u\|_{\H^2(\Omega)},
$$
it follows that the norms $\|\Q_p\u\|_{{\bm L}_2(\Omega)}$ and 
$\|\u\|_{\H^2(\Omega)}$ are equivalent on $D(\Q_p)$;

\item[---] the operators $\Q_p$ are self-adjoint on $D(\Q_p)$:
$$\eqal{
(\Q_p\u,{\bm v})_{{\bm L}_2(\Omega)}&=-\mu_p
(\nabla\u,\nabla{\bm v})_{{\bm L}_2(\Omega)}
-(\lambda_p+\mu_p)(\nabla\cdot\u,\nabla\cdot
{\bm v})_{{\bm L}_2(\Omega)}\cr
&=(\u,\Q_p{\bm v})_{{\bm L}_2(\Omega)}\quad{\rm for}\ \ 
\u,{\bm v}\in D(\Q_p);\cr}
\leqno(1.12)
$$

\item[---] the operators $-\Q_p$ are positive on $D(\Q_p)$:
$$\eqal{
(-\Q_p\u,\u)&=\mu_p\|\nabla\u\|_{{\bm L}_2(\Omega)}^2+(\lambda_p+\mu_p)
\|\nabla\cdot\u\|_{L_2(\Omega)}^2\ge0\cr
&\qquad\qquad{\rm for}\ \ \u\in D(\Q_p).\cr}
\leqno(1.13)
$$
Hence, there exist fractional powers $\Q_p^{1/2}$ with the domains \\
$D(\Q_p^{1/2})=\H_0^1(\Omega)$, satisfying
$$\eqal{
(\Q_p^{1/2}\u,\Q_p^{1/2}{\bm v})_{{\bm L}_2(\Omega)}
&=(-\Q_p\u,{\bm v})_{{\bm L}_2(\Omega)}=(\u,-\Q_p{\bm v})_{{\bm L}_2(\Omega)}
\cr
&\qquad\qquad{\rm for}\ \ \u,{\bm v}\in D(\Q_p).\cr}
\leqno(1.14)
$$
\end{itemize}

Let us also notice that by (1.10) and the Korn inequality
$$
d^{1/2}\|\u\|_{\H^1(\Omega)}\le\|\epsi(\u)\|_{{\bm L}_2(\Omega)}\quad 
{\rm for}\quad\u\in\H_0^1(\Omega),\ \ d>0,
\leqno(1.15)
$$
it follows that
$$\eqal{
&\|\Q_p^{1/2}\u\|_{{\bm L}_2(\Omega)}^2=\mu_p\|\nabla\u\|_{{\bm L}_2(\Omega)}^2+
(\lambda_p+\mu_p)\|\nabla\cdot\u\|_{L_2(\Omega)}^2\cr
&=(\A_p\epsi(\u),\epsi(\u))_{{\bm L}_2(\Omega)}\ge a_{p*}
\|\epsi(\u)\|_{{\bm L}_2(\Omega)}^2\ge a_{p*}d\|\u\|_{\H^1(\Omega)}^2.\cr}
\leqno(1.16)
$$
Thus, the norms $\|\Q_p^{1/2}\u\|_{{\bm L}_2(\Omega)}$ and $\|\u\|_{\H^1(\Omega)}$ 
are equivalent on $D(\Q_p^{1/2})$.

\subsection{Main result}\hfill\break

\begin{thm} (Existence) 
Let the assumptions (A1)--(A3) hold, $S\in C^2$, $T>0$ finite and
$$\eqal{
&\u_0\in\W_{12}^2(\Omega),\quad &\u_1\in\B_{12,12}^{11/6}(\Omega),\quad
&\theta_0\in B_{6,6}^{5/3}(\Omega),\cr
&g\in L_{\infty,12}(\Omega^T),\quad &\b\in{\bm L}_{12}(\Omega^T),\quad
&g\ge0,\ \ \theta_0\ge\underline{\theta}>0,\cr}
$$
where $\underline{\theta}$ is a constant. 
Then there exists a solution to problem (1.1)--(1.4) such that 
$\u\in C([0,T];\W_{12}^2(\Omega))$, $\u_t\in\W_{12}^{2,1}(\Omega^T)$, 
$\theta\in W_6^{2,1}(\Omega^T)$ and 
$\theta(t)\ge\underline{\theta}\exp(-c_0T)\equiv\theta_*>0$, where the 
positive constant $c_0$ depends on $a_{1*}$, $a_2^*$, $|\alphaa|$, $c_v$.\\
Moreover, the following estimates are satisfied
$$\eqal{
&\|\u\|_{C([0,T];\W_{12}^2(\Omega))}\le c\|\u_t\|_{\W_{12}^{2,1}(\Omega^T)},\cr
&\|\u_t\|_{\W_{12}^{2,1}(\Omega^T)}+\|\theta\|_{W_6^{2,1}(\Omega^T)}\le
\varphi(T,\|\u_0\|_{\W_{12}^2(\Omega)}+\|\u_1\|_{\B_{12,12}^{11/6}(\Omega)}\cr
&\quad+\|\theta_0\|_{B_{6,6}^{5/3}(\Omega)}+\|\b\|_{{\bm L}_{12}(\Omega)}+
\|g\|_{L_{\infty,12}(\Omega^T)}),\cr}
$$
where $\varphi$ is an increasing positive function of its arguments.
\end{thm}

\begin{thm} (Uniqueness) 
Let us assume that tensors $\A_p$, $p=1,2$, satisfy (1.10). 
Then any solution $(\u,\theta)$ to problem (1.1)--(1.4) satisfying
$$\eqal{
&\epsi_t\in L_2(0,T;{\bm L}_3(\Omega)),\cr
&\theta\in L_2(0,T;L_\infty(\Omega)),\quad \theta_t\in L_2(0,T;L_3(\Omega)),\cr
&0<\theta_*<\theta,\cr}
\leqno(1.17)
$$
is uniquely defined.
\end{thm}
\vskip12pt

\begin{corollary} 
The regular solution in Theorem A is uniquely defined.
\end{corollary}
\subsection{Relation to other results}

We comment on the connections of our result to the two other global 
existence results in three space dimensions. Firstly, we mention the result 
by Roub\' i\v cek [21] who proved the existence of a very weak solution to the 
thermoviscoelasticity system (1.1)--(1.2) involving monotone viscosity of 
a $p$-Laplacian type, $(\A_1\epsi_t)\cdot\epsi_t\sim|\epsi_t|^p$, and the 
specific heat having $(\omega-1)$-polynomial growth, 
$c_*(\theta)\sim c_v\theta^{\omega-1}$. This result, based on the Galerkin 
method, was obtained for $L^1$-data under the conditions $p\ge2$, $\omega\ge1$ 
and $p>1+{3\over2\omega}$ (in 3-D).
In the case of linear viscosity, $p=2$, the latter condition implies that 
$\omega>3/2$, that is the growth of the specific heat is greater than $1/2$.

Our result concerns the case $p=2$ and $\omega=2$. We have to restrict 
ourselves to the linear viscosity, $p=2$, because the proof relies on the 
results by Krylov [13] and Solonnikov [24] on the solvability of the linear 
problem
$$
\u_{tt}-\nabla\cdot\A_1\epsi(\u_t)=\f
$$
with the boundary and initial conditions (1.3), (1.4) (see Lemma 3.4).

\noindent
Concerning the specific heat growth exponent, $\omega-1$, it seems that 
after some additional technical effort it would be possible to admit 
$\omega<2$. However, in the case of a constant specific heat, i.e., 
$\omega=1$, we have been faced with a serious mathematical obstacle.

\noindent
As already mentioned in Subsection 1.1, the local existence result in such 
a case was obtained by Bonetti and Bonfanti [3].

Secondly, we recall the multidimensional result by Blanchard and Guib\'e [2] 
who addressed problem (1.1)--(1.2) equally in the prototype case $p=2$ and 
$\omega=2$, and in the more general setting involving linear viscosity, 
$p=2$, specific heat with $(\omega-1)$-polynomial growth and a nonlinear 
thermoelastic coupling; more precisely, the term $\nabla\theta$ in (1.1) was 
replaced by $\nabla f(\theta)$ with $f$ hawing an $\alpha$-polynomial growth. 
The existence of solutions in the weak-renormalized sense was proved there by 
the Schauder fixed point theorem. It is worth to remark that in the case of 
the linear thermoelastic coupling, $\alpha=1$, the result in [2] requires the 
specific heat to have growth of the order greater than $1/2$, as in the result 
by Roub\'i\v cek [21].

\subsection{Outline}

In Section 2 we present a~thermodynamic basis of system (1.1), (1.2). 
Section 3 recalls basic results on the Sobolev spaces with a~mixed norm and on 
the solvability of boundary-value problems for linear parabolic equations in 
such spaces. In Section 4 we derive a~priori estimates for problem 
(1.1)--(1.4). The procedure consists in a~recursive improvement of the basic 
energy estimates. The main tool in this procedure are the results on the 
solvability of linear parabolic problems in Sobolev spaces with a~mixed norm. 
Section 5 presents the proof of Theorem A, which is based on the successive 
approximation method.
The proof of the uniqueness, stated in Theorem B, is given in Section 6.

Since a priori estimates in Section 4 are crucial for the proof of the global 
existence we advertise here the main steps of the procedure of deriving such 
estimates. First we prove the energy type estimate (see Lemma 4.2)
$$
\|\u_t\|_{{\bm L}_{2,\infty}(\Omega^T)}+
\|\epsi\|_{{\bm L}_{2,\infty}(\Omega^T})+\|\theta\|_{L_{2,\infty}(\Omega^T)}
\le{\rm data}.
\leqno(1.18)
$$
In Lemma 4.6 we show the estimates
$$\eqal{
&\|\u_t\|_{{\bm L}_{2,\infty}(\Omega^T)}+\|\u\|_{L_\infty(0,T;\H^1(\Omega))}+
\|\u_t\|_{L_2(0,T;\H^1(\Omega))}
&\le c\|\theta\|_{L_2(\Omega^T)}+{\rm data},\cr}
$$
and
$$\eqal{
&\|\u_t\|_{L_\infty(0,T;\H^1(\Omega))}+\|\u\|_{L_\infty(0,T;\H^2(\Omega))}+
\|\u_t\|_{L_2(0,T;\H^2(\Omega))}
&\le c\|\nabla\theta\|_{{\bm L}_2(\Omega^T)}+{\rm data}.\cr}
$$
The norms of $\theta$ will be later removed by some interpolation inequalities 
based on estimate (1.18).

\noindent
In Lemma 4.7 we obtain the estimate
$$
\|\theta\|_{L_\infty(0,T;L_3(\Omega))}+\|\theta\|_{L_2(0,T;H^1(\Omega))}+
\|\epsi_t\|_{\V_2(\Omega^T)}\le{\rm data},
\leqno(1.19)
$$
and next in Lemma 4.8,
$$
\|\theta_t\|_{L_2(\Omega^T)}+\|\nabla\theta\|_{L_\infty(0,T;{\bm L}_2(\Omega))}
\le{\rm data}.
$$
To deduce the boundedness of $\theta$ we first prove in Lemma 4.10 the estimate
$$
\|\theta\|_{L_{r,\infty}(\Omega^T)}\le{\rm data},\quad r<\infty,
\leqno(1.20)
$$
and in Lemma 4.17,
$$
\|\theta\|_{L_\infty(\Omega^T)}\le{\rm data}.
\leqno(1.21)
$$
To get (1.21) we make use of the important inequality (see Lemma 4.9)
$$
\|\epsi_t\|_{{\bm L}_{p,\sigma}(\Omega^T)}\le c(t)
(\|\theta\|_{L_{p,\sigma}(\Omega^T)}+{\rm data}),\quad p,\sigma\in(1,\infty),
$$
and the fact that the coefficient near $\theta_t$ in (1.2) is proportional 
to $\theta$.

\noindent
To establish the continuity of $\theta$ (proved in Lemma 4.10) we need (1.21), 
estimates $\|\epsi_t\|_{L_2(0,T;{\bm L}_\infty(\Omega))}\le{\rm data}$ 
(see Corollary 4.16), and
$$
\|\theta\|_{W_2^{2,1}(\Omega^T)}\le{\rm data}\quad 
{\rm (see\ Corollary\ 4.18)}.
$$
Having the previous estimates for $\epsi_t$ and the continuity of $\theta$ 
we finally prove in Lemma 4.23 that
$$
\theta\in W_6^{2,1}(\Omega^T),\quad \u_t\in\W_{12}^{2,1}(\Omega^T).
$$

\section{Thermodynamic basis}

System (1.1), (1.2) represents balance laws for the linear momentum and energy 
in a~referential description, with the referential mass density assumed 
constant, normalized to unity, $\rho_0=1$:
$$\eqal{
\u_{tt}-\nabla\cdot\S&=\b,\cr
e_t+\nabla\cdot\q-\S\cdot\epsi_t&=g,\cr}
\leqno(2.1)
$$
where $\S$ is the stress tensor, $\q$ -- the referential heat flux, and 
$e$ -- the specific internal energy.

The system is governed by two thermodynamic potentials: the free energy 
$f=\hat f(\epsi,\theta)$, which by a~thermodynamic requirement is strictly 
concave with respect to $\theta$, and the dissipation potential (called 
pseudopotential of dissipation in [10], [3]) 
${\cal D}=\hat{\cal D}(\epsi_t,\nabla\theta;\epsi,\theta)$, which by 
a~thermodynamic requirement is nonnegative, convex in $(\epsi_t,\nabla\theta)$ 
and such that $\hat{\cal D}(\0,\0;\epsi,\theta)=0$.

In the case of (1.1), (1.2) the free energy is specified by
$$
f(\epsi,\theta)=f_*(\theta)+W(\epsi,\theta),
\leqno(2.2)
$$
where
$$
f_*(\theta)=-{1\over2}c_v\theta^2,\quad c_v=\const>0,
\leqno(2.3)
$$
is the caloric energy, and
$$\eqal{
W(\epsi,\theta)&={1\over2}(\epsi-\theta\alphaa)\cdot\A_2(\epsi-\theta\alphaa)-
{\theta^2\over2}\alphaa\cdot(\A_2\alphaa)\cr
&={1\over2}\epsi\cdot(\A_2\epsi)-\theta\epsi\cdot(\A_2\alphaa)\cr}
\leqno(2.4)
$$
is the elastic energy; we recall that $\A_2$ stands for the fourth order 
elasticity tensor and $\alphaa$ for the second order thermal expansion tensor.

The caloric energy (2.3) is associated with temperature-dependent ca\-lo\-ric 
specific heat
$$
c_*(\theta)=-\theta f''_*(\theta)=c_v\theta,
\leqno(2.5)
$$
which gives rise to the term $c_v\theta\theta_t$ in energy equation (1.2).

We remark that in the case $f_*$ is given by the standard formula
$$
f_*(\theta)=-c_v\theta\log{\theta\over\theta_1}+c_v\theta+\tilde c,
\leqno(2.6)
$$
where $c_v$, $\theta_1$, $\tilde c$ are positive constants, the caloric 
specific heat is constant
$$
c_*=-\theta f''_*(\theta)=c_v.
\leqno(2.7)
$$
This gives rise to the usual parabolic term $c_v\theta_t$ in (1.2) in place 
of $c_v\theta\theta_t$. As mentioned in Section 1, in the case of a constant 
caloric heat there are serious mathematical obstacles in the proof of the 
global existence.

A thermoviscoelastic system with the specific heat $c_*(\theta)$ given by 
(2.5) has been considered in [9], [2] and [21], where also more general forms 
of $c_*(\theta)$ have been analysed.
Moreover, we mention that a~fourth order thermoviscoelastic systems with 
tem\-pe\-ra\-tu\-re\--dependent specific heat, arising in shape memory 
materials, have been studied in [27] and [19].

The dissipation potential corresponding to system (1.1), (1.2) is given by
$$\eqal{
{\cal D}&={1\over2\theta}\epsi_t\cdot(\A_1\epsi_1)+{k\over2}\theta^2
\bigg|\nabla{1\over\theta}\bigg|^2
&={1\over2\theta}\epsi_t\cdot(\A_1\epsi_t)+{k\over2}|\nabla\log\theta|^2,\cr}
\leqno(2.8)
$$
where $\A_1$ is the viscosity tensor and $k>0$ the constant heat conductivity.

In accord with the basic thermodynamic relations the internal energy $e$ and 
the entropy $\eta$ are related to the free energy $f$ by the equations
$$
e=f+\theta\eta,\quad \eta=-f_{,\theta}.
\leqno(2.9)
$$
For the free energy $f$ defined by (2.2)--(2.4) this gives
$$
e={1\over2}c_v\theta^2+{1\over2}\epsi\cdot(\A_2\epsi),\quad
\eta=c_v\theta+(\A_2\alphaa)\cdot\epsi.
\leqno(2.10)
$$

As a~consequence of the second law of thermodynamics expressed by the 
Clausius-Duhem inequality, the stress tensor $\S$ and the heat flux $\q$ 
satisfy the following relations:
$$
\S={\partial f\over\partial\epsi}+
\theta{\partial{\cal D}\over\partial\epsi_t},\quad
\q={\partial{\cal D}\over\partial\nabla{1\over\theta}}.
\leqno(2.11)
$$

For $f$ given by (2.2)--(2.4) and $\cal D$ by (2.8) the formulas (2.11) yield 
the standard forms of the stress tensor and the heat flux
$$
\S=\A_2(\epsi-\theta\alphaa)+\A_1\epsi_t,\quad
\q=k\theta^2\nabla{1\over\theta}=-k\nabla\theta.
\leqno(2.12)
$$
Thus, $\S$ consists of two terms: the nondissipative equilibrium term 
determined by $f$, and the dissipative one determined by $\cal D$. The 
dissipative heat flux $\q$ is entirely determined by $\cal D$.

Inserting the relations $(2.10)_1$ and (2.12) into balance laws (2.1) one 
arrives at the system (1.1)--(1.2).

For further purposes (see Lemma 4.2) it is of interest to notice that on 
account of the identity
$$
e_t=(f+\theta\eta)_t=f_t+\theta_t\eta+\theta\eta_t=\theta\eta_t+
{\partial f\over\partial\epsi}\cdot\epsi_t,
$$
along with the relation $(2.11)_1$, the energy balance $(2.1)_2$ admits the 
form
$$
\theta\eta_t+\nabla\cdot\q=\theta{\partial{\cal D}\over\partial\epsi_t}
\cdot\epsi_t+g.
\leqno(2.13)
$$
For $\cal D$ given by (2.8) this leads to the following equivalent form 
of equation (1.2):
$$
\theta\eta_t-k\Delta\theta=(\A_1\epsi_t)\cdot\epsi_t+g.
\leqno(2.14)
$$

Let us also notice that assuming $\theta>0$ and using $(2.11)_2$, the equation 
(2.13) may be expressed as
$$
\eta_t+\nabla\cdot{q\over\theta}=\sigma+{g\over\theta},
\leqno(2.15)
$$
where
$$
\sigma={\partial{\cal D}\over\partial\nabla{1\over\theta}}\cdot\nabla
{1\over\theta}+{\partial{\cal D}\over\partial\epsi_t}\cdot\epsi_t
=k\theta^2\bigg|\nabla{1\over\theta}\bigg|^2+{1\over\theta}(\A_1\epsi_t)\cdot
\epsi_t\ge0
$$
is the specific entropy production. From (2.15) it follows that system (1.1), 
(1.2) complies with the Clausius-Duhem inequality
$$
\eta_t+\nabla\cdot{\q\over\theta}\ge{g\over\theta}.
\leqno(2.16)
$$

\section{Notation and auxiliary results}

\subsection{Notation}

Let $\Omega\subset\R^n$ be an open bounded subset of $\R^n$, $n\ge1$, with 
a~smooth boundary $S$, and $\Omega^T=\Omega\times(0,T)$, $S^T=S\times(0,T)$, 
$T>0$ finite.

\noindent
We introduce the following spaces: $W_p^k(\Omega)$, 
$k\in\N\cup\{0\}$, $p\in[1,\infty)$ -- the Sobolev space on $\Omega$ 
with the finite norm
$$
\|u\|_{W_p^k(\Omega)}=\bigg(\sum_{|\alpha|\le k}\intop_\Omega
|D_x^\alpha u|^pdx\bigg)^{1/p},
$$
where $\alpha=(\alpha_1,\dots,\alpha_n)$ is a~multiindex,
$\alpha_i\in\N\cup\{0\}$, $|\alpha|=\alpha_1+\alpha_2+\dots+\alpha_n$, 
$D_x^\alpha=\partial_{x_1}^{\alpha_1}\dots\partial_{x_n}^{\alpha_n}$; 
$H^k(\Omega)=W_2^k(\Omega)$; 
$L_{p,p_0}(\Omega^T)=L_{p_0}(0,T;L_p(\Omega))$, $p,p_0\in[1,\infty)$ -- 
the space of functions $u:\ (0,T)\to L_p(\Omega)$ with the finite norm
$$
\|u\|_{L_{p,p_0}(\Omega^T)}=\bigg(\intop_0^T\|u(t)\|_{L_p(\Omega)}^{p_0}dt
\bigg)^{1/p_0};
$$
$V_2(\Omega^T)=L_\infty(0,T;L_2(\Omega))\cap L_2(0,T;H^1(\Omega))$ -- the 
space of functions $u:(0,T)\to H^1(\Omega)$ with the finite norm
$$
\|u\|_{V_2(\Omega^T)}=\ess\sup_{t\in[0,T]}\|u(t)\|_{L_2(\Omega)}+
\|\nabla u\|_{{\bm L}_2(\Omega^T)};
$$
$V_2^{1,0}(\Omega^T)=V_2(\Omega^T)\cap C([0,T];L_2(\Omega))$ -- the space 
with the finite norm
$$
\|u\|_{V_2^{1,0}(\Omega^T)}=\max_{t\in[0,T]}\|u(t)\|_{L_2(\Omega)}+
\|\nabla u\|_{L_2(\Omega^T)};
$$
$W_{p,p_0}^{k,k/2}(\Omega^T)$, $k,k/2\in\N\cup\{0\}$, $p,p_0\in[1,\infty)$ -- 
the Sobolev space with a~mixed norm, which is a completion of 
$C^\infty(\Omega^T)$-functions under the finite norm
$$
\|u\|_{W_{p,p_0}^{k,k/2}(\Omega^T)}=\bigg(\intop_0^T\bigg(
\sum_{|\alpha|+2a\le k}\intop_\Omega|D_x^\alpha\partial_t^au|^pdx\bigg)^{p_0/p}
dt\bigg)^{1/p_0};
$$
$W_{p,p_0}^{s,s/2}(\Omega^T)$, $s\in\R_+$, $p,p_0\in[1,\infty)$ -- the 
Sobolev-Slobodecki space with the finite norm
$$\eqal{
&\|u\|_{W_{p,p_0}^{s,s/2}(\Omega^T)}=\sum_{|\alpha|+2a\le[s]}
\|D_x^\alpha\partial_t^au\|_{L_{p,p_0}(\Omega^T)}\cr
&+\bigg[\intop_0^T\bigg(\intop_\Omega\intop_\Omega\sum_{|\alpha|+2a=[s]}
{|D_x^\alpha\partial_t^au(x,t)-D_{x'}^\alpha\partial_t^au(x',t)|^p\over
|x-x'|^{n+p(s-[s])}}dxdx'\bigg)^{p_0/p}dt\bigg]^{1/p}\cr
&+\bigg[\intop_0^T\intop_0^T\bigg(\intop_\Omega\sum_{|\alpha|+2a=[s]}
{|D_x^\alpha\partial_t^au(x,t)-D_x^\alpha\partial_{t'}^au(x,t')|^p\over
|t-t'|^{1+p\left({s\over2}-\left[{s\over2}\right]\right)}}dx\bigg)^{p_0/p}
dtdt'\bigg]^{1/p_0},\cr}
$$
where $a\in\N\cup\{0\}$ and $[s]$ is the integer part of $s$. For $s$ odd the 
last term in the above norm vanishes whereas for $s$ even the two last terms 
vanish.

\noindent
$B_{p,p_0}^l(\Omega)$, $l\in\R_+$, $p,p_0\in[1,\infty)$ -- the Besov space 
with the finite norm
$$
\|u\|_{B_{p,p_0}^l(\Omega)}=\|u\|_{L_p(\Omega)}+\bigg(\sum_{i=1}^n
\intop_0^\infty{\|\Delta_i^m(h,\Omega)\partial_{x_i}^ku\|_{L_p(\Omega)}^{p_0}
\over h^{1+(l-k)p_0}}dh\bigg)^{1/p_0},
$$
where:
$$\eqal{
&k\in\N\cup\{0\},\quad m\in\N,\quad m>l-k>0,\cr
&\Delta_i^j(h,\Omega)u,\quad j\in\N,\quad h\in\R_+,\quad {\rm is\ the\ 
finite\ difference\ of\ the\ order}\cr
&j\ {\rm of\ the\ function}\ u(x)\ {\rm with\ respect\ to}\ x_i,\ {\rm with}\cr
&\Delta_i^1(h,\Omega)u=\Delta_i(h,\Omega)u\cr
&\qquad\qquad\ =u(x_1,\dots,x_{i-1},x_i+
h,x_{i+1},\dots,x_n)-u(x_1,\dots,x_n),\cr
&\Delta_i^j(h,\Omega)u=\Delta_i(h,\Omega)\Delta_i^{j-1}(h,\Omega)u,\cr
&{\rm and}\cr
&\Delta_i^j(h,\Omega)u=0\quad {\rm for}\quad x+jh\not\in\Omega.\cr}
$$

In [12] it has been proved that the norms of the Besov space 
$B_{p,p_0}^l(\Omega)$ are equivalent for different $m$ and $k$ satisfying 
the condition $m>l-k>0$.

By $c$ we denote a generic positive constant which changes its value from 
formula to formula and depends at most on the imbedding constants, constants 
of the considered problem and the regularity of the boundary.

By $\varphi=\varphi(\sigma_1,\dots,\sigma_k)$, $k\in\N$, we denote a generic 
function which is a positive increasing function of its arguments 
$\sigma_1,\dots,\sigma_k$, and may change its form from formula to formula.

\subsection{Auxiliary results}

We need the following interpolation lemma

\begin{lemma}\label{3.1.} [1, Chap. 4, Sect. 18] 
Let $u\in W_{p,p_0}^{s,s/2}(\Omega^T)$, $s\in\R_+$, 
$p,p_0\in[1,\infty]$, $\Omega\subset\R^3$. Let $\sigma\in\R_+\cup\{0\}$, and
$$
\kappa={3\over p}+{2\over p_0}-{3\over q}-{2\over q_0}+|\alpha|+
2a+\sigma<s.
$$
Then $D_x^\alpha\partial_t^au\in W_{q,q_0}^{\sigma,\sigma/2}(\Omega^T)$, 
$q\ge p$, $q_0\ge p_0$, and there exists $\varepsilon\in(0,1)$ such that
$$
\|D_x^\alpha\partial_t^au\|_{W_{q,q_0}^{\sigma,\sigma/2}(\Omega^T)}\le
\varepsilon^{s-\kappa}\|u\|_{W_{p,p_0}^{s,s/2}(\Omega^T)}+c
\varepsilon^{-\kappa}\|u\|_{L_{p,p_0}(\Omega^T)}.
$$
\end{lemma}

We recall from [4] the trace and the inverse trace theorems for Sobolev 
spaces with a mixed norm.

\begin{lemma}\label{3.2.} (Traces in $W_{p,p_0}^{s,s/2}(\Omega^T)$)
\begin{romannum}
\item Let $u\in W_{p,p_0}^{s,s/2}(\Omega^T)$, $s\in\R_+$, 
$p,p_0\in(1,\infty)$.\\ 
Then $u(x,t_0)\equiv u(x,t)|_{t=t_0}$ for $t_0\in[0,T]$, 
belongs to $B_{p,p_0}^{s-2/p_0}(\Omega)$, and
$$
\|u(\cdot,t_0)\|_{B_{p,p_0}^{s-2/p_0}(\Omega)}\le c
\|u\|_{W_{p,p_0}^{s,s/2}(\Omega^T)},
$$
where constant $c$ does not depend on $u$.
\item For a given $\tilde u\in B_{p,p_0}^{s-2/p_0}(\Omega)$, $s\in\R_+$, 
$s>2/p_0$, $p,p_0\in(1,\infty)$, there exists a function 
$u\in W_{p,p_0}^{s,s/2}(\Omega^T)$ such that $u|_{t=t_0}=\tilde u$ for 
$t_0\in[0,T]$, and
$$
\|u\|_{W_{p,p_0}^{s,s/2}(\Omega^T)}\le c
\|\tilde u\|_{B_{p,p_0}^{s-2/p_0}(\Omega)},
$$
where constant $c$ does not depend on $\tilde u$.
\end{romannum}
\end{lemma}

We recall also (see [1]) that if $l>1/p$ then every function from 
$B_{p,p_0}^l(\Omega)$ has a trace on the boundary $S$ belonging to 
$B_{p,p_0}^{l-1/p}(S)$, and
$$
\|u\|_{B_{p,p_0}^{l-1/p}(S)}\le c\|u\|_{B_{p,p_0}^l(\Omega)}.
$$

We apply the following imbeddings between Besov spaces.

\begin{lemma}\label{3.3.} [25, Theorem 4.6.1] 
Let $\Omega\subset\R^n$ be an arbitrary domain.
\begin{itemize}
\item[(a)] Let $s\in\R_+$, $\varepsilon>0$, $p\in(1,\infty)$ and 
$1\le q_1\le q_2\le\infty$. Then
$$
B_{p,\infty}^{s+\varepsilon}(\Omega)\subset B_{p,1}^s(\Omega)\subset
B_{p,q_1}^s(\Omega)\subset B_{p,q_2}^s(\Omega)\subset
B_{p,\infty}^s(\Omega)\subset B_{p,1}^{s-\varepsilon}(\Omega).
$$
\item[(b)] Let $\infty>q\ge p>1$, $1\le r\le\infty$, $0\le t\le s<\infty$, and 
$$
t+{n\over p}-{n\over  q}\le s.
$$
Then
$$
B_{p,r}^s(\Omega)\subset B_{q,r}^t(\Omega).
$$
\end{itemize}
\end{lemma}

We recall now from [19] a result on the solvability of a linear 
parabolic system with elasticity operator $\Q$ in Sobolev space with a mixed 
norm. This result will be repeatedly used in Section 4 in deriving a priori 
estimates for viscoelasticity system (1.1).
It generalizes the result by Krylov [13] from the single parabolic equation 
to the following parabolic system
$$\eqal{
&\u_t-\Q\u=\f\quad &{\rm in}\ \ \Omega^T=\Omega\times(0,T),\cr
&\u=\0\quad &{\rm on}\ \ S^T=S\times(0,T),\cr
&\u|_{t=0}=\u_0\quad &{\rm in}\ \ \Omega,\cr}
\leqno(3.1)
$$
where $\Omega\subset\R^3$, $S=\partial\Omega$, $\f=(f_i)$, and
$$
\Q\u=\mu\Delta\u+\nu\nabla(\nabla\cdot\u)
$$
with $\mu>0$, $\nu>0$. Let us notice that letting
$$
\Q\equiv\Q_1,\quad \mu\equiv\mu_1,\quad \nu\equiv\lambda_1+\mu_1,
$$
the assumption (1.6) implies that $\mu>0$ and $\nu>0$.

\begin{lemma}\label{3.4.} (Parabolic system in $\W_{p,p_0}^{2,1}(\Omega^T)$  
[13, 19, 24])
\begin{romannum}
\item Assume that $\f\in{\bm L}_{p,p_0}(\Omega^T)$, 
$\u_0\in\B_{p,p_0}^{2-2/p_0}(\Omega)$, $p,p_0\in(1,\infty)$, $S\in C^2$. 
If $2-2/p_0-1/p>0$ the compatibility condition $\u_0|_S=\0$ is assumed. Then 
there exists a unique solution to problem (3.1) such that 
$\u\in\W_{p,p_0}^{2,1}(\Omega^T)$ and
$$
\|\u\|_{\W_{p,p_0}^{2,1}(\Omega^T)}\le c(\|\f\|_{{\bm L}_{p,p_0}(\Omega^T)}+
\|\u_0\|_{\B_{p,p_0}^{2-2/p_0}(\Omega)})
\leqno(3.2)
$$
with a constant $c$ depending on $\Omega$, $S$, $p$, $p_0$.
\item Assume that 
$\f=\nabla\cdot\g+\b$, $\g=(g_{ij})$, $\b=(b_i)$, 
$\g,\b\in{\bm L}_{p,p_0}(\Omega^T)$, $\u_0\in\B_{p,p_0}^{2-2/p_0}(\Omega)$. 
Assume the compatibility condition
$$
\u_0|_S=\0\quad \textsl{if}\quad 1-2/p_0-1/p>0.
$$
Then there exists a unique solution to (3.1) such that 
$\u\in\W_{p,p_0}^{1,1/2}(\Omega^T)$ and
$$\eqal{
&\|\u\|_{\W_{p,p_0}^{1,1/2}(\Omega^T)}\le c(\|\g\|_{{\bm L}_{p,p_0}(\Omega^T)}+
\|\b\|_{{\bm L}_{p,p_0}(\Omega^T)}\cr
&\quad+\|\u_0\|_{\B_{p,p_0}^{2-2/p_0}(\Omega)})\cr}
\leqno(3.3)
$$
with a constant $c$ depending on $\Omega$, $S$, $p$, $p_0$.
\end{romannum}
\end{lemma}

In the proof of Theorem A we shall apply also the following regularity 
result for a linear parabolic equation. This result is the special case of 
a more general theorem due to Denk-Hieber-Pr\"uss [7, Theorem 2.3].

\begin{lemma}\label{3.5.} (Parabolic equation in $W_{p,p_0}^{2,1}(\Omega^T)$) 
Let us consider the problem
$$\eqal{
&\theta_t-\varrho\Delta\theta=g\quad &\textsl{in}\ \ \Omega^T,\cr
&\n\cdot\nabla\theta=0\quad &\textsl{on}\ \ S^T,\cr
&\theta|_{t=0}=\theta_0\quad &\textsl{in}\ \ \Omega,\cr}
\leqno(3.4)
$$
where $\varrho(x,t)$ is a continuous function on $\Omega^T$ such that 
$\inf_\Omega\varrho>0$. Assume that 
$g\in L_{p,p_0}(\Omega^T)$, $\theta_0\in B_{p,p_0}^{2-2/p_0}(\Omega)$, 
$p,p_0\in(1,\infty)$, $S\in C^2$, and the corresponding compatibility 
condition is satisfied. Then there exists a unique solution to problem (3.4) 
such that $\theta\in W_{p,p_0}^{2,1}(\Omega^T)$ and
$$
\|\theta\|_{W_{p,p_0}^{2,1}(\Omega^T)}\le c(\|g\|_{L_{p,p_0}(\Omega^T)}+
\|\theta_0\|_{B_{p,p_0}^{2-2/p_0}(\Omega)})
\leqno(3.5)
$$
with a constant $c$ depending on $\Omega$, $T$, $S$, 
$\inf_{\Omega^T}\varrho$ and $\sup_{\Omega^T}\varrho$.
\end{lemma}

\begin{remark}\label{3.6.} 
The constants $c$ in Lemmas 3.4 and 3.5 do not depend on $T$. For $T$ small 
the proof of this fact is evident whereas for $T$ large it can be deduced by 
applying the same arguments as in the proof of Theorem 3.1.1 in [26, Ch. 3].
\end{remark}

\section{A priori estimates}

In this section we derive a priori estimates for solutions of problem 
(1.1)--(1.4) on an arbitrary finite time interval $(0,T)$. The estimates are 
essential for the existence proof by the successive approximation method, 
presented in Section 5.

The procedure of deriving a priori estimates consists in a recursive 
improvement of the basic energy estimates. The main tool used in this 
procedure is Lemma 3.4 which provides the solvability of the viscoelasticity 
system $(1.8)_1$ in the Sobolev space $W_{p,p_0}^{2,1}(\Omega^T)$. The applied 
procedure is aimed to establish the continuity of temperature $\theta$ and 
finally to apply Lemma 3.5 on the solvability of parabolic equation in Sobolev 
space $W_{q,q_0}^{2,1}(\Omega^T)$.

Throughout this section we assume that assumptions (A1)--(A3) (see Sect. 1.3) 
hold, and
$$
\theta_0\ge\underline{\theta}>0\quad {\rm in}\ \ \Omega,\quad g\ge0\quad 
{\rm in}\ \ \Omega^T,
\leqno(4.1)
$$
where $\underline\theta$ is a positive constant.

First we prove the lower bound on $\theta$ by using similar arguments as in 
[16, Lemma 3.7], [27, Lemma 3.3].

\begin{lemma}\label{4.1.} (Lower bound on $\theta$) 
Let us assume that (4.1) holds. Then there exists a positive 
constant $c$ depending only on parameters $a_{1*}$, $a_2^*$ (from (1.10), 
$|\alphaa|$ (see A2), $c_v$ (see 1.2)), such that
$$
\theta(t)\ge{\underline\theta}\exp(-cT)\equiv\theta_*>0\quad \textsl{for}\ \ 
t\in[0,T].
\leqno(4.2)
$$
\end{lemma}

\begin{proof} 
For $m\in\R_+$ let us define the truncation
$$
\theta_m=\max\left\{\theta,{1\over m}\right\}
$$
and
$$\Omega_m(t)=\left\{x\in\Omega:\ \theta(x,t)>{1\over m}\right\}.
$$
Multiplying (1.2) by $-\theta_m^{-\varrho}$ with $\varrho>2$ (admissible 
test function) and integrating over $\Omega_m(t)$ gives
$$\eqal{
&-c_v\intop_{\Omega_m(t)}\theta\theta_t\theta_m^{-\varrho}dx+
k\intop_{\Omega_m(t)}\theta_m^{-\varrho}\Delta\theta dx+\intop_{\Omega_m(t)}
{(\A_1\epsi_t)\cdot\epsi_t\over\theta_m^\varrho}dx\cr
&\quad+\intop_{\Omega_m(t)}{g\over\theta_m^\varrho}dx=\intop_{\Omega_m(t)}
{\theta\over\theta_m^\varrho}(\A_2\alphaa)\cdot\epsi_tdx.\cr}
\leqno(4.3)
$$
The first term on the left-hand side of (4.3) is equal to
$$\eqal{
-c_v\intop_{\Omega_m(t)}\theta_m\theta_{m,t}\theta_m^{-\varrho}dx&=
{c_v\over\varrho-2}\intop_{\Omega_m(t)}\partial_t\theta_m^{2-\varrho}dx\cr
&={c_v\over\varrho-2}\intop_\Omega\partial_t\theta_m^{2-\varrho}dx=
{c_v\over\varrho-2}{d\over dt}\intop_\Omega\theta_m^{2-\varrho}dx,\cr}
\leqno(4.4)
$$
because $\partial_t\theta_m^{2-\varrho}=0$ for 
$x\in\Omega\setminus\Omega_m(t)=\{x\in\Omega:\ \theta_m(t)={1\over m}\}$.

\noindent
The second term on the left-hand side of (4.3) equals
$$
k\intop_{\Omega_m(t)}\theta_m^{-\varrho}\Delta\theta_mdx=k\intop_\Omega
\theta_m^{-\varrho}\Delta\theta_mdx={4k\varrho\over(\varrho-1)^2}\intop_\Omega
\bigg|\nabla\bigg({1\over\theta_m^{\varrho-1\over2}}\bigg)\bigg|^2dx,
\leqno(4.5)
$$
because $\nabla\theta_m=\nabla\theta$ for $x\in\Omega_m(t)$ and 
$\nabla\theta_m=0$ for $x\in\Omega\setminus\Omega_m(t)$. On account of (1.10) 
the third term on the left-hand side of (4.3) is bounded from below by
$$
\intop_{\Omega_m(t)}{(\A_1\epsi_t)\cdot\epsi_t\over\theta_m^\varrho}dx\ge 
a_{1*}\intop_{\Omega_m(t)}{|\epsi_t|^2\over\theta_m^\varrho}dx,
\leqno(4.6)
$$
and the fourth one by
$$
\intop_{\Omega_m(t)}{g\over\theta_m^\varrho}dx\ge0.
\leqno(4.7)
$$
In view of the boundedness of $\A_2$ and $\alphaa$, the integral on the 
right-hand side of (4.3) is estimated as follows:
$$\eqal{
\lc{\intop_{\Omega_m(t)}{\theta\over\theta_m^\varrho}(\A_2\alphaa)\cdot\epsi_tdx=
\intop_{\Omega_m(t)}{\theta_m\over\theta_m^{\varrho/2}}(\A_2\alphaa)\cdot
{\epsi_t\over\theta_m^{\varrho/2}}dx}
\rc{\le{\delta\over2}\intop_{\Omega_m(t)}{|\epsi_t|^2\over\theta_m^\varrho}dx+
{c\over2\delta}\intop_{\Omega_m(t)}{\theta_m^2\over\theta_m^\varrho}dx,\quad 
\delta>0.}}
\leqno(4.8)
$$
Now, setting $\delta/2=a_{1*}$ and incorporating (4.4)--(4.8) into (4.3) 
we arrive at
$$\eqal{
\lc{{c_v\over\varrho-2}{d\over dt}\intop_\Omega\theta_m^{2-\varrho}dx+
{4k\varrho\over(\varrho-1)^2}\intop_\Omega\bigg|\nabla\bigg(
{1\over\theta_m^{\varrho-1\over2}}\bigg)\bigg|^2dx}
\rc{\le c\intop_{\Omega_m(t)}\theta_m^{2-\varrho}dx\le c\intop_\Omega
\theta_m^{2-\varrho}dx,}}
$$
where in the last inequality we used the fact that $\theta_m>0$ in $\Omega$. 
Hence, by the Gronwall inequality, it follows that
$$
\intop_\Omega\theta_m^{2-\varrho}(t)dx\le\intop_\Omega\theta_m^{2-\varrho}(0)
dx\exp\bigg[{c(\varrho-2)\over c_v}t\bigg]\quad {\rm for}\ \ t\in[0,T],
$$
that is,
$$
\|\theta_m^{-1}(t)\|_{L_{\varrho-2}(\Omega)}\le
\|\theta_m^{-1}(0)\|_{L_{\varrho-2}(\Omega)}\exp(cT)
\leqno(4.9)
$$
with a constant $c$ independent of $\varrho$ and $m$. Letting 
$\varrho\to\infty$, (4.9) yields the bound
$$
\theta_m(t)\ge\theta_m(0)\exp(-cT).
$$
Further, letting $m\to\infty$ and noting that for sufficiently large $m$, 
$\theta_m(0)=\break=\max\left\{\theta_0,{1\over m}\right\}\ge{\underline\theta}$, 
we conclude the bound (4.2).
\end{proof}

\begin{lemma}\label{4.2.} (Energy estimates) 
Let us assume that (4.1) holds, $\theta>0$, and
$$
\u_0\in\H^1(\Omega),\ \ \u_1\in{\bm L}_2(\Omega),\ \ \theta_0\in L_2(\Omega),\ \ 
\b\in{\bm L}_{2,1}(\Omega^T),\ \ g\in L_1(\Omega^T).
$$
Then a sufficiently smooth solution $(\u,\theta)$ to (1.1)--(1.4) satisfies 
the estimate
$$\eqal{
&\|\u_{t'}\|_{{\bm L}_{2,\infty}(\Omega^t)}+\|\epsi\|_{{\bm L}_{2,\infty}(\Omega^t)}+
\|\theta\|_{L_{2,\infty}(\Omega^t)}\cr
&\quad+\|\theta^{-1}\nabla\theta\|_{{\bm L}_2(\Omega^t)}+
\|\theta^{-1/2}\epsi_{t'}\|_{{\bm L}_2(\Omega^t)}\cr
&\le c(\|\epsi(\u_0)\|_{{\bm L}_2(\Omega)}+\|\u_1\|_{{\bm L}_2(\Omega)}+
\|\theta_0\|_{L_2(\Omega)}\cr
&\quad+\|\b\|_{{\bm L}_{2,1}(\Omega^t)}+\|g\|_{L_1(\Omega^t)}+1)\equiv cA_0,\cr}
\leqno(4.10)
$$
where $t\le T$.
\end{lemma}

\begin{proof} 
Note that the positivity of $\theta$ is ensured by Lemma 4.1. Multiplying 
(1.1) by $\u_t$, integrating over $\Omega$ and integrating by parts using 
boundary condition $(1.3)_1$, gives
$$\eqal{
&{1\over2}{d\over dt}\intop_\Omega|\u_t|^2dx+\intop_\Omega(\A_1\epsi_t)\cdot
\epsi_tdx+\intop_\Omega\A_2(\epsi-\theta\alphaa)\cdot\epsi_tdx\cr
&=\intop_\Omega\b\cdot\u_tdx.\cr}
\leqno(4.11)
$$
Further, integrating (1.2) over $\Omega$ and by parts using $(1.3)_2$ yields
$$
{c_v\over2}{d\over dt}\intop_\Omega\theta^2dx+\intop_\Omega\theta(\A_2\alphaa)
\cdot\epsi_tdx-\intop_\Omega(\A_1\epsi_t)\cdot\epsi_tdx=\intop_\Omega gdx.
\leqno(4.12)
$$
Next, multiplying (1.2) by $1/\theta$, integrating over $\Omega$ and 
integrating by parts leads to
$$
{d\over dt}\intop_\Omega[c_v\theta+(\A_2\alphaa)\cdot\epsi]dx-k\intop_\Omega
{|\nabla\theta|^2\over\theta^2}dx-\intop_\Omega
{(\A_1\epsi_t)\cdot\epsi_t\over\theta}dx=\intop_\Omega{g\over\theta}dx.
\leqno(4.13)
$$
Adding by sides (4.11) and (4.12) gives
$$
{d\over dt}\intop_\Omega\bigg[{1\over2}c_v\theta^2+{1\over2}(\A_2\epsi)\cdot
\epsi+{1\over2}|\u_t|^2\bigg]dx=\intop_\Omega(\b\cdot\u_t+g)dx.
\leqno(4.14)
$$
Now, multiplying (4.13) by a positive constant $\beta$, and subtracting by 
sides from (4.14), we get
$$\eqal{
&{d\over dt}\intop_\Omega\bigg[{1\over2}c_v\theta^2+{1\over2}(\A_2\epsi)\cdot
\epsi+{1\over2}|\u_t|^2-\beta c_v\theta-\beta(\A_2\alphaa)\cdot\epsi\bigg]dx
\cr
&\quad+k\beta\intop_\Omega{|\nabla\theta|^2\over\theta^2}dx+\beta\intop_\Omega
{(\A_1\epsi_t)\cdot\epsi_t\over\theta}dx=\intop_\Omega\bigg[\b\cdot\u_t+
\bigg(1-{\beta\over\theta}\bigg)g\bigg]dx.\cr}
\leqno(4.15)
$$
Integrating (4.15) with respect to time, noting that by the boundedness of 
$\A_2$ and $\alphaa$,
$$\eqal{
\lc{{1\over2}c_v\theta^2+{1\over2}(\A_2\epsi)\cdot\epsi+{1\over2}|\u_t|^2-
\beta c_v\theta-\beta(\A_2\alphaa)\cdot\epsi}
\rc{\ge{1\over4}c_v\theta^2+{1\over4}(\A_2\epsi)\cdot\epsi+{1\over2}|\u_t|^2-c,}}
$$
and using the coercivity of $\A_1$ (see (1.10)), we get
$$\eqal{
&\|\theta(t)\|_{L_2(\Omega)}^2+\|\epsi(t)\|_{{\bm L}_2(\Omega)}^2+
\|\u_t(t)\|_{{\bm L}_2(\Omega)}^2+\|\theta^{-1}\nabla\theta\|_{{\bm L}_2(\Omega^t)}^2
\cr
&\quad+\|\theta^{-1/2}\epsi_{t'}\|_{{\bm L}_2(\Omega^t)}^2\le c
(\|\theta_0\|_{L_2(\Omega)}^2+\|\epsi(\u_0)\|_{{\bm L}_2(\Omega)}^2+
\|\u_1\|_{{\bm L}_2(\Omega)}^2)\cr
&\quad+\|\u_{t'}\|_{{\bm L}_{2,\infty}(\Omega^t)}
\|\b\|_{{\bm L}_{2,1}(\Omega^t)}+\|g\|_{L_1(\Omega^t)}+c\cr}
$$
for $t\le T$, with a constant $c$ depending only on parameters.

\noindent
Now, applying the Young inequality to the second term on the right-hand 
side of the above inequality yields (4.10).
\end{proof}

\begin{remark}\label{4.3.} 
By integrating the identity (4.14) with respect to time one can immediately 
conclude that
$$
\|\u_{t'}\|_{{\bm L}_{2,\infty}(\Omega^t)}+\|\epsi\|_{{\bm L}_{2,\infty}(\Omega^t)}+
\|\theta\|_{L_{2,\infty}(\Omega^t)}\le cA_0,
\leqno(4.16)
$$
which is a part of inequality (4.10).

\noindent
On the contrary to (4.10) this inequality does not require the assumption 
$\theta>0$. We point out that in deriving further estimates we shall use just 
the bounds in (4.16) loosing the information contained in the two dissipative 
terms of (4.10).
This information may be of importance in the analysis of the long time 
behaviour of solutions.
\end{remark}

\begin{remark}\label{4.4.} 
We complement Lemma 4.2 by some physical interpretations. In view of 
$(2.10)_1$, (2.12) and boundary conditions (1.3), identity (4.14) represents 
the balance equation for the total energy
$$
{d\over dt}\intop_\Omega\bigg(e+{1\over2}|\u_t|^2\bigg)dx+\intop_\Omega
[-(\S\n)\cdot\u_t+\n\cdot\q]dS=\intop_\Omega(\b\cdot\u_t+g)dx.
$$
On the other hand, in view of $(2.10)_2$, $(2.12)_2$ and the boundary 
condition $(1.3)_2$, identity (4.13) represents the balance equation for 
the entropy
$$
{d\over dt}\intop_\Omega\eta dx+\intop_S\n\cdot{\q\over\theta}dS=\intop_\Omega
\sigma dx+\intop_\Omega{g\over\theta}dx
$$
with the entropy production
$$
\sigma={k\over\theta^2}|\nabla\theta|^2+{1\over\theta}(\A_1\epsi_t)\cdot
\epsi_t\ge0.
$$

Equation (4.15) represents the so-called availability identity
$$\eqal{
&{d\over dt}\intop_\Omega\bigg(e+{1\over2}|\u_t|^2-\beta\eta\bigg)dx+
\intop_S\bigg[-(\S\n)\cdot\u_t+\bigg(1-{\beta\over\theta}\bigg)\n\cdot
\q\bigg]dS\cr
&\quad+\beta\intop_\Omega\Sigma dx=\intop_\Omega\bigg[\b\cdot\u_t+
\bigg(1-{\beta\over\theta}\bigg)g\bigg]dx,\cr}
\leqno(4.17)
$$
where $\beta=\const>0$. Hence, since $\sigma\ge0$, it follows that if the 
external sources vanish
$$
\b=\0,\quad g=0,
$$
and if the boundary conditions on $S$ imply that
$$
(\S\n)\cdot\u_t=0,\quad \bigg(1-{\beta\over\theta}\bigg)\n\cdot\q=0,
\leqno(4.18)
$$
then
$$
{d\over dt}\intop_\Omega\bigg(e+{1\over2}|\u_t|^2-\beta\eta\bigg)dx\le0.
$$
This provides the Lyapunov functional 
$\intop_\Omega\left(e+{1\over2}|\u_t|^2-\beta\eta\right)dx$, which is 
nonincreasing on solutions paths.
Let us notice that the boundary conditions (1.3) ensure (4.18).
The identity (4.17) has been used in deriving energy estimates in Lemma 4.2.
\end{remark}
\vskip6pt

Our goal now is to derive further regularity properties from the energy 
estimates (4.10). 
To this purpose we use the regularity results for parabolic systems in 
Sobolev space with a mixed norm, stated in Lemmas 3.4 and 3.5.

Let us consider the viscoelasticity system (1.1) with boundary and initial 
conditions $(1.3)_1$, $(1.4)_1$, expresses in the form:
$$\eqal{
&\u_{tt}-\Q\u_t=\nabla\cdot[\A_2(\epsi-\theta\alphaa)]+\b\quad &{\rm in}\ \ 
\Omega^T,\cr
&\u=\0\quad &{\rm on}\ \ S^T,\cr
&\u|_{t=0}=\u_0,\ \ \u_t|_{t=0}=\u_1\quad &{\rm in}\ \ \Omega,\cr}
\leqno(4.19)
$$
where $\Q=\Q_1$ is the viscosity operator (1.7).

Applying Lemma 3.4 (i), (ii) to system (4.19) we deduce, in view of the 
boundedness of $\A_2$ and $\alphaa$, the following

\begin{corollary}\label{4.5.} 
Let us assume that
$$
\u_1\in\B_{p,\sigma}^{2-2/\sigma}(\Omega),\quad \b\in{\bm L}_{p,\sigma}(\Omega^T),
\quad p,\sigma\in(1,\infty),
$$
and if $2-2/\sigma-1/p>0$ then the compatibility condition $\u_0|_S=\0$ holds.
\begin{romannum}
\item If $\epsi\in{\bm L}_{p,\sigma}(\Omega^T)$ and 
$\theta\in L_{p,\sigma}(\Omega^T)$, $p,\sigma\in(1,\infty)$, then the 
solution $\u$ to problem (4.19) satisfies
$$\eqal{
\lc{\|\epsi_{t'}\|_{{\bm L}_{p,\sigma}(\Omega^t)}\le c
\|\u_{t'}\|_{\W_{p,\sigma}^{1,1/2}(\Omega^T)}\le c
(\|\epsi\|_{{\bm L}_{p,\sigma}(\Omega^t)}+
\|\theta\|_{L_{p,\sigma}(\Omega^t)}}
\rc{\quad+\|\b\|_{{\bm L}_{p,\sigma}(\Omega^t)}+
\|\u_1\|_{\B_{p,\sigma}^{2-2/\sigma}(\Omega)})}}
\leqno(4.20)
$$
for $t\in(0,T]$, with a constant $c$ depending on $\Omega$, $S$, $T$, $p$ 
and $\sigma$.
\item If $\nabla\epsi\in{\bm L}_{p,\sigma}(\Omega^T)$ and 
$\nabla\theta\in{\bm L}_{p,\sigma}(\Omega^T)$, $p,\sigma\in(1,\infty)$, 
then the solution $\u$ to (4.19) satisfies
$$\eqal{
\lc{\|\epsi_{t'}\|_{\W_{p,\sigma}^{1,1/2}(\Omega^t)}\le c
\|\u_{t'}\|_{\W_{p,\sigma}^{2,1}(\Omega^t)}\le c
(\|\nabla\epsi\|_{{\bm L}_{p,\sigma}(\Omega^t)}+
\|\nabla\theta\|_{{\bm L}_{p,\sigma}(\Omega^t)}}
\rc{\quad+\|\b\|_{{\bm L}_{p,\sigma}(\Omega^t)}+
\|\u_1\|_{\B_{p,\sigma}^{2-2/\sigma}(\Omega)})}}
\leqno(4.21)
$$
for $t\in(0,T]$, with a constant $c$ depending on $\Omega$, $S$, $T$, $p$ 
and $\sigma$.
\end{romannum}
\end{corollary}

\noindent
Using (4.10) in (4.20) for $p=2$ and $\sigma$ arbitrary finite we have
$$\eqal{
\lc{\|\epsi_{t'}\|_{{\bm L}_{2,\sigma}(\Omega^t)}\le c(A_0+
\|\u_1\|_{\B_{2,\sigma}^{2-2/\sigma}(\Omega)}+
\|\b\|_{{\bm L}_{2,\sigma}(\Omega^t)})}
\rc{\equiv cA_1(\sigma),\quad \sigma\in(1,\infty),\quad t\le T,}}
\leqno(4.22)
$$
where $A_0$ is defined in (4.10).

For further purposes (see the proof of Lemma 4.7) we prepare now some 
inequalities between the norms of $\u$ and $\theta$.

Let us consider viscoelasticity system (4.19) rewritten in the form
$$\eqal{
&\u_{tt}-\Q\u_t-\Q\u=\nabla\cdot[(\A_2-\A_1)\epsi-\theta\A_2\alphaa]+\b\qquad 
&{\rm in}\ \ \Omega^T,\cr
&\u=\0\quad &{\rm on}\ \ S^T,\cr
&\u|_{t=0}=\u_0,\ \ \u_t|_{t=0}=\u_1\quad &{\rm in}\ \ \Omega,\cr}
\leqno(4.23)
$$
where $\Q\equiv\Q_1$.

\begin{lemma}\label{4.6.} 
\begin{romannum}
\item Let $\u_0\in\H_0^1(\Omega)$, $\u_1\in{\bm L}_2(\Omega)$, 
$\b\in{\bm L}_2(\Omega^T)$ and $\theta\in L_2(\Omega^T)$. Then a solution $\u$ 
to system (4.23) satisfies
$$\eqal{
&\|\u_{t'}\|_{{\bm L}_{2,\infty}(\Omega^t)}+
\|\Q^{1/2}\u\|_{{\bm L}_{2,\infty}(\Omega^t)}+
\|\Q^{1/2}\u_{t'}\|_{{\bm L}_2(\Omega^t)}\cr
&\le c(t)(\|\theta\|_{L_2(\Omega^t)}+\|\Q^{1/2}\u_0\|_{{\bm L}_2(\Omega)}+
\|\u_1\|_{{\bm L}_2(\Omega)}\cr
&\quad+\|\b\|_{{\bm L}_2(\Omega^t)})\qquad \textsl{for}\ \ t\le T,\cr}
\leqno(4.24)
$$
with a constant $c(t)$ exponentially depending on $t$.
\item Let $\u_0\in\H^2(\Omega)\cap\H_0^1(\Omega)$, 
$\u_1\in\H_0^1(\Omega)$, $\b\in L_2(\Omega^T)$ and 
$\nabla\theta\in{\bm L}_2(\Omega)$. Then a solution $\u$ to system (4.23) 
satisfies
$$\eqal{
&\|\Q^{1/2}\u_{t'}\|_{{\bm L}_{2,\infty}(\Omega^t)}+
\|\Q\u\|_{{\bm L}_{2,\infty}(\Omega^t)}+\|\Q\u_{t'}\|_{{\bm L}_2(\Omega^t)}\cr
&\le c(t)(\|\nabla\theta\|_{{\bm L}_2(\Omega^t)}+
\|\Q^{1/2}\u_1\|_{{\bm L}_2(\Omega)}+\|\Q\u_0\|_{{\bm L}_2(\Omega)}\cr
&\quad+\|\b\|_{{\bm L}_2(\Omega^t)})\qquad \textsl{for}\ \ t\le T,\cr}
\leqno(4.25)
$$
with a constant $c(t)$ exponentially depending on $t$.
\end{romannum}
\end{lemma}

\begin{proof}
(i) Multiplying $(4.23)_1$ by $\u_t$, integrating over $\Omega$ and using the 
boundary condition $(4.23)_2$ gives
$$\eqal{
&{1\over2}{d\over dt}\intop_\Omega(|\u_t|^2+|\Q^{1/2}\u|^2)dx+\intop_\Omega
|\Q^{1/2}\u_t|^2dx\cr
&=\intop_\Omega[(\A_1-\A_2)\epsi+\theta\A_2\alphaa]\cdot\epsi_tdx+
\intop_\Omega\b\cdot\u_tdx.\cr}
\leqno(4.26)
$$
Using the estimate
$$
\intop_\Omega[(\A_1-\A_2)\epsi+\theta\A_2\alphaa]\cdot\epsi_tdx\le
\delta_1\intop_\Omega|\Q^{1/2}\u_t|^2dx
+c(1/\delta_1)\intop_\Omega(|\epsi|^2+\theta^2)dx,
$$
which results on account of the Young inequality and (1.16), we conclude that
$$\eqal{
&{d\over dt}\intop_\Omega(|\u_t|^2+|\Q^{1/2}\u|^2)dx+\intop_\Omega
|\Q^{1/2}\u_t|^2dx\cr
&\le c\intop_\Omega(\theta^2+|\b|^2)dx+c_1\intop_\Omega(|\u_t|^2+
|\Q^{1/2}\u|^2)dx,\cr}
\leqno(4.27)
$$
where we distinguished the constant $c_1$.

\noindent
Hence, omitting the last integral on the left-hand side, it follows that
$$
{d\over dt}\bigg[\intop_\Omega(|\u_t|^2+|\Q^{1/2}\u|^2)dxe^{-c_1t}\bigg]\le 
ce^{-c_1t}\intop_\Omega(\theta^2+|\b|^2)dx,
$$
which after integrating with respect to $t'\in(0,t)$ gives
$$\eqal{
&\intop_\Omega(|\u_t|^2+|\Q^{1/2}\u|^2)dx\le ce^{c_1t}\intop_{\Omega^t}
(\theta^2+|\b|^2)dxdt'\cr
&\quad+e^{c_1t}\intop_\Omega(|\u_1|^2+|\Q^{1/2}\u_0|^2)dx.\cr}
\leqno(4.28)
$$
Now, using (4.28) in (4.27) and again integrating the result with respect to 
$t'\in(0,t)$ leads to
$$\eqal{
&\intop_\Omega(|\u_t|^2+|\Q^{1/2}\u|^2)dx+\intop_{\Omega^t}|\Q^{1/2}\u_t|^2
dxdt'\cr
&\le c(te^{c_1t}+1)\intop_{\Omega^t}(\theta^2+|\b|^2)dxdt'+(te^{c_1t}+1)
\intop_\Omega(|\u_1|^2+|\Q^{1/2}\u_0|^2)dx,\cr}
$$
which proves (4.24).
\goodbreak

\noindent
(ii) Multiplying $(4.23)_1$ by $\Q\u_t$, integrating over $\Omega$ and 
integrating by parts yields
$$\eqal{
&{1\over2}{d\over dt}\intop_\Omega(|\Q^{1/2}\u_t|^2+|\Q\u|^2)dx+
\intop_\Omega|\Q\u_t|^2dx\cr
&=\intop_\Omega(\nabla\cdot[(\A_1-\A_2)\epsi+\theta\A_2\alphaa])\cdot
\Q\u_tdx-\intop_\Omega\b\cdot\Q\u_tdx\equiv{\cal R}.\cr}
\leqno(4.29)
$$
In view of the boundedness of $\A_1$, $\A_2$ and $\alphaa$,
$$
{\cal R}\le\delta_2\intop_\Omega|\Q\u_t|^2dx+c(1/\delta_2)\intop_\Omega
(|\nabla\epsi|^2+|\nabla\theta|^2+|\b|^2)dx.
$$
Hence, choosing $\delta_2$ suffciently small and recalling the ellipticity 
of the operator $\Q$, we get
$$\eqal{
&{d\over dt}\intop_\Omega(|\Q^{1/2}\u_t|^2+|\Q\u|^2)dx+\intop_\Omega
|\Q\u_t|^2dx\cr
&\le c\intop_\Omega(|\nabla\theta|^2+|\b|^2)dx+c_2\intop_\Omega
(|\Q^{1/2}\u_t|^2+|\Q\u|^2)dx,\cr}
\leqno(4.30)
$$
where we distinguished the constant $c_2$. Omitting the last integral on the 
left-hand side the latter inequality leads to
$$
{d\over dt}\bigg[\intop_\Omega(|\Q^{1/2}\u_t|^2+|\Q\u|^2)dxe^{-c_2t}\bigg]\le
ce^{-c_2t}\intop_\Omega(|\nabla\theta|^2+|\b|^2)dx,
$$
which after integrating with respect to $t'\in(0,t)$, leads to
$$\eqal{
&\intop_\Omega(|\Q^{1/2}\u_t|^2+|\Q\u|^2)dx\le ce^{c_2t}\intop_{\Omega^t}
(|\nabla\theta|^2+|\b|^2)dxdt'\cr
&\quad+e^{c_2t}\intop_\Omega(|\Q^{1/2}\u_1|^2+|\Q\u_0|^2)dx.\cr}
\leqno(4.31)
$$
Inserting (4.31) into (4.30) and again integrating the result with respect 
to $t'\in(0,t)$ gives
$$\eqal{
&\intop_\Omega(|\Q^{1/2}\u_t|^2+|\Q\u|^2)dx+\intop_{\Omega^t}
|\Q\u_{t'}|^2dxdt'\cr
&\le c(te^{c_2t}+1)\intop_{\Omega^t}(|\nabla\theta|^2+|\b|^2)dxdt'
+(te^{c_2t}+1)\intop_\Omega(|\Q^{1/2}\u_1|^2+|\Q\u_0|^2)dx,\cr}
$$
which proves (4.25).
\end{proof}

>From (4.24) we conclude that
$$\eqal{
&\|\u_{t'}\|_{{\bm L}_{2,\infty}(\Omega^t)}+\|\u\|_{L_\infty(0,t;\H^1(\Omega))}+
\|\u_{t'}\|_{L_2(0,t;\H^1(\Omega))}\cr
&\le c(t)(\|\theta\|_{L_2(\Omega^t)}+\|\u_0\|_{\H^1(\Omega)}+
\|\u_1\|_{{\bm L}_2(\Omega)}+\|\b\|_{{\bm L}_2(\Omega^t)}),\quad t\le T.\cr}
\leqno(4.32)
$$
Similarly, from (4.25) it follows that
$$\eqal{
&\|\u_{t'}\|_{L_\infty(0,t;\H^1(\Omega))}+\|\u_{t'}\|_{L_2(0,t;\H^2(\Omega))}+
\|\u\|_{L_\infty(0,t;\H^2(\Omega))}\cr
&\le c(t)(\|\nabla\theta\|_{{\bm L}_2(\Omega^t)}+\|\u_0\|_{\H^2(\Omega)}+
\|\u_1\|_{\H^1(\Omega)}+\|\b\|_{{\bm L}_2(\Omega^t)}),\quad t\le T.\cr}
\leqno(4.33)
$$
Hence, by the definition of $\epsi$,
$$\eqal{
&\|\epsi_{t'}\|_{\V_2(\Omega^t)}\equiv
\|\epsi_{t'}\|_{{\bm L}_{2,\infty}(\Omega^t)}+
\|\epsi_{t'}\|_{{\bm L}_2(0,t;\H^1(\Omega))}\cr
&\le c(t)(\|\nabla\theta\|_{{\bm L}_2(\Omega^t)}+\|\u_0\|_{\H^2(\Omega)}+
\|\u_1\|_{\H^1(\Omega)}+\|\b\|_{{\bm L}_2(\Omega^t)}),\cr}
\leqno(4.34)
$$
where $t\le T$. With the help of this inequality we prove

\begin{lemma}\label{4.7.} 
Assume that $\u_0\in\H^2(\Omega)$, $\u_1\in\B_{2,10}^{2-1/5}(\Omega)$, 
$\theta_0\in L_3(\Omega)$, $\b\in{\bm L}_{2,10}(\Omega^t)$, 
$g\in L_{2,1}(\Omega^t)$, $t\le T$.\\
Then
$$
\|\theta(t)\|_{L_3(\Omega)}+\|\theta\|_{L_2(0,t;H^1(\Omega))}+
\|\epsi_{t'}\|_{\V_2(\Omega^t)}\le cA_2,\quad t\le T,
\leqno(4.35)
$$
where
$$\eqal{
A_2&=\varphi(A_0,A_1(10),\|\u_0\|_{\H^2(\Omega)},\|\u_1\|_{\H^1(\Omega)},
\|\theta_0\|_{L_3(\Omega)},
\|\b\|_{{\bm L}_2(\Omega^t)},\|g\|_{L_{2,1}(\Omega^t)})\cr
&\le\varphi(\|\u_0\|_{\H^2(\Omega)},\|\u_1\|_{\B_{2,10}^{2-1/5}(\Omega)},
\|\theta_0\|_{L_3(\Omega)},
\|\b\|_{{\bm L}_{2,10}(\Omega^t)},\|g\|_{L_{2,1}(\Omega^t)})\equiv A_3,\cr}
\leqno(4.36)
$$
with $A_1(\cdot)$ defined in (4.22) and $A_0$ in (4.10).
\end{lemma}

\begin{proof} 
Multiplying (1.2) by $\theta$ and integrating over $\Omega^t$ we get
$$\eqal{
&\|\theta(t)\|_{L_3(\Omega)}^3+\intop_0^t
\|\nabla\theta(t')\|_{{\bm L}_2(\Omega)}^2dt'\cr
&\le c\intop_0^t\intop_\Omega\theta^2|\epsi_{t'}|dxdt'+c\intop_0^t
\intop_\Omega\theta|\epsi_{t'}|^2dxdt'
+c\intop_0^t\intop_\Omega\theta|g|dxdt'+
\|\theta_0\|_{L_3(\Omega)}^3.\cr}
\leqno(4.37)
$$
With the use of the H\"older inequality the first term on the right-hand 
side of (4.37) is estimated by
$$\eqal{
&\intop_0^t\|\theta\|_{L_3(\Omega)}\|\theta\|_{L_6(\Omega)}
\|\epsi_{t'}\|_{{\bm L}_2(\Omega)}dt'\cr
&\le\sup_t\|\theta(t)\|_{L_3(\Omega)}\intop_0^t\|\theta(t')\|_{L_6(\Omega)}
\|\epsi_{t'}(t')\|_{{\bm L}_2(\Omega)}dt'\cr
&\le\sup_t\|\theta(t)\|_{L_3(\Omega)}\|\theta\|_{L_{6,2}(\Omega^t)}
\|\epsi_{t'}\|_{{\bm L}_2(\Omega^t)}\cr
&\le\delta_1\sup_t\|\theta(t)\|_{L_3(\Omega)}^3+c(1/\delta_1)
\|\theta\|_{L_{6,2}(\Omega^t)}^{3/2}\|\epsi_{t'}\|_{{\bm L}_2(\Omega^t)}^{3/2}
\equiv I_1.\cr}
$$
Using (4.22) for $\sigma=2$ yields
$$
I_1\le\delta_1\sup_t\|\theta(t)\|_{L_3(\Omega)}^3+c(1/\delta_1)
\|\theta\|_{L_2(0,t;H^1(\Omega))}^{3/2}A_1^{3/2}(2).
$$
The second term on the right-hand side of (4.37) is estimated by
$$\eqal{
&\intop_0^t\|\theta\|_{L_6(\Omega)}\|\epsi_{t'}\|_{{\bm L}_3(\Omega)}
\|\epsi_{t'}\|_{{\bm L}_2(\Omega)}dt'\cr
&\le\|\theta\|_{L_{6,2}(\Omega^t)}\|\epsi_{t'}\|_{{\bm L}_{3,s_1}(\Omega^t)}
\|\epsi_{t'}\|_{{\bm L}_{2,\sigma_1}(\Omega^t)}\cr
&\le\delta_2\|\theta\|_{L_2(0,t;\H^1(\Omega))}^2+c(1/\delta_2)
\|\epsi_{t'}\|_{{\bm L}_{3,s_1}(\Omega^t)}^2A_1^2(\sigma_1),\cr}
$$
where (4.22)was used with $\sigma=\sigma_1$, 
${1\over s_1}+{1\over\sigma_1}={1\over2}$, $s_1>2$ but close to 2, 
because $\sigma_1$ can be an arbitrary positive finite number.

\noindent
Now we examine the integral
$$\eqal{
&\bigg(\intop_0^t\bigg|\intop_\Omega|\epsi_{t'}|^3dx\bigg|^{s_1/3}dt'
\bigg)^{1/s_1}\le\bigg(\intop_0^t\|\epsi_{t'}(t')\|_{{\bm L}_4(\Omega)}^{2s_1/3}
\|\epsi_{t'}(t')\|_{{\bm L}_2(\Omega)}^{s_1/3}dt'\bigg)^{1/s_1}\cr
&\le\bigg(\intop_0^t\|\epsi_{t'}(t')\|_{{\bm L}_4(\Omega)}^{2s_1\lambda_1/3}dt'
\bigg)^{1/s_1\lambda_1}\bigg(\intop_0^t
\|\epsi_{t'}(t')\|_{{\bm L}_2(\Omega)}^{s_1\lambda_2/3}dt'
\bigg)^{1/s_1\lambda_2}\equiv I_2,\cr}
$$
where $1/\lambda_1+1/\lambda_2=1$.
Setting $s_1\lambda_1=3$, $s_1\lambda_2={3s_1\over3-s_1}$, we obtain
$$
I_2=\bigg(\intop_0^t\|\epsi_{t'}(t')\|_{{\bm L}_4(\Omega)}^2dt'\bigg)^{1/3}
\bigg[\bigg(\intop_0^t\|\epsi_{t'}(t')\|_{{\bm L}_2(\Omega)}^{s_1\over3-s_1}dt'
\bigg)^{3-s_1\over s_1}\bigg]^{1/3}.
$$
Now, using (4.22) with $\sigma={s_1\over3-s_1}$ gives
$$
I_2\le c\bigg(\intop_0^t\|\epsi_{t'}(t')\|_{\H^1(\Omega)}^2dt'\bigg)^{1/3}
A_1^{1/3}\bigg({s_1\over3-s_1}\bigg)
$$
for any $s_1\in(2,3)$.

\noindent
Finally, the third term on the right-hand side of (4.37) is estimated by
$$
c\|\theta\|_{L_{2,\infty}(\Omega^t)}\|g\|_{L_{2,1}(\Omega^t)}\le cA_0
\|g\|_{L_{2,1}(\Omega^t)}.
$$
Inserting the above estimates into (4.37) and assuming that $\delta_1$, 
$\delta_2$ are sufficiently small, we arrive at
$$\eqal{
&\|\theta(t)\|_{L_3(\Omega)}^3+\|\theta\|_{L_2(0,t;H^1(\Omega))}^2
\le\|\theta\|_{L_2(\Omega^t)}^2+cA_1^{3/2}(2)
\|\theta\|_{L_2(0,t;H^1(\Omega))}^{3/2}\cr
&\quad+cA_1^2(\sigma_1)A_1^{2/3}\bigg({s_1\over3-s_1}\bigg)
\|\epsi_{t'}\|_{L_2(0,t;\H^1(\Omega))}^{4/3}+cA_0\|g\|_{L_{2,1}(\Omega^t)}
+\|\theta_0\|_{L_3(\Omega)}^3,\cr}
\leqno(4.38)
$$
where ${1\over s_1}+{1\over\sigma_1}={1\over2}$, $2<s_1<3$.

\noindent
Let us choose $s_1={5\over2}$. Then ${s_1\over3-s_1}=5$ and $\sigma_1=10$. 
Since $A_1(\sigma)$ is an increasing function of $\sigma$ it follows from 
(4.38) that
$$\eqal{
&\|\theta(t)\|_{L_3(\Omega)}^3+\|\theta\|_{L_2(0,t;H^1(\Omega))}^2\le
\|\theta\|_{L_2(\Omega^t)}^2\cr
&\quad+c(A_1^6(10)+A_1^{8/3}(10)\|\epsi_{t'}\|_{L_2(0,t;\H^1(\Omega))}^{4/3}
+A_0\|g\|_{L_{2,1}(\Omega^t)})
+\|\theta_0\|_{L_3(\Omega)}^3,\cr}
\leqno(4.39)
$$
where we used the Young inequality in the second term on the right-hand side 
of (4.38).

By virtue of (4.34) and (4.10) we obtain from (4.39) the inequality
$$\eqal{
&\|\epsi_{t'}\|_{\V_2(\Omega^t)}\le c(t)(A_0+A_1^3(10)+A_1^4(10)+A_0^{1/2}
\|g\|_{L_{2,1}(\Omega^t)}^{1/2}\cr
&\quad+\|\theta_0\|_{L_3(\Omega)}^{3/2}+\|\u_0\|_{\H^2(\Omega)}+
\|\u_1\|_{\H^1(\Omega)}+\|\b\|_{{\bm L}_2(\Omega^t)}).\cr}
\leqno(4.40)
$$
Employing (4.40) in (4.39) yields (4.35). This completes the proof.
\end{proof}

Using (4.35) in (4.33) implies
$$\eqal{
&\|\u_{t'}\|_{L_\infty(0,t;\H^1(\Omega))}+\|\u_{t'}\|_{L_2(0,t;\H^2(\Omega))}
+\|\u\|_{L_\infty(0,t;\H^2(\Omega))}\cr
&\le c(t)(A_2+\|\u_0\|_{\H^2(\Omega)}+\|\u_1\|_{\H^1(\Omega)}+
\|\b\|_{{\bm L}_2(\Omega^t)})
\le c(t)A_3,\cr}
\leqno(4.41)
$$
with $A_3$ defined in (4.36).

>From what it has already been proved we deduce

\begin{lemma}\label{4.8.}
Assume that $\u_0\in\H^2(\Omega)$, 
$\u_1\in\B_{2,10}^{2-1/5}(\Omega)\cap\B_{2,\sigma}^{2-2/\sigma}(\Omega)$, 
$\theta_0\in H^1(\Omega)$, 
$\b\in{\bm L}_{2,10}(\Omega^t)\cap{\bm L}_{2,\sigma}(\Omega^t)$, 
$g\in L_2(\Omega^t)$, $\sigma>4$.\\
Then
$$
\|\theta_t\|_{L_2(\Omega^t)}+\|\nabla\theta\|_{L_\infty(0,t;{\bm L}_2(\Omega))}
\le\varphi(A_4(\sigma)),\quad \sigma>4,
\leqno(4.42)
$$
where
$$\eqal{
A_4(\sigma)&=\|\u_0\|_{\H^2(\Omega)}+\|\u_1\|_{\B_{2,10}^{2-1/5}(\Omega)}+
\|\u_1\|_{\B_{2,\sigma}^{2-2/\sigma}(\Omega)}\cr
&\quad+\|\theta_0\|_{H^1(\Omega)}+\|\b\|_{{\bm L}_{2,10}(\Omega^t)}+
\|\b\|_{{\bm L}_{2,\sigma}(\Omega^t)}+\|g\|_{L_2(\Omega^t)}.\cr}
\leqno(4.43)
$$
\end{lemma}

\begin{proof} 
Multiplying (1.2) by $\theta_t$ and integrating over $\Omega$ gives
$$\eqal{
&c_v\intop_\Omega\theta\theta_t^2dx+{k\over2}{d\over dt}\intop_\Omega
|\nabla\theta|^2dx\cr
&\le c\intop_\Omega\theta|\epsi_t|\,|\theta_t|dx+c\intop_\Omega
|\epsi_t|^2|\theta_t|dx+\intop_\Omega|g|\,|\theta_t|dx.\cr}
\leqno(4.44)
$$
Applying the H\"older and the Young inequalities the first term on the 
right-hand side of (4.44) is estimated by
$$\eqal{
&\bigg(\intop_\Omega\theta\theta_t^2dx\bigg)^{1/2}\bigg(\intop_\Omega
\theta|\epsi_t|^2dx\bigg)^{1/2}
\le\delta_1\intop_\Omega\theta\theta_t^2dx+c(1/\delta_1)\intop_\Omega
\theta|\epsi_t|^2dx,\cr}
$$
where on account of (4.35) the second integral is bounded by
$$
c\bigg(\intop_\Omega\theta^3dx\bigg)^{1/3}\bigg(\intop_\Omega
|\epsi_t|^3dx\bigg)^{2/3}\le cA_2\|\epsi_t\|_{{\bm L}_3(\Omega)}^2.
$$
The second term on the right-hand side of (4.44) can be estimated by
$$\eqal{
&c\|\theta_t\|_{L_2(\Omega)}\|\epsi_t\|_{{\bm L}_6(\Omega)}
\|\epsi_t\|_{{\bm L}_3(\Omega)}
\le\delta_2\|\theta_t\|_{L_2(\Omega)}^2+c(1/\delta_2)
\|\epsi_t\|_{{\bm L}_6(\Omega)}^2\|\epsi_t\|_{{\bm L}_3(\Omega)}^2.\cr}
$$
Finally, the third term on the right-hand side of (4.44) is bounded by
$$
\delta_3\|\theta_t\|_{L_2(\Omega)}^2+c(1/\delta_3)\|g\|_{L_2(\Omega)}^2.
$$
Employing the above estimates in (4.44), assuming $\delta_i$, $i=1,2,3$, 
sufficiently small, recalling that $\theta\ge\theta_*>0$, and integrating 
the result with respect to time, we conclude on account of (4.41) that
$$\eqal{
&\|\theta_t\|_{L_2(\Omega^t)}^2+\|\nabla\theta(t)\|_{{\bm L}_2(\Omega)}^2\le 
cA_2A_3^2\cr
&\quad+cA_3^2\|\epsi_{t'}\|_{L_\infty(0,t;{\bm L}_3(\Omega))}^2+c
\|g\|_{L_2(\Omega^t)}^2+c\|\nabla\theta_0\|_{{\bm L}_2(\Omega)}^2.\cr}
\leqno(4.45)
$$

In view (4.41) and (4.45), applying Corollary 4.5 to problem (4.19) we 
conclude that
$$\eqal{
&\|\u_{t'}\|_{\W_{2,\sigma}^{2,1}(\Omega^t)}\le c(A_3+A_2^{1/2}A_3+A_3\|\epsi_{t'}\|_{L_\infty(0,t;{\bm L}_3(\Omega))}\cr
&\quad+\|g\|_{L_2(\Omega^t)}+\|\b\|_{{\bm L}_{2,\sigma}(\Omega^t)}+\|\theta_0\|_{H^1(\Omega)}+
\|\u_1\|_{\B_{2,\sigma}^{2-2/\sigma}(\Omega)})\cr
&\le cA_3\|\epsi_{t'}\|_{L_\infty(0,t;{\bm L}_3(\Omega))}+c\varphi
(A_4(\sigma)),\cr}
\leqno(4.46)
$$
where $\sigma$ is an arbitrary finite number.

\noindent
By the definition of the tensor $\epsi$, inequality (4.46) implies
$$
\|\epsi_{t'}\|_{\W_{2,\sigma}^{1,1/2}(\Omega^t)}\le cA_3
\|\epsi_{t'}\|_{L_\infty(0,t;{\bm L}_3(\Omega))}+cA_4(\sigma).
\leqno(4.47)
$$
In view of the interpolation inequality
$$
\|\epsi_{t'}\|_{L_\infty(0,t;{\bm L}_3(\Omega))}\le\delta
\|\epsi_{t'}\|_{\W_{2,\sigma}^{1,1/2}(\Omega^t)}
+c(1/\delta)\|\epsi_{t'}\|_{{\bm L}_{2,\sigma}(\Omega^t)},
\leqno(4.48)
$$
which holds for $\sigma>4$, it follows from (4.47) and (4.35) that
$$
\|\epsi_{t'}\|_{\W_{2,\sigma}^{1,1/2}(\Omega^t)}\le\varphi(A_3,A_4(\sigma)),
\quad\sigma>4.
\leqno(4.49)
$$
Hence,
$$
\|\epsi_{t'}\|_{L_\infty(0,t;{\bm L}_3(\Omega))}\le\varphi(A_3,A_4(\sigma)),
\quad\sigma>4.
\leqno(4.50)
$$
Applying (4.50) in (4.45) gives (4.42). This completes the proof.
\end{proof}

Let us note that since $t$ is finite, estimate (4.49) in conjunction with the 
H\"older inequality implies that
$$
\|\epsi_{t'}\|_{\W_{2,\sigma_0}^{1,1/2}(\Omega^t)}\le\varphi
(t,A_3,A_4(\sigma)),
\leqno(4.51)
$$
where $\sigma>4$ and $\sigma_0\ge1$.

\noindent
Similarly, by (4.46) and (4.50), it follows that
$$
\|\u_{t'}\|_{\W_{2,\sigma_0}^{2,1}(\Omega^t)}\le\varphi(t,A_3,A_4(\sigma)),
\leqno(4.52)
$$
where $\sigma>4$ and $\sigma_0\ge1$.

\begin{lemma}\label{4.9.} 
Assume that $\theta\in L_{p,\sigma}(\Omega^t)$, 
$\b\in{\bm L}_{p,\sigma}(\Omega^t)$, $\u_0\in\W_p^1(\Omega)$, 
$\u_1\in\B_{p,\sigma}^{2-2/\sigma}(\Omega)$, $p\in(1,\infty)$, 
$\sigma\in(1,\infty)$.\\
Then
$$
\|\u_{t'}\|_{\W_{p,\sigma}^{1,1/2}(\Omega^t)}\le 
c(t)[\|\theta\|_{L_{p,\sigma}(\Omega^t)}+A_5(p,\sigma)],
\leqno(4.53)
$$
where
$$
A_5(p,\sigma)=\|\u_0\|_{\W_p^1(\Omega)}+
\|\u_1\|_{\B_{p,\sigma}^{2-2/\sigma}(\Omega)}+
\|\b\|_{{\bm L}_{p,\sigma}(\Omega^t)},
$$
and the constant $c(t)$ depends exponentially on $t$.
\end{lemma}

\begin{proof} 
Let us consider system (4.19) and apply the inequality (4.20). Representing 
$\epsi$ by
$$
\epsi(t)=\intop_0^t\epsi_{t'}(t')dt'+\epsi(0),
\leqno(4.54)
$$
and using the generalized Minkowski inequality, we obtain
$$\eqal{
&\|\epsi\|_{{\bm L}_{p,\sigma}(\Omega^t)}^\sigma\le c\intop_0^t
\bigg\|\intop_0^{t'}\epsi_{t''}(t'')dt''\bigg\|_{{\bm L}_p(\Omega)}^\sigma dt'
+c\intop_0^t\|\epsi(\u_0)\|_{{\bm L}_p(\Omega)}^\sigma dt'\cr
&\le c\intop_0^t\bigg(\intop_0^{t'}\|\epsi_{t''}(t'')\|_{{\bm L}_p(\Omega)}
dt''\bigg)^\sigma dt'+ct\|\epsi(\u_0)\|_{{\bm L}_p(\Omega)}^\sigma\cr
&\le c\intop_0^t\bigg(\intop_0^{t'}\|\epsi_{t''}(t'')\|_{{\bm L}_p(\Omega)}^\sigma
dt''\bigg)(t')^{\sigma-1}dt'+ct\|\epsi(\u_0)\|_{{\bm L}_p(\Omega)}^\sigma.\cr}
$$
Consequently, denoting 
$$
a(t)=\|\epsi_{t'}\|_{{\bm L}_{p,\sigma}(\Omega^t)}^\sigma=\intop_0^t
\|\epsi_{t'}(t')\|_{{\bm L}_p(\Omega)}^\sigma dt',
$$
we deduce from (4.20) the inequality
$$
a(t)\le\intop_0^t\alpha(t')a(t')dt'+A(t),
$$
where
$$\eqal{
&\alpha(t)=ct^{\sigma-1},\cr
&A(t)=c(\|\theta\|_{L_{p,\sigma}(\Omega^t)}^\sigma+
\|\b\|_{{\bm L}_{p,\sigma}(\Omega^t)}^\sigma+
t\|\epsi(\u_0)\|_{{\bm L}_p(\Omega)}^\sigma
+\|\u_1\|_{\B_{p,\sigma}^{2-2/\sigma}(\Omega)}^\sigma).\cr}
$$
Hence, by the Gronwall inequality, it follows that
$$\eqal{
a(t)&\le A(t)+\intop_0^t\alpha(t')A(t')e^{\intop_{t'}^t\alpha(t'')dt''}dt'\cr
&\le A(t)(1+\alpha_1(t)e^{\alpha_1(t)}),\quad 
\alpha_1(t)=t\alpha(t).\cr}
$$
Thus,
$$\eqal{
\|\epsi_{t'}\|_{{\bm L}_{p,\sigma}(\Omega^t)}^\sigma
&\le c(t)(\|\theta\|_{L_{p,\sigma}(\Omega^t)}^\sigma+t
\|\epsi(\u_0)\|_{{\bm L}_p(\Omega)}^\sigma\cr
&\quad+\|\u_1\|_{\B_{p,\sigma}^{2-2/\sigma}(\Omega)}^\sigma+
\|\b\|_{{\bm L}_{p,\sigma}(\Omega^t)}^\sigma)\cr
&\equiv c(t)(\|\theta\|_{L_{p,\sigma}(\Omega^t)}^\sigma+D^\sigma(t))\cr}
\leqno(4.55)
$$
with $c(t)=c(1+\alpha_1(t)e^{\alpha_1(t)})$.

\noindent
Using (4.55) in (4.54) yields the analogous bound on 
$\|\epsi\|_{{\bm L}_{p,\sigma}(\Omega^t)}^\sigma$. Consequently, on account of 
(4.20) the corresponding bound on 
$\|\u_{t'}\|_{\W_{p,\sigma}^{1,1/2}(\Omega^t)}$ follows as well. The proof 
is completed.
\end{proof}

On the basis of Lemma 4.9 we prove now

\begin{lemma}\label{4.10.} 
Assume that $\u_0\in\W_r^1(\Omega)$, $\u_1\in\B_{r,r}^{2-2/r}(\Omega)$, 
$\theta_0\in L_r(\Omega)$, $\b\in{\bm L}_r(\Omega^T)$, $g\in L_r(\Omega^T)$, 
$r\in(1,\infty)$.\\
Then
$$
\|\theta\|_{L_{r,\infty}(\Omega^t)}\le A_6(r,r,t),\quad r<\infty,
\leqno(4.56)
$$
where
$$
A_6(r,r,t)=c(t)(\|\theta_0\|_{L_r(\Omega)}+\root r \of{r}A_5(r,r)+
\|g\|_{L_r(\Omega^t)}+1),
\leqno(4.57)
$$
and $A_5(p,\sigma)$ is defined in (4.53).
\end{lemma}

\begin{proof} 
Multiplying (1.2) by $\theta^r$, $r>1$, and integrating over $\Omega$ gives
$$\eqal{
&{c_v\over r+2}{d\over dt}\|\theta\|_{L_{r+2}(\Omega)}^{r+2}+
{4kr\over(r+1)^2}\intop_\Omega\left|\nabla\theta^{r+2\over2}\right|^2dx\cr
&=-\intop_\Omega\theta^{r+1}(\A_2\alphaa)\cdot\epsi_tdx+\intop_\Omega
\theta^r(\A_1\epsi_t)\cdot\epsi_tdx+\intop_\Omega\theta^rgdx.\cr}
\leqno(4.58)
$$
Hence, after integrating with respect to time,
$$\eqal{
&{c_v\over r+2}\|\theta(t)\|_{L_{r+2}(\Omega)}^{r+2}+
{4kr\over(r+1)^2}\intop_{\Omega^t}\left|\nabla\theta^{r+2\over2}\right|^2dxdt'
\cr
&\le c\intop_{\Omega^t}\theta^{r+1}|\epsi_{t'}|dxdt'+c\intop_{\Omega^t}
\theta^r|\epsi_{t'}|^2dxdt'+\intop_{\Omega^t}\theta^rgdxdt'
+{c_v\over r+2}\|\theta_0\|_{L_{r+2}(\Omega)}^{r+2}.\cr}
\leqno(4.59)
$$
Here let us recall Lemma 4.9 which provides the inequality
$$
\|\epsi_{t'}\|_{{\bm L}_r(\Omega^t)}\le c(t)(\|\theta\|_{L_r(\Omega^t)}+
A_5(r,r))
\leqno(4.60)
$$
for $r\in(1,\infty)$.

On account of (4.60) the second integral on the right-hand side of (4.59) 
is estimated as follows
$$\eqal{
&\intop_{\Omega^t}\theta^r|\epsi_{t'}|^2dxdt'\le
\|\theta\|_{L_{r+2}(\Omega^t)}^r\|\epsi_{t'}\|_{{\bm L}_{r+2}(\Omega^t)}^2\cr
&\le c(t)\|\theta\|_{L_{r+2}(\Omega^t)}^r
(\|\theta\|_{L_{r+2}(\Omega^t)}^2+A_5^2(r+2,r+2))\cr
&\le c(t)(\|\theta\|_{L_{r+2}(\Omega^t)}^{r+2}+A_5^{r+2}(r+2,r+2)).\cr}
\leqno(4.61)
$$
Now, using (4.61) we estimate the first integral on the right-hand side of 
(4.59) by
$$\eqal{
&\intop_{\Omega^t}\theta^{r+1}|\epsi_{t'}|dxdt'=\intop_{\Omega^t}
\theta^{{r\over2}+1}\theta^{r\over2}|\epsi_{t'}|dxdt'\cr
&\le\left\|\theta^{{r\over2}+1}\right\|_{L_2(\Omega^t)}
\left\|\theta^{r\over2}|\epsi_{t'}|\right\|_{L_2(\Omega^t)}
=\|\theta\|_{L_{r+2}(\Omega^t)}^{r+2\over2}
\|\theta^r|\epsi_{t'}|^2\|_{L_1(\Omega^t)}^{1/2}\cr
&\le c(t)\|\theta\|_{L_{r+2}(\Omega^t)}^{r+2\over2}
(\|\theta\|_{L_{r+2}(\Omega^t)}^{r+2\over2}+A_5^{r+2\over2}(r+2,r+2))\cr
&\le c(t)(\|\theta\|_{L_{r+2}(\Omega^t)}^{r+2}+A_5^{r+2}(r+2,r+2)).\cr}
$$
Finally, with the use of H\"older's and Young's inequalities the third 
integral on the right-hand side of (4.59) is estimated by
$$\eqal{
&\intop_{\Omega^t}\theta^rgdxdt'\le\|\theta\|_{L_{r+2}(\Omega^t)}^r
\|g\|_{L_{(r+2)/2}(\Omega^t)}
\le c\|\theta\|_{L_{r+2}(\Omega^t)}^{r+2}+c(\|g\|_{L_{r+2}(\Omega^t)}^{r+2}
+1).\cr}
$$
Inserting the above estimates into (4.59) leads to
$$\eqal{
\|\theta(t)\|_{L_{r+2}(\Omega)}^{r+2}&\le{r+2\over c_v}c(t)\bigg(\intop_0^t
\|\theta(t')\|_{L_{r+2}(\Omega)}^{r+2}dt'+A_5^{r+2}(r+2,r+2)\cr
&\quad+\|g\|_{L_{r+2}(\Omega^t)}^{r+2}+1\bigg)+
\|\theta_0\|_{L_{r+2}(\Omega)}^{r+2}.\cr}
$$
Hence, by the Gronwall inequality, we conclude that
$$\eqal{
\|\theta(t)\|_{L_{r+2}(\Omega)}^{r+2}
&\le[\|\theta_0\|_{L_{r+2}(\Omega)}^{r+2}+(r+2)c(t)(A_5^{r+2}(r+2,r+2)\cr
&\quad+\|g\|_{L_{r+2}(\Omega^t)}^{r+2}+1)]\exp(c(t)(r+2))\cr
&\le A_6^{r+2}(r+2,r+2,t)\cr}
$$
for $t\in(0,T)$, with $A_6$ defined by (4.57). This gives (4.56). The proof 
is completed.
\end{proof}

\begin{corollary}\label{4.11.} 
Taking into account that $\root r \of{r}$ is bounded let us define
$$\eqal{
A_7(r,r,t)&=c(t)(\|\u_0\|_{\W_r^1(\Omega)}+
\|\u_1\|_{\B_{r,r}^{2-2/r}(\Omega)}+\|\theta_0\|_{L_r(\Omega)}\cr
&\quad+\|\b\|_{{\bm L}_{r,r}(\Omega^t)}+\|g\|_{L_r(\Omega^t)}).\cr}
\leqno(4.62)
$$
Then, according to (4.56) and the definition of $A_5(r,r)$ in Lemma 4.9,
$$
\|\theta\|_{L_{r,\infty}(\Omega^t)}\le A_7(r,r,t),\quad r\in(1,\infty),\quad 
t\le T.
\leqno(4.63)
$$
Moreover, by (4.53) and (4.63),
$$\eqal{
\|\epsi_{t'}\|_{{\bm L}_{r,\sigma}(\Omega^t)}
&\le c\|\u_{t'}\|_{\W_{r,\sigma}^{1,1/2}(\Omega^t)}\cr
&\le c(t)A_7(r,r,t),\quad (r,\sigma)\in(1,\infty).\cr}
\leqno(4.64)
$$
\end{corollary}

Let us consider now the elliptic problem resulting from (1.2) and $(1.3)_2$:
$$\eqal{
&-k\Delta\theta=-c_v\theta\theta_t-\theta(\A_2\alphaa)\cdot\epsi_t+
(\A_1\epsi_t)\cdot\epsi_t+g\quad &{\rm in}\ \ \Omega,\cr
&\n\cdot\nabla\theta=0\quad &{\rm on}\ \ S.\cr}
\leqno(4.65)
$$
We have

\begin{lemma}\label{4.12.}
Asume that $\u_0\in\H^2(\Omega)\cap\W_{2s\over2-s}^1(\Omega)\cap
\B_{2,\sigma}^{2-2/\sigma}(\Omega)$, \\
$\u_1\in\B_{2,10}^{9/5}(\Omega)\cap
\B_{{2s\over2-s},{2s\over2-s}}^{3s-2\over s}(\Omega)$, 
$\theta_0\in H^1(\Omega)\cap L_{2s\over2-s}(\Omega)$, \\
$\b\in{\bm L}_{2,10}(\Omega^T)\cap{\bm L}_{2,\sigma}(\Omega^T)\cap
{\bm L}_{2s\over2-s}(\Omega^T)$, 
$g\in L_2(\Omega^T)\cap L_{2s\over2-s}(\Omega^T)$.\\
Then
$$
\|\theta\|_{W_s^{2,1}(\Omega^t)}\le\varphi(A_8(s,\sigma,t)),\quad t\le T,
\leqno(4.66)
$$
where
$$\eqal{
&A_8(s,\sigma,t)=c(t)(\|\u_0\|_{\H^2(\Omega)}+
\|\u_0\|_{\W_{2s\over2-s}^1(\Omega)}\cr
&\quad+\|\u_1\|_{\B_{2,10}^{9/5}(\Omega)}+
\|\u_1\|_{\B_{2,\sigma}^{2-2/\sigma}(\Omega)}+
\|\u_1\|_{\B_{{2s\over2-s},{2s\over2-s}}^{3s-2\over s}(\Omega)}\cr
&\quad+\|\theta_0\|_{H^1(\Omega)}+\|\theta_0\|_{L_{2s\over2-s}(\Omega)}+
\|\b\|_{{\bm L}_{2,10}(\Omega^t)}+\|\b\|_{{\bm L}_{2,\sigma}(\Omega^t)}\cr
&\quad+\|\b\|_{{\bm L}_{2s\over2-s}(\Omega^t)}+\|g\|_{L_2(\Omega^t)}+
\|g\|_{L_{2s\over2-s}(\Omega^t)}),\cr}
\leqno(4.67)
$$
and $1<s<2$, $\sigma>4$.
\end{lemma}

\begin{proof} 
On account of (4.42), (4.63) and (4.64) it follows from $(4.65)_1$ that
$$\eqal{
\|\Delta\theta\|_{L_s(\Omega^t)}&\le c(\|\theta\|_{L_{2s\over2-s}(\Omega^t)}
\|\theta_t\|_{L_2(\Omega^t)}\cr
&\quad+\|\theta\|_{L_{2s}(\Omega^t)}\|\epsi_{t'}\|_{{\bm L}_{2s}(\Omega^t)}+
\|\epsi_{t'}\|_{{\bm L}_{2s}(\Omega^t)}^2+\|g\|_{L_s(\Omega^t)})\cr
&\le c(t)\bigg(A_7\bigg({2s\over2-s},{2s\over2-s},t\bigg)\varphi
(A_4(\sigma))\cr
&\quad+A_7^2(2s,2s,t)+\|g\|_{L_s(\Omega^t)}\bigg),\cr}
\leqno(4.68)
$$
where $s<2$ is close to 2, $\sigma>4$ and $t\le T$.

\noindent
In view of the estimate
$$
\|\theta\|_{W_s^2(\Omega)}\le c(\|\Delta\theta\|_{L_s(\Omega)}+
\|\theta\|_{L_s(\Omega)}),
\leqno(4.69)
$$
which holds for the Neumann problem (4.65), we deduce from (4.68) that
$$\eqal{
&\|\theta\|_{L_s(0,t;W_s^2(\Omega))}\le c\|\theta\|_{L_s(\Omega^t)}+
c(t)\bigg(A_7\bigg({2s\over2-s},{2s\over2-s},t\bigg)\varphi(A_4(\sigma))\cr
&\quad+A_7^2(2s,2s,t)+\|g\|_{L_s(\Omega^t)}\bigg),\cr}
\leqno(4.70)
$$
where $\sigma>4$, $s<2$ close to 2, $t\le T$.

\noindent
Now, combining (4.70) with (4.42) and using (4.10) we conclude that
$$\eqal{
&\|\theta\|_{W_s^{2,1}(\Omega^t)}\le c(t)(A_0+\varphi(A_4(\sigma)))\cr
&\quad+c(t)\bigg(A_7\bigg({2s\over2-s},{2s\over2-s},t\bigg)
\varphi(A_4(\sigma))\cr
&\quad+A_7^2(2s,2s,t)+\|g\|_{L_s(\Omega^t)}\bigg)\le\varphi
(A_8(s,\sigma,t)),\cr}
\leqno(4.71)
$$
where the latter inequality follows from the definitions of the quantities 
$A_0$ (see (4.10)), $A_4(\sigma)$ (see (4.43)), $A_7(r,r,t)$ (see (4.62)) 
and $A_8(s,\sigma,t)$ (see (4.67)). This proves the lemma.
\end{proof}

\begin{corollary}\label{4.13.} 
In view of the imbedding
$$
\nabla W_s^{2,1}(\Omega^T)\subset{\bm L}_{p',s}(\Omega^T),\quad s\in(1,2),
\leqno(4.72)
$$
where $s<p'\le{3\over3/s-1}$, and the H\"older inequality with respect to 
the integral over $\Omega$ in the case $1<p'<s$, we have
$$
\|\nabla v\|_{{\bm L}_{p',s}(\Omega^T)}\le c\|v\|_{W_s^{2,1}(\Omega^T)}
\leqno(4.73)
$$
for any $v\in W_s^{2,1}(\Omega^T)$, $1<p'\le{3\over3/s-1}$.\\
Applying the imbedding inequality (4.73) to (4.66) we deduce that
$$
\|\nabla\theta\|_{{\bm L}_{p',s}(\Omega^t)}\le c(t)\varphi(A_8(s,\sigma,t))
\leqno(4.74)
$$
where $1<p'\le{3\over3/s-1}$, $s\in(1,2)$, $\sigma>4$, and $A_8(s,\sigma,t)$ 
is defined in (4.67).
\end{corollary}

We shall use (4.74) to get more regularity estimates on $\epsi$. To this 
purpose we return to the viscoelasticity system (4.19) and prove the following 
result analogous to Lemma 4.9.

\begin{lemma}\label{4.14.} 
Assume that $\nabla\theta\in{\bm L}_{p,\sigma}(\Omega^T)$, 
$\b\in{\bm L}_{p,\sigma}(\Omega^T)$,\\ 
$\u_1\in\B_{p,\sigma}^{2-2/\sigma}(\Omega)$, $p,\sigma\in(1,\infty)$.\\
Then
$$\eqal{
&\|\epsi_{t'}\|_{\W_{p,\sigma}^{1,1/2}(\Omega^t)}\le c
\|\u_{t'}\|_{\W_{p,\sigma}^{2,1}(\Omega^t)}\cr
&\le c(t)(\|\nabla\theta\|_{{\bm L}_{p,\sigma}(\Omega^t)}+
\|\b\|_{{\bm L}_{p,\sigma}(\Omega^t)}+
\|\u_1\|_{\B_{p,\sigma}^{2-2/\sigma}(\Omega)})\cr
&\equiv c(t)(\|\nabla\theta\|_{{\bm L}_{p,\sigma}(\Omega^t)}+A_9(p,\sigma)),
\cr}
\leqno(4.75)
$$
where $t\le T$, $p,\sigma\in(1,\infty)$, and
$$
A_9(p,\sigma)=\|\u_1\|_{\B_{p,\sigma}^{2-2/\sigma}(\Omega)}+
\|\b\|_{{\bm L}_{p,\sigma}(\Omega^t)}.
$$
\end{lemma}

\begin{proof} 
Let us consider system (4.19). By inequality (4.21) we have
$$\eqal{
&\|\nabla\epsi_{t'}\|_{{\bm L}_{p,\sigma}(\Omega^t)}\le
\|\u_{t'}\|_{\W_{p,\sigma}^{2,1}(\Omega^t)}\le c
(\|\nabla\epsi\|_{{\bm L}_{p,\sigma}(\Omega^t)}+
\|\nabla\theta\|_{{\bm L}_{p,\sigma}(\Omega^t)}\cr
&\quad+\|\u_1\|_{\B_{p,\sigma}^{2-2/\sigma}(\Omega)}+
\|\b\|_{{\bm L}_{p,\sigma}(\Omega^t)})\cr}
\leqno(4.76)
$$
for $t\le T$. We use now the formula
$$
\nabla\epsi(t)=\intop_0^t\nabla\epsi_{t'}(t')dt'+\nabla\epsi(0),\quad
\nabla\epsi(0)=\nabla\epsi(\u_0),
\leqno(4.77)
$$
and repeat the proof of Lemma 4.9 with $\nabla\epsi$, $\nabla\theta$ in 
place of $\epsi$, $\theta$. In result we conclude that
$$
\|\nabla\epsi_{t'}\|_{{\bm L}_{p,\sigma}(\Omega^t)}\le c(t)
(\|\nabla\theta\|_{{\bm L}_{p,\sigma}(\Omega^t)}+
\|\u_1\|_{\B_{p,\sigma}^{2-2/\sigma}(\Omega)}+
\|\b\|_{{\bm L}_{p,\sigma}(\Omega^t)}).
\leqno(4.78)
$$
Using (4.78) in (4.77) implies the analogous bound on 
$\|\nabla\epsi\|_{{\bm L}_{p,\sigma}(\Omega^t)}$ and then, by (4.76), on 
$\|\u_{t'}\|_{\W_{p,\sigma}^{2,1}(\Omega^t)}$ as well.
\end{proof}

\begin{corollary}\label{4.15.} 
Applying estimate (4.74) in (4.75) yields
$$
\|\epsi_{t'}\|_{\W_{p',s}^{1,1/2}(\Omega^t)}\le c(t)(\varphi(A_8(s,\sigma,t))+
A_9(p',s))
\leqno(4.79)
$$
for $1<p'\le{3\over3/s-1}$, $s\in(1,2)$, $\sigma>4$, $t\le T$.
\end{corollary}

\begin{corollary}\label{4.16.} 
Let us consider the imbedding
$$
\W_{p',s}^{1,1/2}(\Omega^T)\subset{\bm L}_{\infty,2}(\Omega^T),\quad s\in(1,2),
\leqno(4.80)
$$
which holds true provided $p'>{3\over2-{2\over s}}$. 
This condition together with\\ 
$p'\le{3\over{3\over s}-1}$ implies that
$$
{5\over3}<s<2.
$$
Consequently, it follows from (4.79) that
$$
\|\epsi_{t'}\|_{{\bm L}_{\infty,2}(\Omega^t)}\le c(t)(\varphi(A_8(s,\sigma,t))+
A_9(p',s))\equiv A_{10}(s,\sigma,p',t)
\leqno(4.81)
$$
for ${3\over2-{2\over s}}<p'\le{3\over{3\over s}-1}$, 
$s\in\left({5\over3},2\right)$, $\sigma>4$.
\end{corollary}

Estimate (4.81) plays the key role in getting $L_\infty(\Omega^T)$-norm 
bound for~$\theta$.

\begin{lemma}\label{4.17.} 
Let the assumptions of Lemma 4.12 be satisfied. Moreover, let 
$\u_1\in\B_{p',s}^{2-2/s}(\Omega)$, $\theta_0\in L_\infty(\Omega)$, 
$\b\in{\bm L}_{p',s}(\Omega^T)$, $g\in L_{\infty,1}(\Omega^T)$, where
$$
s\in\bigg({5\over3},2\bigg),\quad 
{3\over2-{2\over s}}<p'\le{3\over{3\over s}-1},\quad \sigma>4.
$$
Then
$$\eqal{
\lc{\|\theta\|_{L_\infty(\Omega^t)}\le c(t)[A_{10}(s,\sigma,p',t)+
A_{10}^2(s,\sigma,p',t)}
\rc{\quad+\|\theta_0\|_{L_\infty(\Omega)}+\|g\|_{L_{\infty,1}(\Omega^t)}]
\equiv c(t)A_{11}(s,\sigma,p',t),}}
\leqno(4.82)
$$
where $A_{10}(s,\sigma,p',t)$ is defined in (4.81).
\end{lemma}

\begin{proof} 
Let us consider once more the identity (4.58). Regarding (4.81), we have
$$\eqal{
{1\over r+2}{d\over dt}\|\theta\|_{L_{r+2}(\Omega)}^{r+2}&\le c
\|\theta\|_{L_{r+1}(\Omega)}^{r+1}\|\epsi_t\|_{{\bm L}_\infty(\Omega)}\cr
&\quad+c\|\theta\|_{L_r(\Omega)}^r\|\epsi_t\|_{{\bm L}_\infty(\Omega)}^2+
\|\theta\|_{L_r(\Omega)}^r\|g\|_{L_\infty(\Omega)}.\cr}
\leqno(4.83)
$$
Taking into account that $\theta\ge\theta_*>0$, and using the inequality
$$
\|\theta\|_{L_{r+1}(\Omega)}^{r+1}\le|\Omega|^{1\over r+2}
\|\theta\|_{L_{r+2}(\Omega)}^{r+1}\le c\|\theta\|_{L_{r+2}(\Omega)}^{r+1},
$$
and
$$\eqal{
\|\theta\|_{L_r(\Omega)}^r
&\le|\Omega|^{2\over r+2}\|\theta\|_{L_{r+2}(\Omega)}^r\le
|\Omega|^{2\over r+2}{1\over\|\theta_*\|_{L_{r+2}(\Omega)}}
\|\theta\|_{L_{r+2}(\Omega)}^{r+1}\cr
&\le{1\over\theta_*}|\Omega|^{1\over r+2}\|\theta\|_{L_{r+2}(\Omega)}^{r+1}
\le c\|\theta\|_{L_{r+2}(\Omega)}^{r+1},\cr}
$$
where constant $c$ is independent of $r$, we deduce from (4.83) that
$$
\|\theta\|_{L_{r+2}(\Omega)}^{r+1}{d\over dt}\|\theta\|_{L_{r+2}(\Omega)}
\le c\|\theta\|_{L_{r+2}(\Omega)}^{r+1}(\|\epsi_t\|_{{\bm L}_\infty(\Omega)}+
\|\epsi_t\|_{{\bm L}_\infty(\Omega)}^2+\|g\|_{L_\infty(\Omega)}).
$$
Hence,
$$
{d\over dt}\|\theta\|_{L_{r+2}(\Omega)}\le c
(\|\epsi_t\|_{{\bm L}_\infty(\Omega)}+\|\epsi_t\|_{{\bm L}_\infty(\Omega)}^2+
\|g\|_{L_\infty(\Omega)}),
$$
which in conjunction with (4.81) leads to
$$\eqal{
&\|\theta(t)\|_{L_{r+2}(\Omega)}\le\|\theta_0\|_{L_{r+2}(\Omega)}+
c(\|\epsi_{t'}\|_{{\bm L}_{\infty,1}(\Omega^t)}+
\|\epsi_{t'}\|_{{\bm L}_{\infty,2}(\Omega^t)}^2\cr
&\quad +\|g\|_{L_{\infty,1}(\Omega^t)})\le c(t)(A_{10}(s,\sigma,p',t)+
A_{10}^2(s,\sigma,p',t)\cr
&\quad+\|\theta_0\|_{L_{r+2}(\Omega)}+\|g\|_{L_{\infty,1}(\Omega^t)}).\cr}
\leqno(4.84)
$$
Now, letting $r\to\infty$ in (4.84) we get the assertion.
\end{proof}

\begin{corollary}\label{4.18.} 
Repeating the proof of Lemma 4.12 with the use of the upper bound (4.82) 
allows us to deduce that
$$
\|\theta\|_{W_2^{2,1}(\Omega^t)}\le\varphi(t,A_{11}(s,\sigma,p',t)),
\leqno(4.85)
$$
where $t\le T$, $s\in\left({5\over3},2\right)$, 
${3\over2-2/s}<p'\le{3\over3/s-1}$, $\sigma>4$.
\end{corollary}

\begin{corollary}\label{4.19.} 
Applying (4.82) in (4.53) yields
$$
\|\epsi_{t'}\|_{{\bm L}_{p,\sigma'}(\Omega^t)}\le c(t)[A_{11}(s,\sigma,p',t)+
A_5(p,\sigma')],
\leqno(4.86)
$$
where $t\le T$, $s\in\left({5\over3},2\right)$, 
${3\over2-2/s}<p'\le{3\over3/s-1}$, 
$\sigma>4$, $p,\sigma'\in(1,\infty)$, $A_5(p,\sigma)$ is defined by (4.53) 
and $A_{11}(s,\sigma,p',t)$ by (4.82).
\end{corollary}

The next step consists in obtaining a "better" estimate for $\theta$ by means 
of the parabolic regularity result stated in Lemma 3.5. To apply this result 
to the quasilinear equation (1.2) we need to prove first that $\theta$ is 
a~continuous function on $\Omega^T$. We are able to prove more, namely the 
H\"older continuity by means of the parabolic De Giorgi method in the same 
way as in [18, Lemma 6.1].

Since the above reference concerns much more general situation, we present 
here for reader's convenience a direct, simpler proof.

Following [14, Chap. II. 7] we record the definition of the space\break
${\cal B}_2(\Omega^T,M,\gamma,r,\delta,\kappa)$, where 
$\Omega^T=\Omega\times(0,T)$, $\Omega\subset\R^n$, $n\in{\N}$, and
$M,\gamma,r,\delta,\break \kappa$ are positive numbers.

\noindent
The function $u\in{\cal B}_2(\Omega^T,M,\gamma,r,\delta,\kappa)$ if:
\vskip6pt

\noindent
(i) $u\in V_2^{1,0}(\Omega^T):=C(0,T;L_2(\Omega))\cap 
L_2(0,T;W_2^1(\Omega))$,

\noindent
(ii) $\ess\sup_{\Omega^T}|u|\le M$,

\noindent
(iii) the function $w(x,t)=\pm u(x,t)$ satisfies the inequalities
$$\eqal{
&\max_{t_0\le t\le t_0+\tau}
\|(w-k)_+\|_{L_2(B_{\varrho-\sigma_1\varrho}(x_0))}^2\le
\|(w-k)_+(\cdot,t_0)\|_{L_2(B_\varrho(x_0))}^2\cr
&\quad+\gamma[(\sigma_1\varrho)^{-2}\|(w-k)_+\|_{L_2(Q(\varrho,\tau))}^2+
\mu^{{2\over r}(1+\kappa)}(k,\varrho,\tau)]\cr}
$$
and
$$\eqal{
\lc{\|(w-k)_+\|_{V_2(Q(\varrho-\sigma_1\varrho,\tau-\sigma_2\tau))}^2}
\rc{\le\gamma\{[(\sigma_1\varrho)^{-2}+(\sigma_2\tau)^{-1}]
\|(w-k)_+\|_{L_2(Q(\varrho,\tau))}^2+
\mu^{{2\over r}(1+\kappa)}(k,\varrho,\tau)\}.}}
$$
\vskip6pt
\noindent
Here the following notation is used:
$$\eqal{
\lc{(w-k)_+=\max\{w-k,0\}\ -\ {\rm the\ truncation\ of}\ w,}
\lc{B_\varrho(x_0)=\{x\in\Omega:\ |x-x_0|<\varrho\}\ -\ {\rm a\ ball\ in}\ 
\Omega,}
\lc{Q(\varrho,\tau)=B_\varrho(x_0)\times(t_0,t_0+\tau)=\{(x,t)\in\Omega^T:\ 
|x-x_0|<\varrho,}
\rc{t_0<t<t_0+\tau\}\ -\ {\rm a\ cylinder\ in}\ \Omega^T,}}
$$
where $\varrho$, $\tau$ are arbitrary positive numbers, $\sigma_1$, 
$\sigma_2$ are arbitrary numbers from the interval $(0,1)$, and $k$ is an 
arbitrary number such that
$$
\ess\sup_{Q(\varrho,\tau)}w(x,t)-k<\delta.
$$
Moreover,
$$\eqal{
V_2(\Omega^T)&=L_\infty(0,T;L_2(\Omega))\cap L_2(0,T;W_2^1(\Omega)),\cr
\mu(k,\varrho,\tau)&=\intop_{t_0}^{t_0+\tau}{\rm meas}^{r/q}A_{k,\varrho}(t)dt,
\cr}
$$
where
$$
A_{k,\varrho}(t)=\{x\in B_\varrho(x_0):\ w(x,t)>k\},
$$
and positive numbers $q$, $r$ are linked by the relation
$$
{1\over r}+{n\over2q}={n\over4},
$$
with the admissible ranges
$$\eqal{
&q\in\bigg(2,{2n\over n-2}\bigg],\ \ &r\in[2,\infty)\quad &{\rm for}\ \ n\ge3,
\cr
&q\in(2,\infty),\ \ &r\in(2,\infty)\quad &{\rm for}\ \ n=2,\cr
&q\in(2,\infty],\ \ &r\in[4,\infty)\quad &{\rm for}\ \ n=1.\cr}
$$

\begin{lemma}\label{4.20.} 
Let the assumptions of Lemma 4.17 be satisfied, and there exist constants 
$\theta_*$ and $M$ such that
$$\eqal{
\theta\ge\theta_*>0,\cr
M\equiv\|\theta\|_{L_\infty(\Omega^T)}&\le c(T)A_{11}(s,\sigma,p',T),\cr
\quad\|\theta\|_{W_2^{2,1}(\Omega^T)}&\le\varphi(T,A_{11}(s,\sigma,p',T)),\cr
\quad\|\epsi_t\|_{{\bm L}_p(\Omega^T)}&\le c(T)[A_{11}(s,\sigma,p',T)+A_5(p,p)],\cr}
\leqno(4.87)
$$
where $p\in(1,\infty)$, $A_{11}(s,\sigma,p',T)$ is defined by (4.82) and 
$A_5(p,p)$ by (4.53). Besides, let $k$ be a positive number such that
$$
\sup_\Omega\theta_0(x)<k
\leqno(4.88)
$$
and
$$
M-k<\delta\quad \textsl{with\ some}\quad \delta>0.
$$
Then
$$
\theta\in{\cal B}_2(\Omega^T,M,\gamma,r,\delta,\kappa)
\leqno(4.89)
$$
with
$$\eqal{
&r=q={10\over3},\quad \kappa\in\bigg(0,{2\over3}\bigg),\cr
&\gamma=c(t)\bigg(1+A_{11}^2(s,\sigma,p',T)+A_5^2\bigg({5\over2-3\kappa},
{5\over2-3\kappa}\bigg)\cr
&\quad+A_5^2\bigg({10\over2-3\kappa},
{10\over2-3\kappa}\bigg)+\|g\|_{L_{5\over2-3\kappa}(\Omega^T)}\bigg),\cr}
$$
$s\in\left({5\over3},2\right),{3\over2-2/s}<p'\le{3\over3/s-1}$, $\sigma>4$.
\end{lemma}

\begin{proof} 
Note that the bound $(4.87)_1$ is ensured by Lemma 4.1, the bound $(4.87)_2$ 
by Lemma 4.17 and $(4.87)_{3,4}$ by Corollaries 4.18 and 4.19.

\noindent
We check that $\theta$ satisfies the conditions (i)--(iii) in the definition 
of ${\cal B}_2(\Omega^T,M,\gamma,r,\delta,\kappa)$.

>From [1, Chap. 3, Sec. 10] we have the imbedding \\
$W_2^{2,1}(\Omega^T)\subset V_2^{1,0}(\Omega^T)$. Since 
$\theta\in W_2^{2,1}(\Omega^T)$ the condition (i) is satisfied.

Condition (ii) is automatically satisfied on account of $(4.87)_2$.

Let us check that $\theta$ satisfies the second inequality in condition (iii). 
By virtue of $(4.87)_1$ it sufficies to consider (iii) with 
$w(x,t)=\theta(x,t)$.

Let $Q(\varrho,\tau)=B_\varrho(x_0)\times(t_0,t_0+\tau)$ be an 
arbitrary cylinder in $\Omega^T$, and $\zeta(x,t)$ be a smooth function such 
that $\supp\zeta(x,t)\subset Q(\varrho,\tau)$ and $\zeta(x,t)=1$ for 
$(x,t)\in Q(\varrho-\sigma_1\varrho,\tau-\sigma_2\tau)$, where 
$\sigma_1,\sigma_2\in(0,1)$. Moreover, let
$$
A_{k,\varrho}(t)=\{x\in B_\varrho(x_0):\ \theta(x,t)>k\}.
$$
Multiplying equation (1.2) by $\zeta^2(\theta-k)_+$ and integrating over 
$\Omega$ gives
$$\eqal{
&{c_v\over2}\intop_\Omega\theta\zeta^2\partial_t(\theta-k)_+^2dx+k_0
\intop_\Omega|\nabla(\theta-k)_+|^2\zeta^2dx\cr
&\quad+2k_0\intop_\Omega\zeta(\theta-k)_+\nabla(\theta-k)_+\cdot\nabla\zeta dx=
\intop_\Omega G\zeta^2(\theta-k)_+dx,\cr}
\leqno(4.90)
$$
where, for simplicity, the right-hand side of (1.2) is denoted by
$$
G\equiv-\theta(\A_2\alphaa)\cdot\epsi_t+(\A_1\epsi_t)\cdot\epsi_t+g,
$$
and to avoid the notational collision the letter $k$ for heat conductivity is 
replaced by $k_0$ (for this proof only).

\noindent
Let us rearrange the first integral on the left-hand side of (4.90) to the form
$$\eqal{
&{c_v\over2}\intop_\Omega\theta\zeta^2\partial_t(\theta-k)_+^2dx=
{c_v\over2}{d\over dt}\intop_\Omega\theta(\theta-k)_+^2\zeta^2dx\cr
&\quad-{c_v\over2}\intop_\Omega\theta_t(\theta-k)_+^2\zeta^2dx-c_v
\intop_\Omega\theta(\theta-k)_+^2\zeta\zeta_tdx.\cr}
\leqno(4.91)
$$
Further, the term with $\theta_t$ in (4.91) is rearranged as follows
$$\eqal{
&-{c_v\over2}\intop_\Omega\theta_t(\theta-k)_+^2\zeta^2dx=-{c_v\over2}
\intop_\Omega(\theta-k)_+^2\partial_t(\theta-k)_+\zeta^2dx\cr
&=-{c_v\over6}\intop_\Omega\partial_t(\theta-k)_+^3\zeta^2dx\cr
&=-{c_v\over6}{d\over dt}\intop_\Omega(\theta-k)_+^3\zeta^2dx+{c_v\over3}
\intop_\Omega(\theta-k)_+^3\zeta\zeta_tdx.\cr}
\leqno(4.92)
$$
Inserting (4.91) and (4.92) into (4.90) we obtain the identity
$$\eqal{
&{c_v\over2}{d\over dt}\intop_\Omega\theta(\theta-k)_+^2\zeta^2dx+k_0
\intop_\Omega|\nabla(\theta-k)_+|^2\zeta^2dx\cr
&={c_v\over6}{d\over dt}\intop_\Omega(\theta-k)_+^3\zeta^2dx-{c_v\over3}
\intop_\Omega(\theta-k)_+^3\zeta\zeta_tdx\cr
&\quad+c_v\intop_\Omega\theta(\theta-k)_+^2\zeta\zeta_tdx-2k_0\intop_\Omega
\zeta(\theta-k)_+\nabla(\theta-k)_+\cdot\nabla\zeta dx\cr
&\quad+\intop_\Omega G\zeta^2(\theta-k)_+dx.\cr}
\leqno(4.93)
$$
Integrating (4.93) with respect to time, taking into account $(4.87)_1$ and 
the fact that by (4.88), $(\theta_0-k)_+=0$, we conclude that
$$\eqal{
&{c_v\theta_*\over2}\intop_\Omega(\theta-k)_+^2\zeta^2dx+k_0\intop_{\Omega^t}
|\nabla(\theta-k)_+|^2\zeta^2dxdt'\cr
&\le c\intop_\Omega(\theta-k)_+^3\zeta^2dx+c\intop_{\Omega^t}(\theta-k)_+^3
|\zeta|\,|\zeta_{t'}|dxdt'\cr
&\quad+c\intop_{\Omega^t}\theta(\theta-k)_+^2|\zeta|\,
|\zeta_{t'}|dxdt'
+c\intop_{\Omega^t}(\theta-k)_+|\nabla(\theta-k)_+|\,|\zeta|\,
|\nabla\zeta|dxdt'\cr
&\quad+c\intop_{\Omega^t}|G|(\theta-k)_+\zeta^2dxdt'\equiv\sum_{i=1}^5I_i.\cr}
\leqno(4.94)
$$
Since, by the assumption $M-k<\delta$, it holds
$$
(\theta-k)_+^3\le\delta(\theta-k)_+^2\ \ {\rm with\ arbitrary}\ \delta>0,
\leqno(4.95)
$$
the integral $I_1$ can be absorbed by the left-hand side of (4.94). Further, 
using $(4.87)_2$, (4.95) and the bound $|\zeta|\le1$,
$$
I_2+I_3\le M\intop_{\Omega^t}(\theta-k)_+^2|\zeta_{t'}|dxdt'.
$$
Next, by Young's inequality
$$
I_4\le{k_0\over2}\intop_{\Omega^t}|\nabla(\theta-k)_+|^2\zeta^2dxdt'+
{c\over2k_0}\intop_{\Omega^t}(\theta-k)_+^2|\nabla\zeta|^2dxdt',
$$
so the first integral on the right-hand side of the latter inequality is 
absorbed by the left-hand side of (4.94). In result, incorporating the above 
estimates into (4.94), we arrive at
$$\eqal{
&\intop_\Omega(\theta-k)_+^2\zeta^2dx+\intop_{\Omega^t}|\nabla(\theta-k)_+|^2
\zeta^2dxdt'\cr
&\le c(M+1)\intop_{\Omega^t}(\theta-k)_+^2(|\nabla\zeta|^2+|\zeta_{t'}|)dxdt'
+I_5,\cr}
\leqno(4.96)
$$
where
$$
I_5\equiv c\intop_{\Omega^t}|G|(\theta-k)_+\zeta^2dxdt'.
$$
By the definition of $\zeta$, it holds
$$\eqal{
&M\intop_{\Omega^t}(\theta-k)_+^2(|\nabla\zeta|^2+|\zeta_{t'}|^2)dxdt'\cr
&\le M[(\sigma_1\varrho)^{-2}+(\sigma_2\tau)^{-1}]\intop_{Q(\varrho,\tau)}
(\theta-k)_+^2dxdt'.\cr}
\leqno(4.97)
$$

It remains to estimate the integral $I_5$. Recalling $(4.87)_2$ again and 
applying H\"older's inequality yields
$$\eqal{
I_5&=c\intop_{t_0}^{t_0+\tau}\intop_{A_{k,\varrho}(t')}|G|(\theta-k)_+\zeta^2
dxdt'\cr
&\le M\bigg(\intop_{t_0}^{t_0+\tau}\intop_{A_{k,\varrho}(t')}|G|^{\lambda_1}
dxdt'\bigg)^{1/\lambda_1}\bigg(\intop_{t_0}^{t_0+\tau}{\rm meas}\,
A_{k,\varrho}(t')dt'\bigg)^{1/\lambda_2},\cr}
$$
where $1/\lambda_1+1/\lambda_2=1$. To satisfy the conditions in the definition 
of the space ${\cal B}_2(\Omega^T,M,\gamma,r,\delta,\kappa)$ we set
$$
{1\over\lambda_2}={2\over r}(1+\kappa)\quad {\rm and}\quad r=q\quad 
{\rm with}\quad {1\over r}+{3\over2r}={3\over4}.
$$
Then
$$
\lambda_2={5\over3(1+\kappa)},\quad
\lambda_1={5\over2-3\kappa}\in\left({5\over2},\infty\right)\quad {\rm for}
\quad \kappa\in\left(0,{2\over3}\right).
$$
Consequently,
$$
I_5\le M\|G\|_{L_{\lambda_1}(\Omega^T)}\mu^{1/\lambda_2}(k,\varrho,\tau).
$$
By virtue of $(4.87)_2$ and $(4.87)_4$, we obtain
$$\eqal{
&\|G\|_{L_{5\over2-3\kappa}(\Omega^T)}
\le c(\|\theta\|_{L_\infty(\Omega^T)}
\|\epsi_t\|_{{\bm L}_{5\over2-3\kappa}(\Omega^T)}+
\|\epsi_t\|_{{\bm L}_{10\over2-3\kappa}(\Omega^T)}^2\cr
&\quad+\|g\|_{{\bm L}_{5\over2-3\kappa}(\Omega^T)})\le 
c(T)\bigg(A_{11}^2(s,\sigma,p',T)+A_5^2\bigg({5\over2-3\kappa},
{5\over2-3\kappa}\bigg)\cr
&\quad+A_5^2\bigg({10\over2-3\kappa},{10\over2-3\kappa}\bigg)+
\|g\|_{L_{5\over2-3\kappa}(\Omega^T)}\bigg)\equiv\gamma_1\cr}
$$
\goodbreak
\noindent
Hence,
$$
I_5\le M\gamma_1\mu^{{3\over5}(1+\kappa)}(k,\varrho,\tau).
\leqno(4.98)
$$
In result, applying estimates (4.97) and (4.98) in (4.96) leads to
$$\eqal{
&\|(\theta-k)_+\|_{V_2(Q(\varrho-\sigma_1\varrho,\tau-\sigma_2\tau))}^2
\equiv\ess\sup_{t\in[0,T]}\intop_\Omega(\theta-k)_+^2\zeta^2dx
+\intop_{\Omega^T}|\nabla(\theta-k)_+|^2\zeta^2dxdt\cr
&\le\gamma\{[(\sigma_1\varrho)^{-2}+(\sigma_2\tau)^{-1}]
\|(\theta-k)_+\|_{L_2(Q(\varrho,\tau))}^2+\mu^{{3\over5}(1+\kappa)}
(k,\varrho,\tau)\}.\cr}
$$
This proves the second inequality in condition (iii) with the positive 
number
$$
\gamma=c(T)(A_{11}^2(s,\sigma,p',T)+1+\gamma_1).
\leqno(4.99)
$$

The first inequality in (iii) follows by multiplying (1.2) by 
$\zeta_0^2(\theta-k)_+$, where $\zeta_0(x)$ is a smooth function such that 
$\supp\zeta_0(x)\subset B_\varrho(x_0)$, $\zeta_0(x)=1$ for 
$x\in B_{\varrho-\sigma_1\varrho}(x_0)$, $\sigma_1\in(0,1)$, and next 
integrating over $\Omega\times(t_0,t_0+\tau)$. Then, by repeating the 
presented above estimates we arrive at the following inequality in place of 
(4.96):
$$\eqal{
&\intop_{B_\varrho(x_0)}(\theta-k)_+^2\zeta_0^2dx+\intop_{Q(\varrho,\tau)}
|\nabla(\theta-k)_+|^2\zeta_0dxdt\cr
&\le c(M+1)\bigg[\intop_{B_\varrho(x_0)}(\theta(t_0)-k)_+^2\zeta_0^2dx
+\intop_{Q(\varrho,\tau)}(\theta-k)_+^2|\nabla\zeta_0|^2dxdt\bigg]\cr
&\quad+ c\intop_{Q(\varrho,\tau)}|G|(\theta-k)_+\zeta_0^2dxdt.\cr}
\leqno(4.100)
$$
The last two integrals on the right-hand side of (4.100) are estimated 
respectively by (4.97) with $(\sigma_2\tau)^{-1}=0$, and by (4.98). 
In result, (4.100) provides the first inequality in the condition (iii) with 
$\gamma$ defined in (4.99). The proof is complete.
\end{proof}

\begin{corollary}\label{4.21.} 
By virtue of the imbedding (cf. [14, Thm II. 7.1]) 
$$
{\cal B}_2(\Omega^T,M,\gamma,r,\delta,\kappa)\subset C^{\alpha,\alpha/2}
(\Omega^T),\quad \alpha\in(0,1),
$$
it follows from (4.89) that
$$
\theta\in C^{\alpha,\alpha/2}(\Omega^T)
\leqno(4.101)
$$
with H\"older's exponent $\alpha\in(0,1)$ depending on $M$, $\gamma$, $r$, 
$\delta$ and $\kappa$, where
$$\eqal{
&M=\sup_{\Omega^T}\theta\le c(T)A_{11}(s,\sigma,p',T),\cr
&\gamma(s,\sigma,p',\kappa,T)=c(T)\bigg(1+A_{11}^2(s,\sigma,p',T)+A_5^2
\bigg({5\over2-3\kappa},{5\over2-3\kappa}\bigg)\cr
&\quad+A_5^2\bigg({10\over2-3\kappa},{10\over2-3\kappa}\bigg)+
\|g\|_{L_{5\over2-3\kappa}(\Omega^T)}\bigg),\cr
&s\in\bigg({5\over3},2\bigg),{3\over2-2/s}<p'<{3\over3/s-1},\quad
\sigma>4,\quad \kappa\in\bigg(0,{2\over3}\bigg).\cr}
$$
\end{corollary}

In view of H\"older's continuity of $\theta$ we can apply Lemma 3.5 and 
deduce the final estimates on $\theta$ and $\u$.

\begin{proposition}\label{4.22.} 
Let the quantity
$$\eqal{
A(t)&=\|\u_0\|_{\W_{12/5}^2(\Omega)\cap\W_{p'}^1(\Omega)}+
\|\u_1\|_{\B_{12,12}^{11/6}(\Omega)\cap\B_{p',\sigma'}^{2-2/\sigma'}(\Omega)}
\cr
&\quad+\|\theta_0\|_{W_{p_*}^1(\Omega)}+\|g\|_{L_{\infty,12}(\Omega^t)}+
\|\b\|_{{\bm L}_{12}(\Omega^t)\cap {\bm L}_{p',\sigma'}(\Omega^t)},\cr
&p_*>3,\quad p',\sigma'\in(1,\infty),\cr}
\leqno(4.102)
$$
be finite.\\
Then the following a priori estimate
$$
\|\theta\|_{L_\infty(\Omega^t)}+\|\theta\|_{C^{\alpha,\alpha/2}(\Omega^t)}+
\|\epsi_{t'}\|_{{\bm L}_{p',\sigma'}(\Omega^t)}
\le\varphi(A(t)),\quad t\le T,
\leqno(4.103)
$$
is valid.
\end{proposition}

\begin{proof} 
>From (4.82), (4.86) and (4.101) it follows that the estimates on $\|\theta\|_{L_\infty(\Omega^t)}$, 
$\|\epsi_t\|_{{\bm L}_{p,\sigma'}(\Omega^t)}$ and 
$\|\theta\|_{C^{\alpha,\alpha/2}(\Omega^t)}$ involve the parameters 
$s\in\left({5\over3},2\right),p'\in\Big({3\over2-2/s},{3\over3/5-1}\Big]$, $\sigma>4$, 
$(p,\sigma')\in(1,\infty)$ and $\kappa\in\left(0,{2\over3}\right)$. 
We select these parameters in such a way to express the above mentioned 
estimates in an explicit, sufficiently simple way.

\noindent
Let us recall that the estimates depend on the quantities:
$$\eqal{
&A_8(s,\sigma,t)=c(t)[\|\u_0\|_{\H^2(\Omega)}+
\|\u_0\|_{\W_{2s\over2-s}^1(\Omega)}+\|\u_1\|_{\B_{2,10}^{9/5}(\Omega)}\cr
&\quad+\|\u_1\|_{\B_{2,\sigma}^{2-2/\sigma}(\Omega)}+
\|\u_1\|_{\B_{{2s\over2-s},{2s\over2-s}}^{3s-2\over s}(\Omega)}+
\|\theta_0\|_{H^1(\Omega)}\cr
&\quad+\|\theta_0\|_{L_{2s\over2-s}(\Omega)}+\|\b\|_{{\bm L}_{2,10}(\Omega^t)}+
\|\b\|_{{\bm L}_{2,\sigma}(\Omega^t)}+
\|\b\|_{{\bm L}_{2s\over2-s}(\Omega^t)}\cr
&\quad+\|g\|_{L_2(\Omega^t)}+\|g\|_{L_{2s\over2-s}(\Omega^t)}],\cr
&A_9(p',s)=\|\u_1\|_{\B_{p',s}^{2-2/s}(\Omega)}+
\|\b\|_{{\bm L}_{p',s}(\Omega^t)},\cr
&A_{10}(s,\sigma,p',t)=c(t)(\varphi(A_8(s,\sigma,t))+A_9(p',s)),\cr
&A_{11}(s,\sigma,p',t)=c(t)(A_{10}(s,\sigma,p',t)+A_{10}^2(s,\sigma,p',t)\cr
&\quad+\|\theta_0\|_{L_\infty(\Omega)}+\|g\|_{L_{\infty,1}(\Omega^t)}),\cr
&\gamma(s,\sigma,p',\kappa,t)=c\bigg(1+A_{11}^2(s,\sigma,p',t)+A_5^2
\bigg({5\over2-3\kappa},{5\over2-3\kappa}\bigg)\cr
&\quad+A_5^2\bigg({10\over2-3\kappa},{10\over2-3\kappa}\bigg)+
\|g\|_{L_{5\over2-3\kappa}(\Omega^t)}\bigg),\cr
&A_5\bigg({5\over2-3\kappa},{5\over2-3\kappa}\bigg)
=\|\u_0\|_{\W_{5\over2-3\kappa}^1(\Omega)}+
\|\u_1\|_{\B_{{5\over2-3\kappa},{5\over2-3\kappa}}^{6(1+\kappa)\over5}
(\Omega)}\cr
&\quad+\|\b\|_{{\bm L}_{5\over2-3\kappa}(\Omega^t)},\cr
&A_5\bigg({10\over2-3\kappa},{10\over2-3\kappa}\bigg)
=\|\u_0\|_{\W_{10\over2-3\kappa}^1(\Omega)}+
\|\u_1\|_{\B_{{10\over2-3\kappa},{10\over2-3\kappa}}^{8+3\kappa\over5}
(\Omega)}\cr
&\quad+\|\b\|_{{\bm L}_{10\over2-3\kappa}(\Omega)},\cr
&A_5(p,\sigma')=\|\u_0\|_{\W_p^1(\Omega)}+
\|\u_1\|_{\B_{p,\sigma'}^{2-2/\sigma'}(\Omega)}+
\|\b\|_{{\bm L}_{p,\sigma'}(\Omega^t)},\cr}
\leqno(4.104)
$$
where $s\in\left({5\over3},2\right)$, 
$p'\in\left.\Big({3\over2-{2\over s}},{3\over{3\over s}-1}\Big]\right.$, 
$\sigma>4$, $\kappa\in\left(0,{2\over3}\right)$.

To have the above quantities finite we need
$$\eqal{
\u_0&\in\H^2(\Omega)\cap\W_{2s\over2-s}^1(\Omega)\cap
\W_{10\over2-3\kappa}^1(\Omega)\cap\W_p^1(\Omega),\cr
\u_1&\in\B_{2,10}^{9/5}(\Omega)\cap\B_{2,\sigma}^{2-2/\sigma}(\Omega)\cap
\B_{{2s\over2-s},{2s\over2-s}}^{3s-2\over s}(\Omega)\cap
\B_{p',s}^{2-2/s}(\Omega)\cr
&\quad\cap\B_{{5\over2-3\kappa},{5\over2-3\kappa}}^{6(1+\kappa)\over5}
(\Omega)\cap
\B_{{10\over2-3\kappa},{10\over2-3\kappa}}^{8+3\kappa\over5}(\Omega)
\cap\B_{p,\sigma'}^{2-2/r'}(\Omega),\cr
\theta_0&\in L_\infty(\Omega)\cap H^1(\Omega),\cr
\b&\in{\bm L}_{2,10}(\Omega^t)\cap{\bm L}_{2,\sigma}(\Omega^t)\cap
{\bm L}_{2s\over2-s}(\Omega^t)\cap{\bm L}_{p',s}(\Omega^t)\cap
{\bm L}_{10\over2-3\kappa}(\Omega^t)\cr
&\quad\cap{\bm L}_{p,\sigma'}(\Omega^t),\cr
g&\in L_2(\Omega^t)\cap L_{2s\over2-s}(\Omega^t)\cap L_{\infty,1}(\Omega^t)
\cap L_{5\over2-3\kappa}(\Omega^t),\cr}
\leqno(4.105)
$$
where $s\in\left({5\over3},2\right)$, 
$p'\in\left.\Big({3\over2-{2\over s}},{3\over{3\over s}-1}\Big]\right.$, 
$\sigma>4$, $\kappa\in\left(0,{2\over3}\right)$.

Let us set
$$
s={12\over7},\quad p'=4,\quad \kappa={1\over9},\quad \sigma=10.
$$
Then ${2s\over2-s}=12$, ${5\over2-3\kappa}=3$, ${10\over2-3\kappa}=6$, 
and (4.105) takes the form
$$\eqal{
&\u_0\in\H^2(\Omega)\cap\W_{12}^1(\Omega)\cap\W_p^1(\Omega),\cr
&\u_1\in\B_{2,10}^{9/5}(\Omega)\cap\B_{12,12}^{11/6}(\Omega)\cap
\B_{6,6}^{5/6}(\Omega)\cap\B_{4,12/7}^{5/6}(\Omega)
\cap\B_{p,\sigma'}^{2-2/\sigma'}(\Omega),\cr
&\theta_0\in L_\infty(\Omega)\cap H^1(\Omega),\cr
&\b\in{\bm L}_{12}(\Omega^t)\cap{\bm L}_{p,\sigma'}(\Omega^t),\cr
&g\in L_{\infty,1}(\Omega^t)\cap L_{12}(\Omega^t)=L_{\infty,12}(\Omega^t).\cr}
\leqno(4.106)
$$
Using the imbeddings
$$\eqal{
&\W_{12/5}^2(\Omega)\subset\H^2(\Omega)\cap\W_{12}^1(\Omega),\cr
&\B_{12,12}^{11/6}(\Omega)\subset\B_{2,10}^{9/5}(\Omega)
\cap\B_{4,12/7}^{5/6}(\Omega),\cr
&W_p^1(\Omega)\subset L_\infty(\Omega)\cap H^1(\Omega)\quad{\rm for}\ \ p>3,
\cr}
$$
we replace (4.106) by
$$\eqal{
&\u_0\in\W_{12/5}^2(\Omega)\cap\W_p^1(\Omega),\quad 
&\u_1\in\B_{12,12}^{11/6}(\Omega)\cap\B_{p,\sigma'}^{2-2/\sigma'}(\Omega),\cr
&\theta_0\in W_{p_*}^1(\Omega),\quad p_*>3,\cr
&\b\in{\bm L}_{12}(\Omega^t)\cap{\bm L}_{p,\sigma'}(\Omega^t),\quad 
&g\in L_{\infty,12}(\Omega^t).\cr}
\leqno(4.107)
$$
The assumptions (4.107) ensure that the quantity $A(t)$ in (4.102) is finite. 
Then estimates (4.82), (4.86) and (4.101) imply (4.103). Let us note that $p$ 
is replaced by $p'$. This completes the proof.
\end{proof}

\begin{lemma}\label{4.23.} 
Assume that the data satisfy (4.107), $S\in C^2$, and moreover,
$$\eqal{
&\u_0\in\W_{2q}^1(\Omega),\quad 
&\u_1\in\B_{2q,2q_0}^{2-1/q_0}(\Omega)\cap\B_{p,p_0}^{2-2/p_0}(\Omega),\cr
&\theta_0\in B_{q,q_0}^{2-2/q_0}(\Omega),\quad
&\b\in{\bm L}_{p,p_0}(\Omega^t)\cap{\bm L}_{2q,2q_0}(\Omega^t),\cr
&g\in L_{q,q_0}(\Omega^t),\cr}
\leqno(4.108)
$$
with any numbers $p,p_0,q,q_0\in(1,\infty)$ such that
$$
{3\over q}+{2\over q_0}-{3\over p}-{2\over p_0}\le1.
$$
Then
$$\eqal{
\|\theta\|_{W_{q,q_0}^{2,1}(\Omega^t)}
&\le\varphi(A(t),t)+c(\|\b\|_{{\bm L}_{2q,2q_0}(\Omega^t)}^2\cr
&\quad+\|\u_0\|_{\W_{2q}^1(\Omega)}^2+
\|\u_1\|_{\B_{2q,2q_0}^{2-1/q_0}(\Omega)}+\|g\|_{L_{q,q_0}(\Omega^t)}\cr
&\quad+\|\theta_0\|_{B_{q,q_0}^{2-2/q_0}(\Omega)})\equiv B_1(q,q_0,t)\cr}
\leqno(4.109)
$$
and
$$\eqal{
&\|\u_t\|_{\W_{p,p_0}^{2,1}(\Omega^t)}\le c(t)B_1(q,q_0,t)
+c(\|\b\|_{{\bm L}_{p,p_0}(\Omega^t)}+
\|\u_1\|_{\B_{p,p_0}^{2-2/p_0}(\Omega)}),\cr}
\leqno(4.110)
$$
where $A(t)$ is defined by (4.102), and $t\le T$.
\end{lemma}

We remark that, by Lemmas 3.4 and 3.5, the constant $c$ in (4.109) and (4.110) 
depends on $p$, $p_0$, $q$, $q_0$ and $C^2$-norm of the boundary $S$.

\begin{proof} 
In view of the H\"older continuity of $\theta$ along with its lower and 
upper boundedness we deduce from Lemma 3.5 and Proposition 4.22 that
$$\eqal{
\|\theta\|_{W_{q,q_0}^{2,1}(\Omega^t)}
&\le c(\|\theta(\A_2\alphaa)\cdot\epsi_{t'}\|_{L_{q,q_0}(\Omega^t)}+\|(\A_1\epsi_{t'})\cdot\epsi_{t'}\|_{L_{q,q_0}
(\Omega^t)}\cr
&\quad+\|g\|_{L_{q,q_0}(\Omega^t)}+\|\theta_0\|_{B_{q,q_0}^{2-2/q_0}(\Omega)})\cr
&\le\varphi(A(t),t)\|\epsi_{t'}\|_{{\bm L}_{q,q_0}(\Omega^t)}+c
(\|\epsi_{t'}\|_{{\bm L}_{2q,2q_0}(\Omega^t)}^2\cr
&\quad+\|g\|_{L_{q,q_0}(\Omega^t)}+
\|\theta_0\|_{B_{q,q_0}^{2-2/q_0}(\Omega)}),\quad t\le T.\cr}
\leqno(4.111)
$$
>From (4.103) with the parameters $(r,r_0)=(p,\sigma')$ selected in Proposition 
4.22, we have
$$
\|\epsi_{t'}\|_{{\bm L}_{r,r_0}(\Omega^t)}\le\varphi(A(t),t),
\leqno(4.112)
$$
where $(r,r_0)$ is equal $(2q,2q_0)$. \\
Using (4.112) in (4.111) we conclude (4.109).

Let us consider now the imbedding
$$
\nabla W_{q,q_0}^{2,1}(\Omega^T)\subset L_{q',q'_0}(\Omega^T),
$$
which holds under the condition
$$
{3\over q}+{2\over q_0}-{3\over q'}-{2\over q'_0}\le1,\quad
q'\ge q,\quad q'_0\ge q_0.
\leqno(4.113)
$$
Then it follows from (4.75) and (4.112) that
$$\eqal{
&\|\u_{t'}\|_{\W_{p,p_0}^{2,1}(\Omega^t)}\le c(t)
(\|\nabla\theta\|_{{\bm L}_{p,p_0}(\Omega^t)}+
\|\b\|_{{\bm L}_{p,p_0}(\Omega^t)}+\|\u_1\|_{\B_{p,p_0}^{2-2/p_0}(\Omega)})\cr
&\le c(t)(\|\nabla\theta\|_{{\bm L}_{q',q'_0}(\Omega^t)}+
\|\b\|_{{\bm L}_{p,p_0}(\Omega^t)}+\|\u_1\|_{\B_{p,p_0}^{2-2/p_0}(\Omega)})\cr
&\le c(t)(\|\theta\|_{W_{q,q_0}^{2,1}(\Omega^t)}+
\|\b\|_{{\bm L}_{p,p_0}(\Omega^t)}+\|\u_1\|_{\B_{p,p_0}^{2-2/p_0}(\Omega)}),
\cr}
\leqno(4.114)
$$
where $p\le q'$, $p_0\le q'_0$, and $q',q'_0$ satisfy (4.113).
Finally, we set $p=q'$, $p_0=q'_0$.
This establishes (4.110) and thereby completes the proof.
\end{proof}

\begin{remark}\label{4.24.} 
Since $\u=\intop_0^t\u_{t'}dt'+\u_0$, 
$\u_t\in\W_{p,p_0}^{2,1}(\Omega^T)$, $p,p_0\in(1,\infty)$, and 
$\u_0\in\B_{p,\sigma_0}^{2-2/\sigma_0}(\Omega)$, it follows that 
$\u\in\W_{p,\sigma_0}^{2,1}(\Omega^T)$ where $\sigma_0\in[1,\infty)$.
\end{remark}

\begin{corollary}\label{4.25.}
Comparing assumptions (4.107) with $p=p'$ and (4.108) we deduce that
$$\eqal{
&\u_0\in\W_{12/5}^2(\Omega)\cap\W_{p'}^1(\Omega)\cap\W_{2q}^1(\Omega),\cr
&\u_1\in\B_{12,12}^{11/6}(\Omega)\cap\B_{p',\sigma'}^{2-2/\sigma'}(\Omega)\cap
\B_{2q,2q_0}^{2-1/q_0}(\Omega)\cap\B_{p,p_0}^{2-2/p_0}(\Omega).\cr}
\leqno(4.115)
$$
Setting $p'=2q$, $\sigma'=2q_0$, $q=q_0=6$, $p=p_0=12$, we get
$$
\u_0\in\W_{12/5}^2(\Omega),\quad \u_1\in\B_{12,12}^{11/6}(\Omega),
$$
and
$$\eqal{
&\theta_0\in W_{p_*}^1(\Omega)\cap B_{6,6}^{5/3}(\Omega)=B_{6,6}^{5/3}(\Omega)
\quad {\rm for}\ \ 3<p_*<6,\cr
&\b\in{\bm L}_{12}(\Omega^t)\cap{\bm L}_{p,p_0}(\Omega^t)\cap
{\bm L}_{2q,2q_0}(\Omega^t)={\bm L}_{12}(\Omega^t),\cr
&g\in L_{\infty,12}(\Omega^t)\cap L_{q,q_0}(\Omega^t)=
L_{\infty,12}(\Omega^t).\cr}
$$
Introducing the quantity
$$\eqal{
B(t)&=\|\u_0\|_{\W_{12/5}^2(\Omega)\cap\W_{12}^1(\Omega)}+
\|\u_1\|_{\B_{12,12}^{11/6}(\Omega)}\cr
&\quad+\|\b\|_{{\bm L}_{12}(\Omega^t)}+\|\theta_0\|_{B_{6,6}^{5/3}(\Omega)}+
\|g\|_{L_{\infty,12}(\Omega^t)},\cr}
$$
we conclude in place of (4.109) and (4.110) the estimate
$$
\|\u_t\|_{\W_{12}^{2,1}(\Omega^t)}+\|\theta\|_{W_6^{2,1}(\Omega^t)}\le
\varphi(B(t)).
\leqno(4.116)
$$
\end{corollary}

\section{Existence (proof of Theorem A)}

To prove the existence of solutions to problem (1.1)--(1.4) we use the 
following method of successive approximations:
$$\eqal{
&\u_{tt}^{n+1}-\nabla\cdot(\A_1\epsi(\u_t^{n+1}))=\nabla\cdot[\A_2\epsi(\u^n)-
(\A_2\alphaa)\theta^n]+\b\quad &{\rm in}\ \ \Omega^T,\cr
&c_v\theta^n\theta_t^{n+1}-k\Delta\theta^{n+1}=\quad\cr
&\quad-\thetaa^n(\A_2\alphaa)\cdot\epsi(\u_t^n)+\A_1\epsi(\u_t^n)\cdot
\epsi(\u_t^n)+g\quad &{\rm in}\ \ \Omega^T,\cr
&\u^{n+1}=\0\quad &{\rm on}\ \ S^T,\cr
&\n\cdot\nabla\theta^{n+1}=0\quad &{\rm on}\ \ S^T,\cr
&\u^{n+1}|_{t=0}=\u_0,\ \ \u_t^{n+1}|_{t=0}=\u_1\quad &{\rm in}\ \ \Omega,\cr
&\theta^{n+1}|_{t=0}=\theta_0\quad &{\rm in}\ \ \Omega,\cr}
\leqno(5.1)
$$
where $\u^n$, $\theta^n$, $n\in{\N}\cup\{0\}$, are treated as given.

\noindent
Moreover, the approximation $(\u^0,\theta^0)$ is constructed by an extension 
of the initial data in such a way that
$$
\u^0|_{t=0}=\u_0,\ \ \u_t^0|_{t=0}=\u_1,\ \ \theta^0|_{t=0}=\theta_0\quad 
{\rm in}\ \ \Omega,
\leqno(5.2)
$$
and
$$
\u^0=\0,\ \ \n\cdot\nabla\theta^0=0\quad {\rm on}\ \ S^T.
\leqno(5.3)
$$

First we show that the sequence $\{\u^n,\theta^n\}$ is uniformly bounded.

\begin{lemma}\label{5.1.} (Boundedness of the approximation) 
Assume that
$$\eqal{
D=D(p,p_0,q,q_0)&\equiv\|\u_0\|_{\W_p^2(\Omega)}+
\|\u_1\|_{\B_{p,p_0}^{2-2/p_0}(\Omega)}+
\|\theta_0\|_{B_{q,q_0}^{2-2/q_0}(\Omega)}\cr
&\quad+\|\b\|_{{\bm L}_{p,p_0}(\Omega^\tau)}+
\|g\|_{L_{q,q_0}(\Omega^\tau)}<\infty,\cr}
$$
where $p,p_0,q,q_0\in(1,\infty)$ satisfy the conditions
$$\eqal{
&{3\over q}+{2\over q_0}-{3\over p}-{2\over p_0}<1,\cr
&{3\over p}+{2\over p_0}-{3\over2q}-{2\over2q_0}<1.\cr}
$$
Moreover, assume that $\theta_0\ge{\underline\theta}>0$, $\tau>0$ is 
sufficiently small, and $n\in{\N}\cup\{0\}$. Then there exists a constant 
$A$ independent of $n$ but depending on $D$, $p,$ $p_0$, $q$, $q_0$ such that 
solutions to problem (5.1) satisfy the estimates
$$\eqal{
&\|\u_t^n\|_{\W_{p,p_0}^{2,1}(\Omega^\tau)}+
\|\theta^n\|_{W_{q,q_0}^{2,1}(\Omega^\tau)}\le A,\cr
&{1\over2}{\underline\theta}\le\sup_{\Omega^\tau}\theta^n\le A\cr}
\leqno(5.4)
$$
for $n\in{\N}\cup\{0\}$.
\end{lemma}
\goodbreak

\begin{proof} 
Let $\u^n\in\W_{p,p_0}^{2,1}(\Omega^\tau)$, 
$\theta^n\in W_{q,q_0}^{2,1}(\Omega^\tau)$. Then Lemma 3.4 ensures the 
existence of solutions to problem $(5.1)_{1,3,5}$ and the estimate
$$\eqal{
\|\u_t^{n+1}\|_{\W_{p,p_0}^{2,1}(\Omega^\tau)}&\le c
(\|\nabla^2\u^n\|_{{\bm L}_{p,p_0}(\Omega^\tau)}+
\|\nabla\theta^n\|_{{\bm L}_{p,p_0}(\Omega^\tau)}\cr
&\quad+\|\b\|_{{\bm L}_{p,p_0}(\Omega^\tau)}+
\|\u_1\|_{\B_{p,p_0}^{2-2/p_0}(\Omega)}),\cr}
\leqno(5.5)
$$
where constant $c$ does not depend on $\tau$.

\noindent
With the use of the formula
$$
\u^n(x,t)=\intop_0^t\u_{t'}^n(x,t')dt'+\u_0(x),
$$
the first term on the right-hand side of (5.5) is estimated by
$$\eqal{
&\|\nabla^2\u^n\|_{{\bm L}_{p,p_0}(\Omega^\tau)}\le\bigg\|\intop_0^t
\nabla^2\u_{t'}^ndt'\bigg\|_{{\bm L}_{p,p_0}(\Omega^\tau)}
+\tau^{1/p_0}\|\nabla^2\u_0\|_{{\bm L}_p(\Omega)}\cr
&\le\tau\|\nabla^2\u_t^n\|_{{\bm L}_{p,p_0}(\Omega^\tau)}+\tau^{1/p_0}
\|\nabla^2\u_0\|_{{\bm L}_p(\Omega)},\cr}
\leqno(5.6)
$$
where in the last line the H\"older inequality was used.

\noindent
The second term on the right-hand side of (5.5) is estimated with the help 
of the interpolation inequality:
$$\eqal{
&\|\nabla\theta^n\|_{{\bm L}_{p,p_0}(\Omega^\tau)}\le\delta
\|\theta^n\|_{W_{q,q_0}^{2,1}(\Omega^\tau)}+c(1/\delta)
\|\theta^n\|_{L_{q,q_0}(\Omega^\tau)}\cr
&\le\delta\|\theta^n\|_{W_{q,q_0}^{2,1}(\Omega^\tau)}+c(1/\delta)
\bigg\|\intop_0^t\theta_{t'}^ndt'+\theta_0\bigg\|_{L_{q,q_0}(\Omega^\tau)}\cr
&\le\delta\|\theta^n\|_{W_{q,q_0}^{2,1}(\Omega^\tau)}+c(1/\delta)
(\tau\|\theta^n\|_{W_{q,q_0}^{2,1}(\Omega^\tau)}+\tau^{1/q_0}
\|\theta_0\|_{L_q(\Omega)}),\cr}
\leqno(5.7)
$$
where $c(1/\delta)\sim\delta^{-a}$, $a>0$. The interpolation inequality in the 
first line holds under the restriction
$$
{3\over q}+{2\over q_0}-{3\over p}-{2\over p_0}<1,\quad p\ge q,\ \ p_0\ge q_0,
\leqno(5.8)
$$
where the last two restrictions can be relaxed on account of the H\"older 
inequality
$$
\|v\|_{L_{r,r_0}(\Omega^\tau)}\le|\Omega|^{1/r-1/\sigma}\tau^{1/r_0-1/\sigma_0}
\|v\|_{L_{\sigma,\sigma_0}(\Omega^\tau)}
$$
with $r\le\sigma$, $r_0\le\sigma_0$.

\noindent
Applying estimates (5.6) and (5.7) in (5.5) leads to
$$\eqal{
&\|\u_t^{n+1}\|_{\W_{p,p_0}^{2,1}(\Omega^\tau)}\le c
[\tau\|\u_t^n\|_{\W_{p,p_0}^{2,1}(\Omega^\tau)}\cr
&\quad+(\delta+c(1/\delta)\tau)\|\theta^n\|_{W_{q,q_0}^{2,1}(\Omega^\tau)}+
\tau^{1/p_0}\|\nabla^2\u_0\|_{{\bm L}_p(\Omega)}\cr
&\quad+c(1/\delta)\tau^{1/q_0}\|\theta_0\|_{L_q(\Omega)}+
\|\b\|_{{\bm L}_{p,p_0}(\Omega^\tau)}+
\|\u_1\|_{\B_{p,p_0}^{2-2/p_0}(\Omega)}]\cr
&\le c[\tau^a(\|\u_t^n\|_{\W_{p,p_0}^{2,1}(\Omega^\tau)}+
\|\theta^n\|_{W_{q,q_0}^{2,1}(\Omega^\tau)})\cr
&\quad+\|\u_0\|_{\W_p^2(\Omega)}+\|\u_1\|_{\B_{p,p_0}^{2-2/p_0}(\Omega)}+
\|\theta_0\|_{L_q(\Omega)}+\|\b\|_{{\bm L}_{p,p_0}(\Omega^\tau)}],\cr}
\leqno(5.9)
$$
where $p$, $p_0$, $q$, $q_0$ satisfy (5.8) with $p,q,p_0,q_0\in(1,\infty)$, 
$\tau$ is sufficiently small, and $a>0$.

By virtue of Lemma 3.3 there exists a solution to problem $(5.1)_{2,4,6}$, 
and the following estimate holds
$$\eqal{
&\|\theta^{n+1}\|_{W_{q,q_0}^{2,1}(\Omega^\tau)}\le\varphi
\Big(\sup_{\Omega^\tau}{1\over\theta^n},\sup_{\Omega^\tau}\theta^n\Big)
[\|\theta^n\epsi(\u_t^n)\|_{{\bm L}_{q,q_0}(\Omega^\tau)}\cr
&\quad+\|\,|\epsi(\u_t)|^2\|_{L_{q,q_0}(\Omega^\tau)}+\|g\|_{L_{q,q_0}(\Omega^\tau)}+
\|\theta_0\|_{B_{q,q_0}^{2-2/q_0}(\Omega)}],\cr}
\leqno(5.10)
$$
where $\varphi$ is an increasing positive function of its arguments and 
$\tau$ is sufficiently small.

\noindent
By the H\"older inequality it follows from (5.10) that
$$\eqal{
&\|\theta^{n+1}\|_{W_{q,q_0}^{2,1}(\Omega^\tau)}\le\varphi
\Big(\sup_{\Omega^\tau}{1\over\theta^n},\sup_{\Omega^\tau}\theta^n\Big)
[\|\theta^n\|_{L_{2q,2q_0}(\Omega^\tau)}\|\nabla\u_t^n\|_{{\bm L}_{2q,2q_0}(\Omega^\tau)}\cr
&\quad+\|\nabla\u_t^n\|_{{\bm L}_{2q,2q_0}(\Omega^\tau)}^2+\|g\|_{L_{q,q_0}(\Omega^\tau)}+
\|\theta_0\|_{B_{q,q_0}^{2-2/q_0}(\Omega)}].\cr}
\leqno(5.11)
$$
We proceed now to estimate the norms on the right-hand side of (5.11). First 
we examine
$$\eqal{
&\sup_{t\le\tau}\|\theta^n\|_{L_\infty(\Omega)}\le\sup_{t\le\tau}(\delta
\|\theta^n\|_{B_{q,q_0}^{2-2/q_0}(\Omega)}+c(1/\delta)
\|\theta^n\|_{L_q(\Omega)})\cr
&\le\sup_{t\le\tau}\bigg[\delta\|\theta^n\|_{B_{q,q_0}^{2-2/q_0}(\Omega)}+
c(1/\delta)\bigg\|\intop_0^t\theta_{t'}^ndt'+\theta_0\bigg\|_{L_q(\Omega)}
\bigg]\cr
&\le\sup_{t\le\tau}\bigg[\delta\|\theta^n\|_{B_{q,q_0}^{2-2/q_0}(\Omega)}+
c(1/\delta)t^{1/q'_0}\bigg(\intop_0^t\|\theta_{t'}^n\|_{L_q(\Omega)}^{q_0}
dt'\bigg)^{1/q_0}\cr
&\quad+c(1/\delta)\|\theta_0\|_{L_q(\Omega)}\bigg]\cr
&\le(\delta+c(1/\delta)\tau^{1/q'_0})
\|\theta^n\|_{W_{q,q_0}^{2,1}(\Omega^\tau)}+c(1/\delta)
\|\theta_0\|_{B_{q,q_0}^{2-2/q_0}(\Omega)},\cr}
\leqno(5.12)
$$
where $c(1/\delta)\sim\delta^{-a}$, $a>0$, $1/q_0+1/q'_0=1$. The first 
inequality in (5.12) expresses the interpolation inequality which holds for
$$
3/q<2-2/q_0.
\leqno(5.13)
$$

The last inequality in (5.12) follows from the imbedding
$$
\sup_{t\le\tau}\|u(t)\|_{B_{q,q_0}^{2-2/q_0}(\Omega)}\le c
(\|u\|_{W_{q,q_0}^{2,1}(\Omega^\tau)}+
\|u(0)\|_{B_{q,q_0}^{2-2/q_0}(\Omega)})
\leqno(5.14)
$$
with a constant $c$ independent of $\tau$, which holds true for any solution 
of a parabolic equation with vanishing boundary conditions and sufficiently 
smooth coefficients.

To estimate $\theta^n$ from below by a positive constant we use the 
assumption $\theta_0\ge{\underline\theta}>0$. Since by (5.13), 
$\theta^n\in W_{q,q_0}^{2,1}(\Omega^\tau)\subset C^\alpha(\Omega^\tau)$ 
with some $\alpha\in(0,1)$, it follows that
$$
\theta^n=\theta^n-\theta_0+\theta_0\ge{\underline\theta}-|\theta^n-\theta_0|
\ge{\underline\theta}-\sup_{x\in\Omega}
\|\theta^n\|_{C^{\alpha/2}(0,\tau)}\tau^{\alpha/2}.
$$
Hence for $\tau$ so small that
$$
\|\theta^n\|_{L_\infty(\Omega;C^{\alpha/2}(0,\tau))}\tau^{\alpha/2}\le c
\|\theta^n\|_{W_{q,q_0}^{2,1}(\Omega^\tau)}\tau^{\alpha/2}\le{1\over2}
{\underline\theta},
\leqno(5.15)
$$
we have
$$
\theta^n\ge{1\over2}{\underline\theta}.
\leqno(5.16)
$$

To estimate the first two terms in the square bracket on the right-hand\break
side \hskip-1pt of \hskip-1pt \hskip-1pt (5.11) \hskip-1pt we \hskip-1pt 
consider \hskip-1pt the \hskip-1pt factors \hskip-1pt 
$\|\theta^n\|_{L_{2q,2q_0}}\!(\Omega^\tau)$ \hskip-1pt and \hskip-1pt 
$\|\nabla\u_t^n\|_{{\bm L}_{2q,2q_0}(\Omega^\tau)}$.\break 
We have
$$\eqal{
&\|\theta^n\|_{L_{2q,2q_0}(\Omega^\tau)}\le\delta
\|\theta^n\|_{W_{q,q_0}^{2,1}(\Omega^\tau)}+c(1/\delta)
\|\theta^n\|_{L_{q,q_0}(\Omega^\tau)}\cr
&\le\delta\|\theta^n\|_{W_{q,q_0}^{2,1}(\Omega^\tau)}+c(1/\delta)
\bigg\|\intop_0^t\theta_{t'}^ndt'+\theta_0\bigg\|_{L_{q,q_0}(\Omega^\tau)}\cr
&\le(\delta+c(1/\delta)\tau)\|\theta^n\|_{W_{q,q_0}^{2,1}(\Omega^\tau)}+
\tau^{1/q_0}\|\theta_0\|_{L_q(\Omega)},\cr}
\leqno(5.17)
$$
where the first inequality is the interpolation inequality which holds under 
the condition
$$
{3\over q}+{2\over q_0}-{3\over2q}-{2\over2q_0}<2,
$$
that is
$$
{3\over q}+{2\over q_0}<4.
\leqno(5.18)
$$
Next,
$$\eqal{
&\|\nabla\u_t^n\|_{{\bm L}_{2q,2q_0}(\Omega^\tau)}\le\delta
\|\u_t^n\|_{\W_{p,p_0}^{2,1}(\Omega^\tau)}+c(1/\delta)
\|\u_t^n\|_{{\bm L}_{p,p_0}(\Omega^\tau)}\cr
&\le\delta\|\u_t^n\|_{\W_{p,p_0}^{2,1}(\Omega^\tau)}+c(1/\delta)
\bigg\|\intop_0^t\u_{tt'}^ndt'+\u_t(0)\bigg\|_{{\bm L}_{p,p_0}(\Omega^\tau)}\cr
&\le(\delta+c(1/\delta)\tau)\|\u_t^n\|_{\W_{p,p_0}^{2,1}(\Omega^\tau)}+
c(1/\delta)\tau^{1/p_0}\|\u_1\|_{{\bm L}_p(\Omega)},\cr}
\leqno(5.19)
$$
where the first inequality holds under the condition
$$
{3\over p}+{2\over p_0}-{3\over2q}-{2\over2q_0}<1,\quad 
p,p_0,q,q_0\in(1,\infty).
\leqno(5.20)
$$
Using (5.12), (5.16), (5.17), (5.19) and choosing appropriately $\delta$ 
so that to get the factors of $\|\theta^n\|_{W_{q,q_0}^{2,1}(\Omega^T)}$ and 
$\|\u_t^n\|_{\W_{p,p_0}^{2,1}(\Omega^T)}$ proportional to $\tau^a$ with 
$a>0$, we infer from (5.11) that
$$\eqal{
&\|\theta^{n+1}\|_{W_{q,q_0}^{2,1}(\Omega^\tau)}
\le\varphi(\tau^a\|\theta^n\|_{W_{q,q_0}^{2,1}(\Omega^\tau)}+c
\|\theta_0\|_{B_{q,q_0}^{2-2/q_0}(\Omega)},{\underline\theta})
\cdot[\tau^a\|\theta^n\|_{W_{q,q_0}^{2,1}(\Omega^\tau)}^2\cr
&\quad+\tau^a
\|\u_t^n\|_{\W_{p,p_0}^{2,1}(\Omega^\tau)}^2+
\|\theta_0\|_{L_q(\Omega)}^2+\|\u_1\|_{{\bm L}_p(\Omega)}^2+\|g\|_{L_{q,q_0}(\Omega^\tau)}\cr
&\quad+\|\theta_0\|_{B_{q,q_0}^{2-2/q_0}(\Omega)}]\cr}
\leqno(5.21)
$$
for $q,q_0,p,p_0\in(1,\infty)$, restricted by (5.13) and (5.20) ((5.18) is 
inactive).

\noindent
Let us denote
$$
X_n(\tau)=\|\u_t^n\|_{\W_{p,p_0}^{2,1}(\Omega^\tau)}+
\|\theta^n\|_{W_{q,q_0}^{2,1}(\Omega^\tau)}
\leqno(5.22)
$$
and
$$\eqal{
D&=\|\u_0\|_{\W_p^2(\Omega)}+\|\u_1\|_{\B_{p,p_0}^{2-2/p_0}(\Omega)}+
\|\theta_0\|_{B_{q,q_0}^{2-2/q_0}(\Omega)}+\|\b\|_{{\bm L}_{p,p_0}(\Omega^\tau)}\cr
&\quad+\|g\|_{L_{q,q_0}(\Omega^\tau)}.
\cr}
\leqno(5.23)
$$
Then inequalities (5.9) and (5.21) give
$$
\|\u_t^{n+1}\|_{\W_{p,p_0}^{2,1}(\Omega^\tau)}\le c[\tau^aX_n(\tau)+D],
\leqno(5.24)
$$
and
$$\eqal{
&\|\theta^{n+1}\|_{W_{q,q_0}^{2,1}(\Omega^\tau)}\le\varphi(\tau^aX_n(\tau),D)
[(\tau^aX_n(\tau))^2+D+D^2],\cr}
\leqno(5.25)
$$
where $p,p_0,q,q_0\in(1,\infty)$ satisfy (5.8), (5.13) and (5.20), $a>0$, 
$\tau$ is sufficiently small, and $c$ does not depend on $\tau$.

It follows from (5.24) and (5.25) that
$$\eqal{
&X_{n+1}(\tau)\le\varphi(\tau^aX_n(\tau),D)[\tau^aX_n(\tau)+(\tau^aX_n(\tau))^2+
D+D^2]\cr
&{\rm for}\ \ n\in{\N}\cup\{0\},\cr}
\leqno(5.26)
$$
where
$$
X_0(\tau)=\|\u_t^0\|_{\W_{p,p_0}^{2,1}(\Omega^\tau)}+
\|\theta^0\|_{W_{q,q_0}^{2,1}(\Omega^\tau)},
$$
and $p,p_0,q,q_0\in(1,\infty)$ satisfy (5.8), (5.13) and (5.20).

\noindent
Regarding the fact that $\varphi$ stands for a generic, positive, increasing 
function of its arguments, (5.26) can be expressed in the following 
simpler form
$$
X_{n+1}(\tau)\le\varphi(\tau^aX_n(\tau),D)\quad {\rm for}\ \ 
n\in{\N}\cup\{0\}.
\leqno(5.27)
$$
For $\tau$ sufficiently small there exists a positive constant $A$ such that
$$
X_0(\tau)\le A\quad {\rm and}\quad \varphi(\tau^aA,D)\le A.
\leqno(5.28)
$$
Then (5.27) implies that
$$
X_n(\tau)\le A\quad {\rm for\ any}\ \ n\in{\N}.
$$
This establishes $(5.4)_1$. Moreover, in view of (5.15) and the imbedding 
$W_{q,q_0}^{2,1}(\Omega^\tau)\subset C^\alpha(\Omega^\tau)$, the bounds 
$(5.4)_2$ hold. Thereby the proof is completed.
\end{proof}

To show the convergence of the sequence $\{\u^n,\theta^n\}$ we introduce the 
differences
$$
\U^n(t)=\u^n(t)-\u^{n-1}(t),\quad \vartheta^n(t)=\theta^n(t)-\theta^{n-1}(t),
\ \ n\in{\N}\cup\{0\},
\leqno(5.29)
$$
which are solutions to the following problems:
$$\eqal{
&\U_{tt}^{n+1}-\nabla\cdot(\A_1\epsi(\U_t^{n+1}))=\nabla\cdot
(\A_2\epsi(\U^n))-\nabla\cdot(\A_2\alphaa\vartheta^n)\quad &{\rm in}\ \ 
\Omega^T,\cr
&\U^{n+1}=\0\quad &{\rm on}\ \ S^T,\cr
&\U^{n+1}|_{t=0}=\0,\ \ \U_t^{n+1}|_{t=0}=\0\quad &{\rm in}\ \ \Omega,\cr}
\leqno(5.30)
$$
and
$$\eqal{
&c_v\theta^n\vartheta_t^{n+1}-k\Delta\vartheta^{n+1}=-c_v\vartheta^n\theta_t^n-\theta^n(\A_2\alphaa)\cdot
\epsi(\U_t^n)\quad\cr
&\quad-\vartheta^n(\A_2\alphaa)\cdot\epsi(\u_t^{n-1})
+\A_1\epsi(\U_t^n)\cdot\epsi(\u_t^n)+\A_1\epsi(\u_t^{n-1})\cdot
\epsi(\U_t^n)\quad &{\rm in}\ \ \Omega^T,\cr
&\bar\n\cdot\nabla\vartheta^n=0\quad &{\rm on}\ \ S^T,\cr
&\vartheta^{n+1}|_{t=0}=0\quad &{\rm in}\ \ \Omega,\cr}
\leqno(5.31)
$$
where
$$\eqal{
&\epsi(\u^{-1})=\0,\quad \epsi(\u_t^{-1})=\0,\quad \theta^{-1}=0,\cr
&\epsi(\U^0)=\epsi(\u_0),\quad \epsi(\U_t^0)=\epsi(\u_1),\quad 
\vartheta^0=\theta_0.\cr}
$$

Let us introduce the quantity
$$
Y^n(\tau)=\|\U_t^n\|_{\W_{\bar p,\bar p_0}^{2,1}(\Omega^\tau)}+
\|\vartheta^n\|_{W_{\bar q,\bar q_0}^{2,1}(\Omega^\tau)},
\leqno(5.32)
$$
where $\bar p,\bar p_0,\bar q,\bar q_0\in(1,\infty)$.

\begin{lemma}\label{5.2.} (Convergence of the approximation) 
Let the assumptions of Lemma 5.1 hold. Moreover, let the numbers 
$\bar p,\bar p_0,\bar q,\bar q_0\in(1,\infty)$ satisfy the following 
conditions:
$$\eqal{
&q=2\bar q,\quad q_0=2\bar q_0,\cr
&{3\over\bar q}+{2\over\bar q_0}-{3\over\bar p}-{2\over\bar p_0}<1,\cr
&{3\over\bar p}+{2\over\bar p_0}-{3\over2\bar q}-{2\over2\bar q_0}<1,\cr
&{3\over p}+{2\over p_0}-{3\over2\bar q}-{2\over2\bar q_0}<1.\cr}
$$
Then there exists a positive constant $d$ which depends on $A$, $D$, 
${\underline\theta}$, $p,p_0,q,\break
q_0,\bar p,\bar p_0$ such that
$$
Y^{n+1}(\tau)\le d\tau^aY^n(\tau),
\leqno(5.33)
$$
where $n\in{\N}\cup\{0\}$, $a>0$, and
$$
Y^0=\|\u^0\|_{\W_{\bar p,\bar p_0}^{2,1}(\Omega^\tau)}+
\|\theta^0\|_{W_{\bar q,\bar q_0}^{2,1}(\Omega^\tau)}.
$$
\end{lemma}

\begin{proof} 
By Lemma 3.4 the solutions to problem (5.30) satisfy the estimate
$$
\|\U_t^{n+1}\|_{\W_{\bar p,\bar p_0}^{2,1}(\Omega^\tau)}\le c_1
(\|\nabla^2\U^n\|_{{\bm L}_{\bar p,\bar p_0}(\Omega^\tau)}+
\|\nabla\vartheta^n\|_{{\bm L}_{\bar p,\bar p_0}(\Omega^\tau)}),
\leqno(5.34)
$$
where $c_1$ does not depend on $\tau$.

\noindent
The first term on the right-hand side of (5.34) is estimated by
$$
\|\nabla^2\U^n\|_{{\bm L}_{\bar p,\bar p_0}(\Omega^\tau)}=\bigg\|\intop_0^t
\nabla^2\U_{t'}^ndt'\bigg\|_{{\bm L}_{\bar p,\bar p_0}(\Omega^\tau)}\le c\tau
\|\U_t^n\|_{\W_{\bar p,\bar p_0}^{2,1}(\Omega^\tau)},
$$
and the second one by
$$\eqal{
&\|\nabla\vartheta^n\|_{{\bm L}_{\bar p,\bar p_0}(\Omega^\tau)}\le\delta
\|\vartheta^n\|_{W_{\bar q,\bar q_0}^{2,1}(\Omega^\tau)}+c(1/\delta)
\|\vartheta^n\|_{L_{\bar q,\bar q_0}(\Omega^\tau)}\cr
&\le(\delta+c(1/\delta)\tau)
\|\vartheta^n\|_{W_{\bar q,\bar q_0}^{2,1}(\Omega^\tau)},\cr}
$$
where $c(1/\delta)=\delta^{-a}$, $a>0$, and the interpolation inequality 
in the first line holds under the restriction
$$
{3\over\bar q}+{2\over\bar q_0}-{3\over\bar p}-{2\over\bar p_0}<1,\quad 
\bar p,\bar p_0,\bar q,\bar q_0\in(1,\infty).
\leqno(5.35)
$$
Hence,
$$\eqal{
&\|\U_t^{n+1}\|_{\W_{\bar p,\bar p_0}^{2,1}(\Omega^\tau)}\le c\tau^a
(\|\U_t^n\|_{\W_{\bar p,\bar p_0}^{2,1}(\Omega^\tau)}
+\|\vartheta^n\|_{W_{\bar q,\bar q_0}^{2,1}(\Omega^\tau)})\cr
&\le c\tau^aY_n(\tau).\cr}
\leqno(5.36)
$$
For the solutions to problem (5.31) we have
$$\eqal{
&\|\vartheta^{n+1}\|_{W_{\bar q,\bar q_0}^{2,1}(\Omega^\tau)}\le\varphi
\bigg(\sup_{\Omega^\tau}\theta^n,\sup_{\Omega^\tau}{1\over\theta^n}\bigg)
\Big[\|\vartheta^n\theta_t^n\|_{L_{\bar q,\bar q_0}(\Omega^\tau)}\cr
&\quad+\|\theta^n|\nabla\U_t^n|\,\|_{L_{\bar q,\bar q_0}(\Omega^\tau)}
\|\vartheta^n|\nabla\u_t^{n-1}|\,\|_{L_{\bar q,\bar q_0}(\Omega^\tau)}\cr
&\quad+\|\,|\nabla\U_t^n|\,|\nabla\u_t^n|\,\|_{L_{\bar q,\bar q_0}(\Omega^\tau)}+
\|\,|\nabla\U_t^n|\,|\nabla\u_t^{n-1}|\,
\|_{L_{\bar q,\bar q_0}(\Omega^\tau)}\Big].\cr}
\leqno(5.37)
$$
By virtue of Lemma 5.1 the arguments of the function $\varphi$ are estimated 
by $A$ and ${\underline\theta}/2$.
The successive terms under the square bracket on the right-hand side of 
(5.37) are estimated as follows.

The first term by
$$
\|\vartheta^n\|_{L_{2\bar q,2\bar q_0}(\Omega^\tau)}
\|\theta_t^n\|_{L_{2\bar q,2\bar q_0}(\Omega^\tau)}\le\varphi(A)\tau^a
\|\vartheta^n\|_{W_{\bar q,\bar q_0}^{2,1}(\Omega^\tau)}.
$$
This estimate follows on account of $(5.4)_1$ and similar arguments as in 
(5.7) under the restrictions
$$\eqal{
&q=2\bar q,\quad q_0=2\bar q_0,\quad {\rm and}\cr
&{3\over\bar q}+{2\over\bar q_0}-{3\over2\bar q}-{2\over2\bar q_0}<2,\quad 
{\rm so}\quad {3\over2\bar q}+{2\over2\bar q_0}<2.\cr}
\leqno(5.38)
$$
The second term is estimated by
$$\eqal{
&\|\theta^n\|_{L_{2\bar q,2\bar q_0}(\Omega^\tau)}
\|\nabla\U_t^n\|_{{\bm L}_{2\bar q,2\bar q_0}(\Omega^\tau)}\le
\|\theta^n\|_{L_{q,q_0}(\Omega^\tau)}\cdot
\tau^a\|\U_t^n\|_{\W_{\bar p,\bar p_0}^{2,1}(\Omega^\tau)}\cr
&\le\varphi(A)\tau^a\|\U_t^n\|_{\W_{\bar p,\bar p_0}^{2,1}(\Omega^\tau)},\cr}
$$
where we assumed that
$$
{3\over\bar p}+{2\over\bar p_0}-{3\over2\bar q}-{2\over2\bar q_0}<1.
\leqno(5.39)
$$
The third term is estimated by
$$
\|\vartheta^n\|_{L_{2\bar q,2\bar q_0}(\Omega^\tau)}
\|\nabla\u_t^{n-1}\|_{{\bm L}_{2\bar q,2\bar q_0}(\Omega^\tau)}\le
\varphi(A)\tau^a\|\vartheta^n\|_{W_{\bar q,\bar q_0}^{2,1}(\Omega^\tau)},
$$
where we assumed (5.38) and applied (5.4) under the condition
$$
{3\over p}+{2\over p_0}-{3\over2\bar q}-{2\over2\bar q_0}<1.
\leqno(5.40)
$$
In view of (5.39) and (5.40) the fourth term is estimated by
$$
\|\nabla\U_t^n\|_{{\bm L}_{2\bar q,2\bar q_0}(\Omega^\tau)}
\|\nabla\u_t^n\|_{{\bm L}_{2\bar q,2\bar q_0}(\Omega^\tau)}\le\varphi(A)\tau^a
\|\U_t^n\|_{\W_{\bar p,\bar p_0}^{2,1}(\Omega^\tau)}.
$$
A similar estimate holds for the last fifth term.

\noindent
Using the above estimates in (5.37) yields
$$
\|\vartheta^{n+1}\|_{W_{\bar q,\bar q_0}^{2,1}(\Omega^\tau)}\le\varphi
(A,D,{\underline\theta})\tau^aY^n(\tau).
\leqno(5.41)
$$
>From (5.36) and (5.41) we conclude (5.33).
The proof is complete.
\end{proof}

Lemmas 5.1 and 5.2 imply the existence of a local solution.

\begin{lemma}\label{5.3.} (Local existence) \\
Assume that $\u_0\in\W_p^2(\Omega)$, 
$\u_1\in\B_{p,p_0}^{2-2/p_0}(\Omega)$, 
$\theta_0\in B_{q,q_0}^{2-2/q_0}(\Omega)$,\\ 
$\b\in{\bm L}_{p,p_0}(\Omega^\tau)$, 
$g\in L_{q,q_0}(\Omega^\tau)$, $\theta_0\ge{\underline\theta}>\theta_*>0$, 
where $\theta_*$ is given by (4.2), and $p,p_0,q,q_0\in(1,\infty)$ satisfy
$$\eqal{
&{3\over q}+{2\over q_0}-{3\over p}-{2\over p_0}<1,\quad 
&{3\over p}+{2\over p_0}<2,\cr
&{3\over p}+{2\over p_0}-{3\over2q}-{2\over2q_0}<1,\quad 
&{3\over q}+{2\over q_0}<2.\cr}
\leqno(5.42)
$$
Then for $\tau$ sufficiently small there exists a solution to problem 
(1.1)--(1.4) such that
$$
\u_t\in\W_{p,p_0}^{2,1}(\Omega^\tau),\quad
\theta\in W_{q,q_0}^{2,1}(\Omega^\tau),\quad
\u\in C([0,\tau];\W_p^2(\Omega)),
$$
and
$$\eqal{
&\|\u_t\|_{\W_{p,p_0}^{2,1}(\Omega^\tau)}+
\|\theta\|_{W_{q,q_0}^{2,1}(\Omega^\tau)}+
\|\u\|_{C([0,\tau];\W_p^2(\Omega))}\cr
&\le\varphi(D,\theta_*,p,p_0,q,q_0),\cr}
\leqno(5.43)
$$
where
$$\eqal{
D=&D(p,p_0,q,q_0,\tau)\equiv\|\u_0\|_{\W_p^2(\Omega)}+
\|\u_1\|_{\B_{p,p_0}^{2-2/p_0}(\Omega)}\cr
&\quad+\|\theta_0\|_{B_{q,q_0}^{2-2/q_0}(\Omega)}+
\|\b\|_{{\bm L}_{p,p_0}(\Omega^\tau)}+\|g\|_{L_{q,q_0}(\Omega^\tau)}.\cr}
\leqno(5.44)
$$
\end{lemma}

\begin{proof} 
To satisfy the assumptions of Lemmas 5.1 and 5.2 we choose $q=2\bar q$, 
$q_0=2\bar q_0$ and 
${3\over\bar p}+{2\over\bar p_0}\ge{3\over p}+{2\over p_0}$.
Moreover, the assumptions of Lemma 5.1 imply that ${3\over q}+{2\over q_0}<4$, 
${3\over p}+{2\over p_0}<2$, and Lemma 5.2 gives ${3\over q}+{2\over q_0}<2$, 
${3\over\bar p}+{2\over\bar p_0}<1$. Hence, we choose 
${3\over p}+{2\over p_0}<2$, 
${3\over\bar p}+{2\over\bar p_0}<1$, ${3\over q}+{2\over q_0}<2$.
Then the sequence $(\u_t^n,\theta^n)$, $n\in{\N}$, of solutions to the 
approximate problem (5.1) is for sufficiently small $\tau$ uniformly 
bounded in the space 
$\W_{p,p_0}^{2,1}(\Omega^\tau)\times W_{q,q_0}^{2,1}(\Omega^\tau)$ and 
convergent in 
$\W_{\bar p,\bar p_0}^{2,1}(\Omega^\tau)\times 
W_{\bar q,\bar q_0}^{2,1}(\Omega^\tau)$.

\noindent
Hence, by virtue of the well known result (see the proof of Lemma 2.1 
from [15, Ch. 2]) we conclude that the limit
$$
(\u_t=\lim_{n\to\infty}\u_t^n,\theta=\lim_{n\to\infty}\theta^n)\in
\W_{p,p_0}^{2,1}(\Omega^\tau)\times W_{q,q_0}^{2,1}(\Omega^\tau)
$$
and satisfies problem (1.1)--(1.4).
\end{proof}

\begin{remark}\label{5.4.} 
Under the assumptions
$$\eqal{
&\u\in C([0,T];\W_p^2(\Omega)),\quad 
&\u_t\in C([0,T];\B_{p,p_0}^{2-2/p_0}(\Omega)),\cr
&\theta\in C([0,T];B_{q,q_0}^{2-2/q_0}(\Omega)),\qquad
&\b\in{\bm L}_{p,p_0}(\Omega\times(t,t+\tau)),\cr
&g\in L_{q,q_0}(\Omega\times(t,t+\tau)),\quad
&\theta(t)\ge\theta_*>0,\quad t\in[0,T],\cr}
$$
with $p,p_0,q,q_0$ satisfying (5.42), the assertion of Lemma 5.3 holds true 
for the interval $[t,t+\tau]$, where $\tau$ is a sufficiently small number 
depending on the data.
\end{remark}

\begin{lemma}\label{5.5.} (Global existence) 
Assume that
$$\eqal{
&\u_0\in\W_r^2(\Omega),\quad \u_1\in\B_{r,r_0}^{2-2/r_0}(\Omega),\quad 
\theta_0\in B_{\sigma,\sigma_0}^{2-2/\sigma_0}(\Omega),\cr
&g\in L_{\sigma,\sigma_0}(\Omega^T)\cap L_{\infty,r_0}(\Omega^T),\quad
\b\in{\bm L}_{r,r_0}(\Omega^T),\cr}
$$
with $r\ge\max\{12,p\}$, $r_0\ge\max\{12,p_0\}$, $\sigma\ge\max\{6,q\}$, 
$\sigma_0\ge\max\{6,q_0\}$, where $(p,p_0)$ and $(q,q_0)$ satisfy (5.42).
Then there exists a solution to problem (1.1)--(1.4) such that
$$
\u_t\in\W_{r,r_0}^{2,1}(\Omega^T),\quad
\theta\in W_{\sigma,\sigma_0}^{2,1}(\Omega^T),\quad
\u\in C([0,T];\W_r^2(\Omega)),
$$
satisfying the estimate
$$\eqal{
&\|\u_t\|_{\W_{r,r_0}^{2,1}(\Omega^T)}+\|\u\|_{C([0,T];\W_r^2(\Omega))}+
\|\theta\|_{W_{\sigma,\sigma_0}^{2,1}(\Omega^T)}\cr
&\le\varphi(D(r,r_0,\sigma,\sigma_0,T),\theta_*,r,r_0,\sigma,\sigma_0,T),\cr}
\leqno(5.45)
$$
where $D(r,r_0,\sigma,\sigma_0,T)$ is defined by (5.44) and $\varphi$ is 
a generic function.
\end{lemma}

\begin{proof} 
To extend the local solution from Lemma 5.3 step by step we choose 
parameters $p,p_0,q,q_0$ in such a way to satisfy simultaneously the 
restriction (4.115) required by a priori estimate (4.116) and the assumptions 
of Lemma 5.3. Then, in accord with (4.115) we assume that
$$\eqal{
&\|\u_0\|_{\W_{12/5}^2(\Omega)\cap\W_p^2(\Omega)}\le c
\|\u_0\|_{\W_r^2(\Omega)}\cr
&{\rm for}\ \ r\ge\max\bigg\{{12\over5},p\bigg\},\cr
&\|\u_1\|_{\B_{p,p_0}^{2-2/p_0}(\Omega)\cap\B_{12,12}^{11/6}(\Omega)}\le c
\|\u_1\|_{\B_{r,r_0}^{2-2/r_0}(\Omega)}\cr
&{\rm for}\ \ r\ge\max\{12,p\},\ \ r_0\ge\max\{12,p_0\},\cr}
\leqno(5.46)
$$
$$\eqal{
&\|\theta_0\|_{B_{q,q_0}^{2-2/q_0}(\Omega)\cap B_{6,6}^{5/3}(\Omega)}\le 
c\|\theta_0\|_{B_{\sigma,\sigma_0}^{2-2/\sigma_0}(\Omega)}\cr
&{\rm for}\ \ \sigma\ge\max\{6,q\},\quad \sigma_0\ge\max\{6,q_0\}.\cr}
$$
For $r$ and $\sigma$ specified in (5.46), the conditions (5.42) are satisfied, 
because ${3\over r}+{2\over r_0}\le{5\over12}$ and 
${3\over\sigma}+{2\over\sigma_0}\le{5\over6}$.

Let $t\in[0,T]$. By (4.116) in Corollary 4.25 and by the direct trace theorem 
we have
$$\eqal{
&\|\u_t(t)\|_{\B_{r,r_0}^{2-2/r_0}(\Omega)}+
\|\theta(t)\|_{B_{\sigma,\sigma_0}^{2-2/q_0}(\Omega)}\cr
&\le\varphi(D,\theta_*,r,r_0,\sigma,\sigma_0,T).\cr}
\leqno(5.47)
$$
Moreover,
$$\eqal{
\|\u(t)\|_{\W_r^2(\Omega)}&\le\intop_0^t\|\u_{t'}(t')\|_{\W_r^2(\Omega)}dt'+
\|\u_0\|_{\W_r^2(\Omega)}\cr
&\le T^{1-1/r_0}
\bigg(\intop_0^T\|\u_{t'}(t')\|_{\W_r^2(\Omega)}dt'\bigg)^{1/r_0}+
\|\u_0\|_{\W_r^2(\Omega)}\cr
&\le T^{1-1/r_0}\|\u_t\|_{\W_{r,r_0}^{2,1}(\Omega^T)}+
\|\u_0\|_{\W_r^2(\Omega)}\cr
&\le\varphi(D_0,\theta_*,r,r_0,\sigma,\sigma_0,T)\cr}
\leqno(5.48)
$$
for any $t\in[0,T]$.
Let $N\in{\N}$ be given so large that ${T\over N+1}\le\tau$, where $\tau$ 
is determined by Remark 5.4 and Lemma 5.3.
Using (5.47), (5.48) for $t=k\tau$ and taking into account that
$$\eqal{
&\sup_{0\le k\le N}\|\b\|_{{\bm L}_{r,r_0}(\Omega\times(k\tau,(k+1)\tau))}\le
\|\b\|_{{\bm L}_{r,r_0}(\Omega^T)},\cr
&\sup_{0\le k\le N}\|g\|_{L_{\sigma,\sigma_0}(\Omega\times(k\tau,(k+1)\tau))
\cap L_{\infty,r_0}(\Omega\times(k\tau,(k+1)\tau))}\cr
&\ \le\|g\|_{L_{\sigma,\sigma_0}(\Omega^t)\cap L_{\infty,r_0}(\Omega^T)},\cr}
$$
it follows from Remark 5.4 that there exists a solution $(\u,\theta)$ on 
the interval $[k\tau,(k+1)\tau]$, $0\le k\le N$, where $\tau$ does not 
depend on $k$.

\noindent
This implies the existence of the solution of problem (1.1)--(1.4) on the 
whole interval $[0,T]$, which ends the proof.
\end{proof}

Choosing $r=r_0=12$, $\sigma=\sigma_0=6$ in Lemma 5.5 and recalling 
Lemma 4.1 we complete the proof of Theorem A.

\section{Uniqueness (Proof of Theorem B)}

Assume that we have two solutions $(\u_1,\theta_1)$, $(\u_2,\theta_2)$ to 
problem (1.1)\---(1.4). Introducing the differences
$$
\U=\u_1-\u_2,\quad \vartheta=\theta_1-\theta_2,\quad \E=\epsi_1-\epsi_2,
\leqno(6.1)
$$
we see that they satisfy the problem
$$
\U_{tt}-\nabla\cdot[\A_1\E_t+\A_2(\E-\vartheta\alphaa)]=\0,
\leqno(6.2)
$$
$$\eqal{
&c_v(\theta_2\vartheta_t+\vartheta\theta_{1t})-k\Delta\vartheta=
-[\vartheta(\A_2\alphaa)\cdot\epsi_{1t}+\theta_2(\A_2\alphaa)\cdot\E_t]\cr
&\quad+\A_1\E_t\cdot\epsi_{1t}+\A_1\epsi_{2t}\cdot\E_t\quad {\rm in}\ \ 
\Omega^T,\cr}
\leqno(6.3)
$$
$$
\U|_{S^T}=\0,\quad \n\cdot\nabla\vartheta|_{S^T}=0\quad {\rm on}\ \ S^T,
\leqno(6.4)
$$
$$
\U|_{t=0}=\0,\quad \U_t|_{t=0}=\0,\quad \vartheta|_{t=0}=0\quad {\rm in}\ \ 
\Omega.
\leqno(6.5)
$$
Multiplying (6.2) by $\U_t$, integrating over $\Omega$ and integrating by 
parts we obtain
$$
{1\over2}{d\over dt}\bigg(\intop_\Omega\U_t^2dx+\intop_\Omega(\A_2\E)\cdot
\E dx\bigg)+\intop_\Omega\A_1\E_t\cdot\E_t dx=\intop_\Omega
(\A_2\vartheta\alphaa)\cdot\E_tdx.
$$
In view of (1.10) this implies
$$
{1\over2}{d\over dt}\bigg(\intop_\Omega\U_t^2dx+\intop_\Omega
(\A_2\E)\cdot\E dx\bigg)+a_{1*}\intop_\Omega|\E_t|^2dx\le\intop_\Omega
(\A_2\vartheta\alphaa)\cdot\E_tdx
$$
which, after applying the H\"older and the Young inequalities to the right 
hand side, yields
$$
{1\over2}{d\over dt}\intop_\Omega(\U_t^2+(\A_2\E)\cdot\E)dx+{1\over2}a_{1*}
\intop_\Omega|\E_t|^2dx\le c_1\intop_\Omega\vartheta^2dx.
\leqno(6.6)
$$

Further, multiplying (6.3) by $\vartheta$, integrating over $\Omega$, 
integrating by parts and using the boundary conditions implies
$$\eqal{
&{1\over2}\intop_\Omega c_v\bigg(\theta_2{\partial\over\partial t}\vartheta^2+
\theta_{1t}\vartheta^2)dx+k\intop_\Omega|\nabla\vartheta|^2dx\cr
&\le c_2\intop_\Omega(\vartheta^2|\epsi_{1t}|+\theta_2|\E_t|\,|\vartheta|)dx.
\cr}
$$
Continuing, we have
$$\eqal{
&{1\over2}{d\over dt}\intop_\Omega c_v\theta_2\vartheta^2dx+k\intop_\Omega
|\nabla\vartheta|^2dx\le c_3\intop_\Omega(|\theta_{1t}|+|\theta_{2t}|)
\vartheta^2dx\cr
&\quad+c_2\intop_\Omega(\vartheta^2|\epsi_{1t}|+\theta_2|\E_t|\,|\vartheta|)dx.
\cr}
\leqno(6.7)
$$
>From (6.7) we get
$$\eqal{
&{1\over2}{d\over dt}\intop_\Omega c_v\theta_2\vartheta^2dx+k
\|\vartheta\|_{H^1(\Omega)}^2\le\varepsilon_1\|\vartheta\|_{L_6(\Omega)}^2\cr
&\quad+c(1/\varepsilon_1)
(\|\theta_{1t}\|_{L_3(\Omega)}^2+\|\theta_{2t}\|_{L_3(\Omega)}^2)
\|\vartheta\|_{L_2(\Omega)}^2\cr
&\quad+\varepsilon_2\|\vartheta\|_{L_6(\Omega)}^2+c(1/\varepsilon_2)
\|\epsi_{1t}\|_{{\bm L}_3(\Omega)}^2\|\vartheta\|_{L_2(\Omega)}^2\cr
&\quad+\varepsilon_3\|\E_t\|_{{\bm L}_2(\Omega)}^2+c(1/\varepsilon_3)
\|\theta_2\|_{L_\infty(\Omega)}^2\|\vartheta\|_{L_2(\Omega)}^2+k
\|\vartheta\|_{L_2(\Omega)}^2.\cr}
\leqno(6.8)
$$
Adding (6.6) and (6.8), assuming that $\varepsilon_1$, $\varepsilon_2$, 
$\varepsilon_3$ are sufficiently small we obtain
$$\eqal{
&{d\over dt}\intop_\Omega(\U_t^2+(\A_2\E)\cdot\E+c_v\theta_2\vartheta^2)dx+
{1\over2}a_{1*}\intop_\Omega|\E_t|^2dx\cr
&\quad+k\|\vartheta\|_{H^1(\Omega)}^2\le c(c_1+\|\theta_{1t}\|_{L_3(\Omega)}^2+
\|\theta_{2t}\|_{L_3(\Omega)}^2+\|\epsi_{1t}\|_{{\bm L}_3(\Omega)}^2\cr
&\quad+\|\theta_2\|_{L_\infty(\Omega)}^2+k)\|\vartheta\|_{L_2(\Omega)}^2.\cr}
\leqno(6.9)
$$
By virtue of (4.2) it holds
$$
\theta_2\ge\theta_*>0.
$$
Thus, introducing
$$\eqal{
&X(t)=\intop_\Omega(\U_t^2+(\A_2\E)\cdot\E+c_v\theta_2\vartheta^2)dx,\cr
&A(t)=c_1+k+\|\theta_{1t}\|_{L_3(\Omega)}^2+\|\theta_{2t}\|_{L_3(\Omega)}^2+
\|\epsi_{1t}\|_{{\bm L}_3(\Omega)}^2+\|\theta_2\|_{L_\infty(\Omega)}^2,\cr}
$$
we conclude from (6.9) the inequality
$$
{d\over dt}X\le AX\quad {\rm for}\ \ t\in(0,T).
\leqno(6.10)
$$
Hence,
$$
X(t)\le X(0)\exp\bigg(\intop_0^tA(t')dt'\bigg).
\leqno(6.11)
$$
Since $X(0)=0$ and, by the assumption (1.17), $\intop_0^tA(t')dt'<\infty$, 
it follows that $X(t)=0$ for $t\in(0,T)$. This proves the uniquenes of the 
solution.

\end{document}